\documentclass[a4paper,12pt]{amsart}

\usepackage{amsmath,amssymb,amsthm}
\usepackage{amsfonts,latexsym} 

\theoremstyle{plain}
\newtheorem{theorem}{Theorem}[section]
\newtheorem{corollary}{Corollary}[section]
\newtheorem{lemma}{Lemma}[section]
\newtheorem{example}{Example}[section]
\newtheorem{proposition}{Proposition}[section]
\newtheorem{remark}{Remark}[section]
\numberwithin{equation}{section}

\begin{document}

\title[variation formulas (II)]
{Variation formulas for principal functions (II)\\
 \small{
Applications to variation for   harmonic spans}}

\author{Sachiko HAMANO, Fumio MAITANI, Hiroshi YAMAGUCHI}
\address{ Department of Mathematics, 
Faculty of Human Development and Culture, 
Fukushima University, Fukushima 960-1296
 JAPAN;  \quad 
2-7-7 Hiyoshidai, Ohtsu, 520-0112 JAPAN; \quad  2-6-20-3 
Shiromachi, Hikone, Shiga, 522-0068 JAPAN}
\email{ hamano@educ.fukushima-u.ac.jp; \quad hadleigh\_bern@ybb.ne.jp; 
\quad 
h.yamaguchi@s2.dion.ne.jp}
\thanks{This paper is dedicated to Professor Mitsuru Nakai on the 
occasion of his 77th birthday.} 

\begin{abstract} \ 
 For a domain $D$ in $\mathbb{ C}_z$ 
with smooth boundary and for $a,b\in D, 
 a\ne b$, we have the circular (radial) slit mapping $P(z)(Q(z))$ on $D$ 
 such that  $P(z)- \frac{ 1}{z-a}\ (Q(z)- \frac{1}{z-a})$ is regular at $a$ and 
 $P(b)(Q(b))=0$, and we call $p(z)=\log |P(z)|\ (q(z)=\log|Q(z)|)$  the 
 $L_1$-($L_0$-)principal function; \ 
$\alpha =\log|P'(b)|$  $(\beta =\log|Q'(b)|)$
 the $L_1$-($L_0$-)constant, and \  $s=\alpha - \beta 
 $ the harmonic span, for $D$. 
S.\,Hamano in \cite{hamano-2} showed the variation formula of the second order  
 for the $L_1$-const. $\alpha (t)$  for the moving 
domain $D(t)$ in $\mathbb{ C}_z$  with $t \in B:=\{t\in 
 \mathbb{ C}: |t|<\rho\}$. 
 We show  the corresponding formula for the
$L_0$-const. $\beta (t)$ for $D(t)$, and combine these formulas to obtain, 
if the total space ${\mathcal  
 D}=\cup_{t\in B}(t, D(t)) $ is  
 pseudoconvex in $ B \times \mathbb{ C}_z$, then $s(t)$ is subharmonic 
 on $B$. Since the geometric meaning of $s(t)$ is showed, this fact gives one of the 
 relations between the conformal mappings on each fiber $D(t), t\in B$ 
and the pseudoconvexity of ${\mathcal D}$. 
As a simple application we obtain the subharmonicity of 
$\log \cosh d(t)$ on $B$, where  $d(t)$ is the Poincar\'e  distance between $a$ and 
 $b$. 
\end{abstract}

\maketitle

\vspace{-5mm}
\section{Introduction }
\noindent Let $R$ be a 
bordered  Riemann surface with  boundary $\partial R=
C_1+ \cdots + C_\nu$ in a larger Riemann surface 
$\widetilde R$, where $C_j, j=1,\ldots , 
\nu$ is a $C^\omega$ smooth contour in $\widetilde R$.
Fix two points $a,\ b$ 
with local coordinates 
  $|z|<\rho$ and $|z-\xi|<\rho$ where $a\,(b)$ 
corresponds to $0\,(\xi)$. Among 
all harmonic functions $u$ on $R \setminus \{0, \xi\}$ 
with logarithmic 
singularity $\log \frac{ 1}{|z|}$ at $0$ and $\log |z-\xi|$ 
at $\xi$ normalized  
$\lim_{ z\to 0}(u(z)-\log \frac{ 1}{|z|})=0$,  we uniquely have two special ones 
$p$ and 
$q$ with the following boundary conditions: for each $C_j$, 
$p$ satisfies $p(z)= {\rm const.\ } c_j$ on $C_j$ and $\int_{ C_j}\frac{\partial 
p(z)}{\partial n_z}ds_z=0$ (where $\frac{\partial }{\partial n_z}$ is the outer 
normal derivative  and $ds_z$ is the arc length element 
at $z$ of $C_j$),
while $q$ does $\frac{\partial q(z)}{\partial n_z}=0$ on $C_j$.
We call $p $ and $q$ the $L_1$- and the 
$L_0$-principal function for $(R,0,\xi)$, respectively.  The constant terms 
$\alpha :=\lim_{ z\to \xi}(p(z)-\log|z-\xi|)$ and $\beta :=\lim_{ 
z\to \xi}(q(z)-\log|z-\xi|)$ are called the $L_1$- and the $L_0$-constant for 
$(R, 0, \xi)$ (see \cite{ahlfors-sario} and  \cite{nakai-sario}).

Now let $B=\{t\in \mathbb{ C}: |t|<\rho\}$  and let 
${\mathcal  R}: t\in B \to R(t) \Subset 
\widetilde R$ be a variation of Riemann surface $R(t)$ 
with $t\in B$ such that   each $R(t), t\in B$ contains the origin $0$; 
$\partial R(t)=C_1(t)+ \cdots +C_\nu(t)$ is $C^\omega$ smooth in 
$\widetilde R$, and   
$\partial R(t)$ varies $C^\omega$  smoothly on $\widetilde R$ with 
$t\in B$. Let $\xi (t) \in R(t), \ t\in B$ vary holomorphically in $\widetilde R$ 
with $t\in B$.
Then each $R(t), t\in B$ admits the 
$L_1$-($L_0$-)principal function $p(t,z )(q(t,z))$ and  
 $L_1$-($L_0$-)constant $\alpha (t)(\beta (t))$ for $(R(t),0, \xi(t))$.
  S.\,Hamano \cite{hamano-2} established 
 the  variation formula
of  the second order for $\alpha (t)$ (see Lemma \ref{vari-form:u} in 
this paper), which  implied that,
if the total space ${\mathcal  R}=\cup_{t\in B}(t, R(t))$ is  a
pseudoconvex domain in $B \times \widetilde R$, then $\alpha (t)$ is 
subharmonic on $B$. 
Continuing  on  \cite{hamano-2}  
we  show the variation formula for $\beta(t)$ 
(Lemma \ref{vari-form:L_0}) in this paper, which continues on \cite{hamano-3}.
 To prove the formula for $\beta (t)$ we 
 add  a new idea to her proof for $\alpha (t)$.
In fact, the formula for $\alpha (t)$ does not concern 
 the genus of $R(t)$ but that for $\beta(t)$ does 
concern it. The formula for $\beta (t)$ 
implies that, if ${\mathcal  R}$ is 
pseudoconvex in $B \times 
\widetilde R$ and if  $R(t), 
t\in B$ is planar, then $\beta (t)$ is 
superharmonic on $B$. 
This  contrast  
between the subharmincity of $\alpha (t)$
 and  the superharmonicity of $\beta (t)$
 are unified with  the notion of 
the harmonic span $s(t):=\alpha(t)-\beta(t)$ 
 for $(R(t), 0, \xi(t))$ introduced by M.\,Nakai\ (see (\ref{nakai}) in 
 \S 3):
 {\it if  ${\mathcal  R}$ is pseudoconvex in $B \times \widetilde R$ and 
 $R(t), t\in B$ 
is planar, then   $s(t)$ is subharmonic on $B$.}
This implies (Corollary \ref{cor:mdp}): assume moreover that each $R(t), 
t\in $  is simply connected. Let 
$\xi_i:=\cup_{t\in B}(t, \xi_i(t)), i=1,2$ 
be two holomorphic sections of 
${\mathcal  R}$ over $B$  and 
 let $d(t)$ denote the Poincar\'e distance between $\xi_1(t)$ and 
 $\xi_2(t)$ on $R(t)$.
Then 
 $\delta(t):=\log\,\cosh { d(t)}$ is subharmonic on $B$. Further,  
$\delta(t)$ is harmonic on $B$ iff ${\mathcal R}$ is biholomorphic to the product $B \times R(0)$. 

We thank Professor M. Nakai for his helpful advice in our study of 
harmonic spans. We also thank Professor M. Brunella for his kind comment to Corollary 
\ref{cor:mdp}.

\section {Variation formulas for $L_0$-principal functions} 

\noindent Let $B=\{t\in \mathbb{ C}: |t|<\rho\}$ and let $\widetilde {\mathcal  R}$ be an
 unramified (Riemann) domain over $B \times \mathbb{C}_z$. We write 
$\widetilde {\mathcal  R}=\cup_{t\in B}(t, \widetilde R(t))$, where 
$\widetilde R(t)$ is the  fiber of $\widetilde {\mathcal R}$  over  $t\in B$,\,i.e., 
$\widetilde R(t)=\{z: (t,z)\in \widetilde {\mathcal  R}\}$. We assume 
$\widetilde R(t)\ne \emptyset$ for any $t\in B$, so that  
$\widetilde R(t)$ is  Riemann surfaces sheeted over 
$\mathbb{C}_z$ without ramification points.
Consider a subdomain ${\mathcal  R}$ in $ \widetilde {\mathcal  R}$ 
such that, putting ${\mathcal  R}=\cup _{t\in B} (t, R(t))$, where 
$R(t)$ is the  fiber of ${\mathcal  R} $ over $t\in B$, 
\begin{enumerate}
 \item [{\bf  1.}] $\widetilde R(t) \Supset R(t) \ne \emptyset
$, $t\in B$ and $R(t)$ is a connected Riemann surface 
of genus $g\ge 0$ such that  $\partial R(t)$ in $\widetilde R(t)$ 
consists of a finite number of $C^\omega$ smooth contours 
$C_j(t),\,j=1,\ldots ,\nu$;
 \item [{\bf  2.}]  the boundary 
$\partial {\mathcal  R}=\cup_{t\in B}(t, \partial R(t))$ of 
${\mathcal  
R}$ in $\widetilde {\mathcal R}  $ is $C^\omega$ smooth.
\end{enumerate}
Note that $g$ and $\nu$ are independent of $t\in B$.
We give the orientation of $C_j(t)$ such that  $\partial R(t)= C_1(t)+
\cdots +C_\nu(t)$.
\noindent We regard the two-dimensional unramified   domain ${\mathcal  
R}$ over $B \times \mathbb{C}_z$ as a $C^\omega$ smooth 
variation of Riemann surfaces $R(t)$ (sheeted over $\mathbb{ C}_z$ without ramification 
points and  with $C^\omega$ 
smooth boundary $\partial R(t)$) with complex parameter $t\in B$,
$$
{\mathcal  R}: t\in B \to R(t) \  \Subset \widetilde R(t).
$$
{We denote by $\Gamma (B, {\mathcal  R})$ the set of all 
holomorphic sections of ${\mathcal  R}$ over $B$.}
Assume that there exist $\Xi _0, \ \Xi_\xi \in \Gamma (B, 
{\mathcal  R})$ such that  $
\Xi _0:\, z=0$ and $\Xi_\xi : \, z= \xi (t)$ with 
$ \Xi_0\cap\Xi_ \xi =\emptyset$.
Let $t\in B$ be fixed. It is known (cf: \S\,3, Chap.\,III in \cite{ahlfors-sario})
that $R(t)$ carries 
the 
real-valued  
functions $p(t,z)$ and $q(t,z)$ such that  both functions 
are continuous on $\overline{ R(t)}$ and 
harmonic on $R(t) \setminus \{0, \xi(t)\}$ with poles 
$\log \frac{  
1}{|z|}$ at $z=0$ and $\log |z-\xi(t)|$ at $z=\xi(t)$ normalized 
$
\lim_{ z\to 0}(\,p(t,z)-\log \frac{   1}{|z|})
=\lim_{ z\to 0}
(\,
\,q(t,z)-\log \frac{ 1}{|z|}\,)=0
$
at $z=0$, and 
$p(t,z)$ and $q(t,z)$ satisfy the following boundary condition $(L_1)$
 and $(L_0)$, respectively: for $j=1,\ldots , \nu,$ 
\begin{align*}
 &(L_1) \quad \  p(t,z)= \mbox{const.}\,c_j(t) \ \mbox{  on }  C_j(t) 
\ \ \ \mbox{ and}  \ \ \  
 \int_{ C_j(t) } \frac{\partial p(t,z)}{\partial n_z}ds_z=0\,; \\
&(L_0) \quad \ \frac{\partial q(t,z)}{\partial n_z}=0\  \mbox{  on }  C_j(t). 
\end{align*}
We call $p(t,z)$ and $q(t,z)$ 
the $L_1$- and the  $L_0$-{\it principal function} for $(R(t),0,\xi(t))$.
 We  find a neighborhood $U_0(t)
$ of  $z=0$ such that 
\begin{align}\label{eqn:beta-at-0}
 p(t,z)&=\log\,\mbox{$\frac{1}{|z|}$} + h_0(t,z)
 \quad \mbox{\ on\ }U_0(t); \nonumber \\
 q(t,z)&=\log\,\mbox{$\frac{1}{|z|}$} + \mathfrak{h}_0(t,z) \quad \mbox{\ on\ }U_0(t), 
\end{align}
where $h_0(t,z)$,  $\mathfrak{ h}_0(t,z)$ are  harmonic for $z$ on $U_0(t)
$ and  
$$
h_0(t,0), \,\mathfrak{ h}_0(t,0) \equiv 0  \quad \mbox{  on $B$}.
$$
We also find a neighborhood $U_\xi(t)$ of $z=\xi(t)$ such that 
\begin{align}\label{eqn:beta}
 p(t,z)&=\log {|z-\xi(t)|} + \alpha(t)+ h_\xi(t,z) 
\quad \mbox{\ on\ }U_\xi(t); \nonumber \\ 
 q(t,z)&=\log {|z-\xi(t)|} + \beta(t)+ 
\mathfrak{h}_\xi(t,z) \quad \mbox{\ on\ }U_\xi(t), 
\end{align}
where $\alpha(t), \,\beta (t)$ are  the constant terms,  and 
$h_\xi(t,z), \,\mathfrak{h}_\xi(t,z)$ are harmonic for $z$ on 
$U_\xi(t)$ and   
\begin{align}
 \label{eqn:def-hxi}
h_\xi (t,\xi (t)), \ \mathfrak{ h}_\xi (t, \xi(t)) \ \equiv  0 
\quad \mbox{ on } B.
\end{align}
We call $\alpha(t)$ and $\beta (t)$ 
the $L_1$- and the $L_0$-{\it constant}
for $(R(t),0,\xi(t))$.

\vspace{1mm}  The following variation formula of the second order 
for  $\alpha(t)$ is showed: 
\begin{lemma}[Lemma 1.3 in \cite{hamano-2}] \ 
\label{vari-form:u} 
\begin{align*}
    \frac{\partial \alpha (t)}{\partial t}&=
\frac{1}{\pi} 
\int_{\partial R(t)} \!\!\! {k_1(t,z)} 
\bigl|\frac{\partial p(t,z)}{\partial z}
\bigr|^2 ds_z +2\  \frac{\partial h_\xi }{\partial z}\bigl|_{(t, \xi (t))}\cdot 
\  \xi '(t)\ ;\\[3mm]
 \frac{\partial^2 \alpha(t)}{\partial t \partial \bar{t}}&=
\frac{1}{\pi} 
\int_{\partial R(t)} \!\!\! {k_2(t,z)} 
\bigl|\frac{\partial p(t,z)}{\partial z}
\bigr|^2 ds_z +\frac{4}{\pi} \iint_{R(t)} 
\bigl|\frac{\partial^2 p(t,z)}{\partial \bar{t} \partial z}\bigr|^2
 dxdy.
\end{align*}
Here 
\begin{align*}
  k_1(t,z)&= \frac{\partial \varphi }{\partial t}/ \frac{\partial 
 \varphi }{\partial z }\ ; \\[2mm]
k_2(t,z)&=
\bigl(
\frac{\partial^2 \varphi}{\partial t \partial \bar{t}}
\bigl|
\frac{\partial \varphi}{\partial z}
\bigr|^2 -2{\rm Re} \bigl\{
\frac{\partial^2 \varphi}{\partial \bar{t} \partial z}
\frac{\partial \varphi}{\partial t}
\frac{\partial \varphi}{\partial \bar{z}}
\bigr\}
+\bigr|
\frac{\partial \varphi}{\partial t}
\bigr|^2 
\frac{\partial^2 \varphi}{\partial z \partial \bar{z}}
\bigr)
/\bigl|\frac{\partial \varphi}{\partial z}\bigr|^3
\end{align*}
on $\partial\mathcal R$,  where $\varphi(t,z) $ is a $C^2$ defining function of 
 $\partial {\mathcal  R}$.
\end{lemma}
Note that $k_i(t,z), i=1,2$ on $\partial {\mathcal  R}$ 
does not depend on the choice 
of defining functions $ \varphi(t,z)$ of $\partial{\mathcal R}$, where 
$k_1(t,z)$ is due to Hadamard and $k_2(t,z)$  
is called the {\it Levi curvature} for $\partial {\mathcal  R}$
((1.3) in \cite{l-y} and  (7) in \cite{M-Y}). 
The first formula in the  lemma is  proved by the similar 
method to that in Lemma \ref{vari-form:L_0}  below.

\vspace{2mm} 
 We shall 
give the variation formulas 
 for  $\beta(t)$.  In case 
when $R(t)$ is of  positive genus $g\ge1$ we need the following consideration,
 which was not necessary for the variation formulas for
$\alpha (t)$. 
We draw  as usual  $A, B$ cycles  $\{A_k(t), B_k(t)\}_{1\le k\le g}$ on 
$R(t)$ which vary continuously in ${\mathcal R}$ with $t\in B$ without passing through $0, \, \xi(t)$:
\begin{eqnarray}
  \label{akbk} 
\begin{array}{llll}
 \empty{}\quad &   A_k(t)\cap B_l(t)
               &=\emptyset \mbox{ for $k\ne l$}; 
                A_k \times B_k =1 \mbox{ for $k=1,\ldots , g$;}\\[2mm]
\empty{}\quad  & A_k(t)\cap A_l(t)
               &=B_k (t)\cap B_l(t)=\emptyset \mbox{ for  $k\ne l .$}
\end{array}
\end{eqnarray}
Here $A_k(t)\times B_k(t)=1$ 
means that  $A_k(t)$ once crosses $B_k(t) $ from the left-side to the right-side of 
the direction $B_k(t)$.
On $R(t), t\in B$ we put $*dq(t,z)
= -\frac{\partial q(t,z)}{\partial y }dx 
+ \frac{\partial q(t,z)}
{\partial x}dy$,  
the conjugate differential of 
$dq(t,z)$.
\begin{lemma}
\label{vari-form:L_0} 
\begin{align*}
    \frac{\partial \beta (t)}{\partial t}&=
-\frac{1}{\pi} 
\int_{\partial R(t)} \!\!\! {k_1(t,z)} 
\bigl|\frac{\partial q(t,z)}{\partial z}
\bigr|^2 ds_z +2\  \frac{\partial \mathfrak{ h}_\xi }{\partial z}\bigl|_{(t, \xi (t))}\cdot 
\  \xi '(t)\ ;\\[3mm]
 \frac{\partial^2 \beta(t) }{\partial t \partial\overline{t}} 
&= - \frac{1}{\pi} \int_{\partial R(t)}k_2 (t,z) 
\bigl|
\frac{\partial q(t,z)}{\partial z} 
\bigr|^2 ds_z 
- \frac{4}{\pi} \iint_ {R(t)} 
\bigl|
\frac{\partial^2 q(t,z)}{\partial\overline{t} \partial{z}}
\bigr|^2dxdy\\
& \quad -\, \frac{ 2}{\pi}\, \Im\, { \,\sum_ {k=1}^g
 \bigl(\frac{\partial }{\partial t} 
\int_{ A_k(t)}\!\! * dq(t,z)\bigr)\cdot\bigl(
\frac{\partial }{\partial \overline{ t}}
 \int_{ B_k(t)}\!\! * dq(t,z) \bigr) }.
\end{align*}
\end{lemma}
\noindent {\it  Proof.} \   It suffices to prove 
the lemma at $t=0$. If necessary, take a smaller disk $B$ of center  $0$. 
Since both $\partial {\mathcal  R}$ in $\widetilde {\mathcal  R}$ and 
$\partial R(t)$ in $\widetilde R(t)$ are  
$C^\omega$ smooth, we find  a  neighborhood $V=\cup _{j=1}^\nu V_j $ (disjoint union) of $\partial 
R(0)=\cup_{j=1}^\nu C_j(0)$  such that $(B\times V)\cap (\Xi_0 \cup  \Xi_\xi)=\emptyset$;  $V_j$ is a thin 
tubular neighborhood of $C_j(0)$ with  $V_j \supset C_j(t)$ for any $t\in B$,   
and $q(t,z)$ is harmonic on 
$[R(0)\cup V] \setminus \{0, \xi(t)\}$. We  write 
$\widehat R(0):\,={ R(0)}\cup V$, so that 
$q(t,z)$ is defined in the product $B  
\times \widehat R(0)$.

We divide the proof into two steps.

 {\it  $1^{st}$ step. \  Lemma 
\ref{vari-form:L_0} is true in the special case when 
$\Xi_\xi$ is a constant section, say, for example,   $\Xi_1: z= 1$ 
on $B$.}

In fact, formula (\ref{eqn:beta}) becomes
\begin{align}\label{eqn:beta2}
 q(t,z)=\log {|z-1|} + 
\beta(t)+ \mathfrak{ h}_1(t,z) \qquad \mbox{\ on\ }U_1(t),
\end{align}
where $\mathfrak{ h}_1(t,1)\equiv 0$ on $B$. 
 For $t\in B$ we put $u(t,z):= q(t,z)-q(0,z)$ on $\widehat{ R}(0) 
 \setminus \{0,1\}$. 
By putting $u(t,0)=0$ and $u(t,1)= \beta (t)-\beta (0)$, 
$u(t,z)$ is harmonic on  $\widehat{ R}(0)$. 
Let $\varepsilon : 0<\varepsilon \ll 1$, 
$\gamma _\varepsilon (0)=\{|z|< 
\varepsilon \}$ and $\gamma _\varepsilon (1)=\{|z-1|<\varepsilon \}$.
Then  Green's formula implies 
$$\int_{ \partial R(0)-\partial \gamma _\varepsilon (0)- \partial \gamma _\varepsilon (1)} 
 u(t,z) \frac{\partial q(0,z)}{\partial n_z}ds_z
- q(0,z) \frac{\partial u(t,z)}{\partial n_z}ds_z=0.
$$
Letting $\varepsilon \to 0$, we have from $\frac{\partial q(0,z)}{
\partial n_z }=0 $ on $C_j(0), j=1,\ldots , \nu$,
\begin{align}\label{eqn:kihon}
  \beta(t) -\beta(0)&= \frac{ -1}{2\pi} \sum_ {j=1}^\nu \int_{ C_j(0)} 
 q(0,z) \frac{\partial q(t,z)}{\partial n_z}ds_z
=: \frac{ -1}{2\pi} \sum_ {j=1}^\nu I_j(t). 
\end{align}
We take a point $z_j^0(t)$ on each $C_j(t), t\in B$ such that $z_j^0(t)$ 
continuously moves in $\partial {\mathcal  R} $ with $t\in B$, 
and  choose a harmonic conjugate function  $q^*_j(t,z)$ of $q(t,z)$ 
in $V_j$  such that  $q^*_j(t, z_j^0(t))=0$.
Since $\frac{ \partial q(t,z)}{\partial n_z }=0 $ on $C_j(t)$, 
 $q^*_j(t,z)$ is single-valued  in $V_j$ and 
\begin{align}
 \label{eqn:q-star-zero}
q^*_j(t,z) =0  \quad \mbox{ for $z\in C_j(t)$, $t\in B$}.
\end{align}
Since $dq^*_j(t,z)= \frac{\partial q(t,z)}{\partial n_z }ds_z $, 
 $dq(0,z)= - \frac{\partial q^*_j(0,z) }{\partial n_z 
}ds_z$ along $C_j(0)$,
we have 
\begin{alignat*}{3}
I_j(t)&= \int_{ C_j(0)} q(0,z)d q^*_j(t,z)\\
&=\int_{ C_j(0)} d[q(0,z) q_j^*(t,z)] - q^*_j(t,z)dq(0,z)\\
&=\int_{ C_j(0)}  q^*_j(t,z) \frac{\partial q^*_j(0,z) }{\partial n_z 
}ds_z.
\end{alignat*}
Differentiating both sides by $t$ and $\overline{t}$ at $t=0$,
 we have 
\begin{align}
 \label{eqn:Ij-2-d1}
\frac{\partial I_j}{\partial t }(0)
&=\int_{ C_j(0)} \frac{\partial q^*_j }{\partial t} (0,z) \frac{\partial 
 q^*_j(0,z)}{\partial n_z} ds_z;
\\ \label{eqn:Ij-2-d2}
\frac{\partial^2 I_j}{\partial t \partial \overline{ t} }(0)
&= \int_{ C_j(0)} \frac{\partial ^2q^*_j }{\partial t \partial \overline{ 
t}}(0,z) 
\frac{\partial q^*_j(0,z) }{\partial n_z }ds_z.
\end{align}

We  recall the following 
 \begin{proposition}[(1.2) in \cite{hamano-2}] \label{lem:hamano}
\  Let $u(t,z)$ be a $C^2$ function for $(t,z)$ in a neighborhood 
 ${\mathcal    V}_j=\cup_{t\in B}(t, V_j(t))$ of 
 ${\mathcal   C}_j=\cup _{t\in B} (t, C_j(t)) $ over  $B \times \mathbb{C}_z$ 
  such that  each $u(t,z), t\in B$ is harmonic for $z$ in $V_j(t)$ and 
 $u(t,z)= \mbox{a certain  const.} \ c_j(t)$ on $C_j(t)$.
Then 
\begin{align*}
& (i) \ \frac{\partial u}{\partial t } \frac{\partial 
 u}{\partial n_z}ds_z = 2\, 
 k_1(t,z) \bigl|\frac{\partial u }{\partial z}\bigr|^2ds_z \quad 
 \mbox{ along $C_j(t)$};\\[1mm]
 & (ii) \ \frac{\partial^2 u}{\partial t \partial \overline{ t}} 
\frac{\partial u}{\partial n_z}ds_z 
= 2\,k_2(t,z)|\frac{\partial u }{\partial z }|^2 ds_z +
 \frac{\partial ^2 c_j(t)}{\partial t \partial \overline{ t}} \frac{\partial u}{\partial 
n_z}ds_z\\
& \qquad \qquad \quad   
+4\, \Im \, \bigl\{ \frac{\partial u}{\partial t} \frac{\partial^2 u}
{\partial \overline{ t}\partial z}dz  \bigr\} 
 - 4 \Im\, \bigl\{ \frac{\partial c_j(t)}{\partial t} \frac{\partial 
^2 u }{\partial \overline{ t} \partial z } dz \bigr\} \   \mbox{ along} \  
C_j(t). 
  \end{align*}
\end{proposition}
We apply (i)  for $u(t,z)=q^*_j(t,z)$ with
 (\ref{eqn:q-star-zero}) to  
(\ref{eqn:Ij-2-d1}) and obtain
\begin{align*}
 \frac{\partial I_j}{\partial t }(0) &= 2 \int_{C_j(0)}  k_1(0,z) 
\bigl|\frac{\partial q_j^*(0,z)}{\partial z } \bigr|^2 ds_z.\\
\therefore \ \ 
 \frac{\partial \beta }{\partial t}(0)&= -\frac{ 1}{\pi} 
\int_{\partial 
R(0)} 
k_1(0,z) \bigl| \frac{\partial q(0,z)}{\partial z}\bigr|^2ds_z 
\quad \mbox{ by  (\ref{eqn:kihon}),}
\end{align*}
which proves the first formula in Lemma \ref{vari-form:L_0} in 
the 1st step.

To prove the second one, we apply (ii) for $u(t,z)=q^*_j(t,z)$ 
with (\ref{eqn:q-star-zero}) to (\ref{eqn:Ij-2-d2}) and obtain 
\hspace{-1cm}
\begin{align*}
 \frac{\partial^2I_j }{\partial t \partial \overline{ t}}(0)= 
2 \int_{ C_j(0)}\!\! k_2(0,z)\bigl|
\frac{\partial q^*_j(0,z)}{\partial z }\bigr|^2ds_z
+ 4\, \Im\,{
 \int_{ C_j(0)}\! \frac{\partial q^*_j}{\partial t}(0,z) \frac{\partial ^2 
 q^*_j
}{\partial \overline{ t}\partial z }(0,z)dz}.
\end{align*}
We put 

$
\qquad { {\bf  a}_k(t)= \int_{A_k(t)} *dq(t,z), \quad  \quad 
{\bf  b}_k(t)=\int_
{ B_k(t)} *dq(t,z).}
$

\vspace{1.5mm} 
\noindent  We fix a point $z^0 \,(\ne 0,1)$ such that  $B \times \{z^0\} \subset {\mathcal  
R}$. On each $R(t), t\in B$ we  choose a branch $q^*(t,z)$ of harmonic 
conjugate function of 
 $q(t,z)$ on $\widehat R(0)\setminus \{0, 1\}$ such that  
$ q^*(t, z^0)=0$. Since $\int_{ C_j(0)}*dq(t,z)=0$, we 
have 

\vspace{1.5mm} $
\qquad \quad q^*(t,z')= q^*(t,z'') \ \ \mbox{ mod } \{2\pi, \ {\bf  a}_k(t), 
\ {\bf  b}_k(t) \ 
(k=1,\ldots ,g)\}
$

\vspace{1.5mm} 
\noindent 
for any $z',z''$ over the 
same point $z\in \widehat R(0) \setminus \{0,1\}$. 
 We also  have  

\vspace{1.5mm} 
$
\qquad \qquad q^*_j(t,z)-q^*(t,z) =  c_j(t) \quad \mbox{ on $V_j$},
$

\vspace{1.5mm} 
\noindent where $c_j(t)$ is a certain constant for $z\in V_j$. 
 It follows that
\hspace{-1cm}
 \begin{alignat*}{3}
&\int_{ C_j(0)}\frac{\partial q_j^* }{\partial t }(0,z)
 \frac{\partial^2 q^*_j}{ \partial \overline{ t}\partial z}(0,z)dz\\
&\quad = \int_{ C_j(0)}\frac{\partial q^* }{\partial t }(0,z)
 \frac{\partial^2 q^*}{ \partial \overline{ t}\partial z}(0,z)dz
+ \frac{\partial c_j}{\partial t }(0) \int_{ C_j(0)}
\frac{\partial ^2q^* }{\partial \overline{ t}\partial z }(0,z) dz.
\end{alignat*}
If we put $f(t,z):=q^*(t,z)- iq(t,z)$ for $(t,z)\in B \times V_j$, 
then $f\in C^\omega(B \times V_j)$ and each $f(t,z), t\in B$ is
single-valued and holomorphic for $z$ in $V_j$, so that
\begin{alignat*}{3}
 \hspace{-3mm}
&\quad \int_{ C_j(0)} \frac{\partial ^2 q^*}{\partial \overline{ t}\partial 
z}(0,z) dz
= \frac{ 1}2 \, \bigl[\,\frac{\partial }{\partial \overline{ t }}\
 \bigl( \int_{ C_j(0)}
f'_z(t,z) dz \bigr) \
\bigr]_{t=0} =0.\\[3mm]
& \therefore \ 
 \frac{\partial^2 I_j }{\partial t \partial \overline{  t}}(0)= 
2\int_{ C_j(0)}\!\! k_2(0,z)\bigl|\frac{\partial q^*(0,z) }{\partial z } 
\bigr|^2 ds_z  + 4 \,\Im \,\bigl\{
\int_{ C_j(0)} \frac{\partial q^*}{\partial t } (0,z) 
\frac{\partial ^2
q^*}{\partial \overline{ t} \partial z}(0,z) dz\bigr\}.
\end{alignat*}
It follows from (\ref{eqn:kihon}) that
\begin{align}
\label{eqn:kihon-2} 
\frac{\partial^2 \beta}{\partial t \partial \overline{ t} }(0)&=
- \frac{ 1}\pi 
\int_{ \partial R(0)} k_2(0,z)\bigl|\frac{\partial q^*(0,z) }{\partial z } 
\bigr|^2 ds_z
-\  \frac{ 2}\pi\,  
\Im\,  \bigl\{
\int_{\partial R(0)} \frac{\partial q^*}{\partial t } (0,z) \frac{\partial ^2
q^*}{\partial \overline{ t} \partial z}(0,z) dz
\bigr\}.
\end{align}

We shall divide the proof into two cases.

\vspace{1mm} 
{\it Case when $R(t)$ is planar,\,i.e., $g=0$.} \ In this case, each $q^*(t,z), t\in B$ 
is determined up to additive constants mod $2\pi$.
It follows from 
(\ref{eqn:beta-at-0}) and (\ref{eqn:beta2})  that, for any fixed $t\in B$,  
$\frac{\partial q^*(t,z)}{\partial t}$ is a single-valued 
harmonic function on $\widehat R(0)$, and  
 $\frac{\partial ^2q^*(t,z)}{\partial 
\overline{ t} \partial z}$ is a single-valued  holomorphic function 
 on $\widehat R(0)$. 
We have by Green's formula
\begin{align*}
\int_{\partial R(0) }  \ \frac{\partial q^*}{\partial t } (0,z) 
\frac{\partial ^2
q^*}{\partial \overline{ t} \partial z} (0,z) dz
&=2i \iint_{ R(0)} \bigl| \frac{\partial^2 q^* }{\partial \overline{ t} 
\partial z} (0,z)\bigr|^2dx\,dy.
\end{align*}
$$\therefore \ \ 
\frac{\partial^2 \beta }{\partial t \partial \overline{ t} }(0)=
- \frac{ 1}\pi \ \int_{ \partial R(0)} k_2(0,z) \bigl|
\frac{\partial q(0,z)}{\partial z }\bigr|^2ds_z
 - \frac{ 4}\pi \ \iint_{ R(0)}
\bigl| \frac{\partial^2 q }{\partial \overline{ t} 
\partial z} (0,z)\bigr|^2dx\,dy,
$$
which is desired.

\vspace{1mm} {\it Case when $R(t)$ is of genus $g\ge 1$.}\quad 
We put $R'(0)= R(0 ) \setminus \cup_{k=1}^g (A_k(0)\cup B_k(0))$ 
and $\widehat R'(0)= R'(0)\cup V$,
so that $R'(0)$ and $\widehat R'(0)$ are  planar Riemann surfaces such that  
$$
\partial R'(0)=\partial R(0)+ \sum_ {k=1}^g (A_k^+(0)+ A_k^-(0))
+\sum_ {k=1}^g (B_k^+(0) + B_k^-(0)).
$$
Here $A_k^+(0)(A_k^-(0))$ is the same(opposite) direction of $A_k(0)$, 
and $B^+_k(0)$
\,$(B^-_k(0))$ is similar. For $t\in B$, if we restrict 
the branch $q^*(t,z)$ (with $q^*(t,z^0)=0$) to 
 $R'(0) \setminus \{0,1\}$, then $
q^*(t,z')= q^*(t, z'')$ mod $2\pi$ for $z', z''$ over the same 
point $z\in \widehat{ R}'(0)$. Hence 
$\frac{\partial q^* }{\partial t}(0,z), \ \frac{\partial ^2 q^* 
}{\partial \overline{ t}\partial z }(0,z)$ are single-valued harmonic 
functions on $\widehat R'(0)$, so that
\begin{alignat*}{3}
&\int_{ \partial R(0)} \frac{\partial q^*}{\partial t}(0,z) 
 \frac{\partial ^2 q^*}{\partial \overline{ t} \partial z }(0,z)dz\\
&=
\iint_{ R'(0)} 
d\bigl(\,\frac{\partial q^*}{\partial t}(0,z) 
 \frac{\partial ^2 q^*}{\partial \overline{ t} \partial z }(0,z)dz\,\bigr)
- \sum_ {k=1}^g  \int_{ A^\pm_k(0)+ B^\pm _k(0)} \frac{\partial q^*}{\partial t}(0,z) 
 \frac{\partial ^2 q^*}{\partial \overline{ t} \partial z }(0,z)dz
\\
&=: J_1-J_2.
\end{alignat*}
Since $\frac{\partial q^*}{\partial \overline{ t}\partial z }(0,z) $ is 
holomorphic on ${ R}'(0) $,
 we have by Green's formula 
\begin{align*}
 J_1
&=2i\,\iint_{ R(0)} \bigl|  \frac{\partial^2 q }{\partial t \partial 
\overline{ z} }(0,z) \bigr|^2 dx \, dy; \\
 J_2(A_k)&:= 
\int_{ A_k^\pm(0)} \frac{\partial q^*}{\partial t } (0,z)
 \frac{\partial ^2q^*}{\partial \overline{ t } \partial z }  (0,z) dz\\
&=\int_{ A_k(0)} \bigl(\,\frac{\partial q^*}{\partial t } (0,z^+)- 
\frac{\partial q^*}{\partial t } (0,z^-) \bigr)\,
 \frac{\partial ^2q^*}{\partial \overline{ t } \partial z }  (0,z) dz.
\end{align*}
By (\ref{akbk}) and $\int_{ C_j(0)}*dq(t, z)=0, j=1,\ldots , q$, 
we have, for $z^\pm$ over any  $z\in A_k(0)$,
\begin{align*}
q^*(t,z^+)- q^*(t,z^-)&=  \int_{ B_k(0)} *d q(t,\zeta)   \quad \mbox { mod} \ 2\pi. \\
\therefore  \ \  \frac{\partial q^*}{\partial t}(t,z^+)-\frac{\partial q^*}{\partial t}(t,z^-)&= \frac{\partial}{\partial t} \int _{B_k(0)} * dq(t,\zeta),
\end{align*}
which is independent of $z \in A_k(0)$. It follows from 
$\frac{\partial q^*(t,z)}{\partial z}dz=
\frac{ 1}2 (*dq(t,z)-i\,dq(t,z)) $  that 
\begin{align*}
J_2(A_k)&=
\bigl[
\frac{\partial }{\partial t}\bigl(\int_{ B_k(0)} * dq(t,\zeta)\bigr)\bigr]_{t=0}
\cdot \bigl[
 \frac{\partial }{\partial \overline{ t} } \bigl(\int_{ A_k(0)} 
 \frac{\partial q^*(t,z) }{\partial z}dz\bigr)\bigr]_{t=0}\\
&= 
\frac{ 1}2\  \frac{\partial {\bf  b}_k}{\partial t}(0) \cdot \frac{\partial 
 {\bf  a}_k}{\partial \overline{ t}}(0).
\end{align*}
 By $B_k(0) \times A_k(0)= -1$, it similarly holds 
$ J_2(B_k)
=-\frac{ 1}2 \frac{\partial {\bf  a}_k}{\partial t}(0) \cdot \frac{\partial 
 {\bf  b}_k}{\partial \overline{ t}}(0)$,
so that 
$
J_2(A_k)+J_2(B_k)= -i \ \Im\ \bigl\{ \frac{\partial {\bf  a}_k}{\partial t}(0)\cdot
\frac{\partial {\bf  b}_k}{\partial \overline{ t}}(0)  \bigr\}. 
$
We thus have
\begin{alignat*}{3}  
&\Im\,  \bigl\{
\int_{\partial R(0)} \frac{\partial q^*}{\partial t } (0,z) \frac{\partial ^2
q^*}{\partial \overline{ t} \partial z}(0,z) dz
\bigr\}\\ & \ \ =\Im\ \bigl\{ \, J_1 - \sum_ {k=1}^g (J_2(A_k)+ J_2(B_k))
\ \bigr\}  \\
&\quad = 2 \iint_{ R(0)} \bigl| \frac{\partial^2 q }{\partial 
 \overline{ t } \partial z }(0,z) \bigr|^2\,dx\,dy +
 \Im\, \bigl\{\sum_ {k=1}^g
  \frac{\partial {\bf  a}_k}{\partial t}(0)
\cdot
 \frac{\partial {\bf  b}_k }{\partial \overline{ t}}
(0)\bigr\}.
\end{alignat*}
This with (\ref{eqn:kihon-2}) 
completes the second formula in the 1st step.
 
\vspace{1mm} 
{\it  $2^{nd}$ step.  \  
Lemma 
\ref{vari-form:L_0} is true in the general case .}

In fact, it suffices to prove 
Lemma \ref{vari-form:L_0} at $t=0$. 
If necessary, take a smaller disk $B$ of center $0$. Then we find a 
linear transformation
$T: \, (t,z)\in B \times \mathbb{ P}_z \mapsto (t,w)= (t, f(t,z)) \in B 
\times \mathbb{P}_w$ such that  
$f(t,0)= 0$; $\frac{\partial f}{\partial z 
}(t,0)=1 $; $f(t,\xi(t))=  \ \mbox{const.} \ c$ for $t\in B$, and 
${\mathcal  D}: \, =T({\mathcal  R})$
is an unramified  domain over $B\times \mathbb{C}_w$.
We write $D(t)=f(t,R(t)),\,t\in B$, so that  ${\mathcal  D}
=\cup_{t\in B}(t, D(t))$ and ${\mathcal  D}$ has two 
constant sections $\Theta_0:\,w=0 \mbox{ and } \Theta_c: \ w=c$ 
(the pull backs of $\Xi_0$ and $\Xi_\xi$ by 
$T$), hence the variation
${\mathcal  D}: t\in B \to D(t) $
is a case in the $1$st step. 
For $t\in B$, we consider the $L_0$-principal function 
$\widetilde q (t,w)$ 
and the $L_0$-constant $\widetilde {\beta}(t)$ for $(D(t), 0,c)$, so that 
\begin{align} \nonumber 
 &\widetilde  q(t,w)= \log \mbox{$\frac{ 1}{|w|}$}+  \widetilde
{\mathfrak h}_0(t,w)
 \ \ \quad \quad \qquad \ \  \mbox{ in $U_0(t)$};
\\
&\widetilde q(t,w)= \log |w-c|+  \widetilde {\beta}(t) + \widetilde 
{\mathfrak h}_c(t,w) 
 \quad \mbox{in $U_c(t)$},  \nonumber 
\end{align}
where $\widetilde {\mathfrak h}_0(t,0), \, \widetilde{\mathfrak h}_c(t,c) \equiv 0$ on $B$.
 We put 
$
{\widetilde   A}_k(t)= f(t, A_k(t))$ and ${\widetilde  B}_k(t)= f(t, 
B_k(t))
$ on $D(t)$ 
which continuously vary
 in ${\mathcal  D}$ with $t\in B$ without passing through $w=0,c$.
Since 
\begin{align*}
\nonumber 
w=f(t,z)= \left\{
\begin{array}{lll} 
 z + b_2(t)z^2 + \cdots 
&\mbox{at $z=0$}; \\[3mm]
 c+ a_1(t) (z- \xi(t)) + a_2(t)(z- \xi (t))^2 + \cdots &
\mbox{at $z=\xi(t)$},
\end{array}
\right.
\end{align*}
where $a_1(t)\ne 0, a_2(t), \ldots;  b_2(t), \ldots $ are holomorphic on $B$,
we have $q(t,z)= \widetilde q(t, f(t,z))$ in ${\mathcal  R}$, i.e., 

\vspace{1mm} 
$\qquad q(t,z)= \log |f(t,z)-c| + \widetilde{\beta}(t) + \widetilde {
\mathfrak{ h}}_c(t, f(t,z)) 
 \quad \mbox{ at $z=\xi(t)$},
$

\vspace{1mm} 
\noindent 
so that
\begin{align} \label{eqn:bt-tbt}
\beta(t)&= \widetilde { \beta} (t)+ \log |a_1(t)|;\\[-2mm] 
 \mathfrak{ h}_{\xi}(t,z)&= \widetilde{ \mathfrak{ h}}_c(t, f(t,z)) + \log\bigl|\ 1+ \frac{ a_2(t)}{a_1(t)}
(z-\xi(t)) + \ldots\ \bigr|.\nonumber 
\end{align}
Let $\psi (t,w)$ be  a $C^\omega$ defining function of  
$\partial {\mathcal  D}$. Then $ \varphi (t,z):=\psi(t, f(t,z))$
is that of $\partial {\mathcal  R}$, so that  
 we 
have for $w=f(t,z)$
\begin{align*}
& k_1(t,z)= \frac{ \frac{\partial \varphi(t,z) }{\partial t }}{|\frac{\partial 
\varphi (t,z) }{\partial z}| }= \frac{ \widetilde k_1(t,w)}{
|\frac{\partial f(t,z)}{\partial z}|} + 
\frac{ \frac{\partial f(t,z)}{\partial t}}{
|
\frac{\partial f(t,z)}{\partial z}|}  \cdot 
\frac{  \frac{\partial \psi }{\partial w}(t, w)}{|
\frac{\partial \psi }{\partial w}(t, w)|}, \ \ \quad (t,z)\in \partial {\mathcal  R}.\\[3mm]
\therefore \ \ 
 \ & \int_{ \partial R(0)} k_1(0,z) |\frac{\partial q 
(0,z)}{\partial z }|^2 ds_z \\ \nonumber 
&=
\int_{ \partial R(0)} \frac{ \widetilde k_1(0,w)}{|\frac{\partial 
f(0,z)}{\partial z}| } 
\ \bigl|\frac{\partial q(0,z)}{\partial z}\bigr|^2\ ds_z 
+ 
 \int_{ \partial R(0)}
\frac{ \frac{\partial f}{\partial t}(0,z)}
{|\frac{\partial 
 f(0,z)}{\partial z}|}\cdot 
\frac{ \frac{\partial \psi}{\partial w}(0,w) }
{|
\frac{\partial \psi}{\partial w} (0, w)|}
\ \bigl|\frac{\partial q(0,z)}
{\partial z}\bigr|^2\ ds_z\\[3mm]
&=:J_1+J_2. \nonumber 
\end{align*}
Since $\frac{\partial \widetilde q(0,w)}{\partial w} \frac{f(0,z)}{dz}= \frac{\partial q 
(0,z)}{\partial z}dz$, we have by  the 1st step and 
(\ref{eqn:bt-tbt}) 
\begin{align*}
J_1&= \int_{ \partial D(0)} \widetilde k_1(0,w) 
\bigl|\frac{\partial \widetilde q(0,w)}
{\partial w}\bigr|^2\ ds_w
=-\pi \frac{\partial \widetilde \beta }{\partial t}(0)
=-\pi\,( \frac{\partial \beta}{\partial  t}(0) - \frac{ 1}2 \frac{ a_1'(0)}{a_1(0)}).
\end{align*}
For a fixed $t\in B$, we write $z=g(t,w):=f^{-1}(t,w)$.
We put $\widetilde C_j(0)=f(0, C_j(0))$ and $\widetilde V_j=f(0, V_j),
 j=1,\ldots , \nu$, and 
consider the  single-valued conjugate
 harmonic function  $\widetilde q^{\,*}_j(0,w)$  of 
$\widetilde q(0,w)$ in $\widetilde V_j$ which vanishes 
on $\widetilde C_j(0)$. Then we 
find a function  $k(w)\in C^\omega(V_j)$ such that  
$\widetilde q_j^{\,*}(0,w)= k(w) \psi(0,w)$ in 
$\widetilde V_j$. This and the residue theorem imply
 \begin{align*}
J_2&= - \int_{ \partial D(0)} 
\frac{
\frac{\partial g}{\partial t}(0,w)
}
{ \frac{\partial g(0,w)}{\partial w} }\
\frac{
 \frac{\partial \psi(0,w)}{\partial w}
}
{
|\frac{\partial \psi(0,w)}{\partial w} |
}
\ \bigl|\frac{\partial \widetilde q^*(0,w)}
{\partial w}\bigr|^2\ ds_w\\
  &=
 i \int_{ \partial D(0)} 
\frac{
\frac{\partial g}{\partial t}(0,w)
}
{ \frac{\partial g(0,w)}{\partial w} }\
(\frac{\partial \widetilde q_j^{\,*}(0,w)}
{\partial w})^2\ dw\\
&= 2\pi
\mbox{ Res}\,_{w=0,\,c}\  \bigl\{
\ \frac{
\frac{\partial g}{\partial t}(0,w)
}
{ 
\frac{\partial g(0,w)}{\partial w} 
}\
(
\frac{\partial \widetilde q(0,w)}
{\partial w}
)^2 \bigr\} \\
&=
2\pi\, (\,\frac{\partial  \mathfrak{ h}_\xi}{\partial z}(0, \xi(0)) 
 \, \xi'(0) - \frac{ 1}4
 \frac{ a_1'(0)}{a_1(0)}\,).\\
\therefore \quad & J_1+J_2=-\pi\, \bigl(\, \frac{\partial \beta}{\partial t}(0)  -2\ 
  \frac{\partial \mathfrak{ h}_\xi}{\partial z}(0,\xi(0)) \ \xi'(0)\,\bigr),
\end{align*}
 which is identical with the first formula in the 2nd step.

To prove the second one, we have from the 1st step
\begin{alignat}{3} \nonumber
\frac{\partial ^2 \widetilde  \beta}{\partial t \partial \overline{ t} }(0)&=
-\frac{ 1}{\pi} \ \int_{\partial D(0)} \widetilde k_2(0,w)
\bigl|
\frac{\partial \widetilde  q(0,w)}{\partial w}\bigr|^2 ds_w -
 \frac{ 4}\pi \ \iint_{ D(0)} \bigl| \frac{\partial^2 \widetilde  q }{\partial 
 \overline{ t } \partial w }(0,w) \bigr|^2\,du\,dv   \\[0mm]  
 \nonumber 
\nonumber  & \quad -
 \frac{ 2}{\pi}\ \Im\,   \sum_ {k=1}^g
 \bigl[\frac{\partial }{\partial t}
\int_{ {\widetilde A}_k(t)} * d \widetilde  q(t,w) \bigr]_{t=0} \
\cdot \bigl[
\frac{\partial }{\partial \overline{ t}} 
\int_{ {\widetilde  B}_k(t)} * d \widetilde q(t,w) \bigr] _{t=0},
\end{alignat}
where $\widetilde k_2(t,w)$ is the Levi curvature of $\partial {\mathcal  
D}$. It suffices
to show that each term of the above formula  is invariant 
for  $T:(t,z)\in {\mathcal  R} \to (t,w)=(t,f(t,z))\in 
{\mathcal  D}$,\,i.e., it holds for $t\in B$ 
 \begin{alignat*}{3}
&i.\ \ \frac{\partial ^2 \widetilde \beta(t)}{\partial t \partial \overline{ t}}
=  \frac{\partial ^2 \beta (t)}{\partial t \partial \overline{ t}};\\
&ii. \ \ \int_{\partial D(t)} \widetilde k_2(t,w)\bigl|
\frac{\partial \widetilde  q(t,w)}{\partial w}\bigr|^2 ds_w=
\int_{\partial R(t)}  k_2(t,w)\bigl|
\frac{\partial  q(t,z)}{\partial z}\bigl|^2 ds_z;\\
&iii. \ \ 
 \iint_{ D(t)} \bigl| \frac{\partial^2 \widetilde  q }{\partial 
 \overline{ t } \partial w }(t,w) \bigr|^2\,du\,dv=
 \iint_{ R(t)} \bigl| \frac{\partial^2 q }{\partial 
 \overline{ t } \partial z }(t,z) \bigr|^2\,dx\,dy;\\
&iv.  \ \ \frac{\partial }{\partial t} 
\int_{ {\widetilde  A}_k(t)} * d \widetilde  q(t,w)=
\frac{\partial }{\partial t} 
\int_{ {A}_k(t)} * d q(t,z),  \  \mbox{and similar  
  for  ${\widetilde  B}_k(t)$ and $B_k(t)$.} 
\end{alignat*}

In fact, $i.$ is clear from (\ref{eqn:bt-tbt}). 
Since $\widetilde 
q(t, w)= q(t,z)$ (where $w=f(t,z)$ for $(t,z)\in {\mathcal  R}$) is 
harmonic  on each $R(t), \, t\in B$, we have 
$iii.$ and $iv.$ 
Further, by the simple calculation, we see, in general,  
that {\it  $k_2(t,z)\, \frac{1}{|dz|} $ on 
$\partial {\mathcal  R}$ is invariant under all holomorphic transformations $T$ of the form }
$T: (t,z)\in {\mathcal  R}\cup \partial {\mathcal  R }
\mapsto (t,w)=(t,f(t,z))\in \widetilde 
{\mathcal R}\cup \partial \widetilde {\mathcal  R} $,
\,i.\,e.,
$
\widetilde  k_2(t,w)= k_2(t,z) \bigl|\frac{\partial f(t,z)}{\partial z}
 \bigr|$. It follows that 
\hspace{-1cm}\begin{align*}
  \widetilde k_2(t,w)\bigl|
\frac{\partial \widetilde  q(t,w)}{\partial w}\bigr|^2
 |dw|
=k_2(t,z)\bigl|
\frac{\partial q(t,z)}{\partial z}\bigr|^2 |dz|
\end{align*}
for $ z\in \partial R(t)$ and $w= f(t,z)$. This implies $ii.$ We complete the proof of  Lemma \ref{vari-form:L_0}.
\hfill $\Box$ 

\vspace{2mm}  As noted in \cite{hamano-2},
since ${\mathcal  R}$  is pseudoconvex in 
$\widetilde {\mathcal  R}$ iff  $k_2(t,z) \ge 0$ on $\partial 
{\mathcal  R}$,  Lemma \ref{vari-form:u} implies that,  if
${\mathcal  R}$  is pseudoconvex in 
$\widetilde {\mathcal  R}$, then the   
$L_1$-constant $\alpha (t)$ 
for $(R(t), 
0, \xi(t))$
 is $C^\omega$ subharmonic on $B$, while
 Lemma 
\ref{vari-form:L_0} makes the following  contrast with it:
 \begin{theorem}
\label{th:sup-L_0} If ${\mathcal R}$ is pseudoconvex in $\widetilde {\mathcal R}$ 
and $R(t), 
  t\in B$ is planar,
then the $L_0$-constant $\beta (t)$ for $(R(t), 0, \xi(t))$
is $C^\omega$ superharmonic on $B$. 
\end{theorem}

\section{Harmonic span and its geometric meaning}

\noindent 
{We  a little recall the slit mapping theory  in  one complex variable. 
Let $R$ be a planar Riemann surface sheeted over $\mathbb{C} _z$ bounded by a finite number of 
smooth contours $C_j, 
j=1,\ldots , \nu$. 

Let $R\ni 0$ and let ${\mathcal  U} (R)$ denote the set of all 
univalent functions $f$ on $R$ 
such that 
$f(z)- \frac{ 1}z$ is regular at $0$. For $w=f(z) \in {\mathcal  
U}(R)$ we consider the Euclidean area $E(f)$ of $\mathbb{ C}_w \setminus f(R)$ 
and put 
$$
{\mathcal  E}(R)= \mbox{ sup}\,  \{E(f): f\in {\mathcal  U}(R)\}.
$$
Due to P.\,Koebe (see Chap.\,X 
in \cite{ford}),  we have two special ones $w=f_i(z), 
 i=1,0$ in ${\mathcal  U}(R)$ such that  $f_1(R) (f_0(R))$
 is 
 a vertical (horizontal) slit  univalent domain in $\mathbb{ P}_w$.
 In his pioneering work \cite{grunsky},  H. Grunsky  showed in p.\,139-140: if we consider  
$$
g:=\displaystyle{ \frac{ 1}2\, (f_1+f_0)} \quad \mbox{ on  $R$}
$$
and $K_j=-g(C_j), j=1,\ldots ,\nu$, then 
each $K_j$ bounds an unramified domain $G_j$ over $\mathbb{ 
C}_w$ such that, if we denote by  $E_j(g)$ the 
Euclidean (multivalent) area of $G_j$ and put $E(g)=\sum_{j=1}^\nu E_j(g)$, then 
$E(g)\ge {\mathcal  E}(R)$.  Then,  in his substantial work \cite{schiffer},  
M. Schiffer  in p.\,209 introduced  
the following quantity 
 $S(R)$, which he named the {\it  
 span } 
 for $R$,
$$
S(R):=
\Re\, \{a_1 -b_1\},
$$
where $a_1$ and $b_1$ are the coefficients of $z$ (the first degree) of 
the Taylor expansions of $f_1(z)-\frac{1}z$ and $f_0(z)- \frac{ 1}z$ at 
$0$, respectively, and showed the following beautiful results 
(p.\,216 in \cite{schiffer}):\, $g \in 
{\mathcal  U}(R)$;
 {\it  each $G_j, j=1,\ldots ,\nu$ is
a convex domain in $\mathbb{ C}_w$, and }
$$
E(g)={\mathcal  E}(R)={ \frac{ \pi}{2} S(R)}.
$$
His proofs  were rather intuitive and short.  The precise proofs 
are found in \S 12,  Chap.\,III in \cite{ahlfors-sario}. 

Let $\xi \in R, \ \xi \ne 0$ and let  ${\mathcal  S}(R)$ denote the set of all univalent functions 
$f$ on $R$ 
such 
 that   $f(z)- \frac{ 1}z$ is regular at $0$ and $f(\xi)=0$, say
\begin{alignat*}{3}
&f(z)= c_1(z-\xi)+ c_2(z-\xi)^2+\ldots \quad &\mbox
{ at }\  \xi.
\end{alignat*}
We then put $c(f)=c_1 (\ne 0)$. We draw a simple curve $l$ on $R$ from  
$\xi$ to  $0$. Let $f\in {\mathcal  S}(R)$ and  $w=f(z)$ on $R$. 
 Then   $f(l)$ 
is a simple curve from $0$ 
to $\infty$ in $\mathbb P_w$, and 
 each branch of $\log f(z)$ on $R \setminus l$ is 
single-valued {\it univalent}. Fix one of them, say $\tau=\log f(z)$. 
Consider the  Euclidean area $E_{\log} (f) \,( \ge 0)$ 
of the complement of $\log f(R \setminus 
l)$ in $\mathbb{ C}_{\tau}$ and put 
\begin{align*}
{\mathcal  E}_{\log}(R)= \sup\,\{ E_{\log}(f): f\in {\mathcal  S}(R)\}.
\end{align*}

Now let $p(z)\,(\,q(z)\,)$ and $\alpha \,(\beta)$ be 
the  $L_1$-($L_0$-)principal function and $L_1$-($L_0$-)constant for $(R, 
0, \xi)$. We choose the harmonic conjugate $p^*(z)$ $(q^*(z))$ on $R$
such that, if we put   $
P(z)= e^{p(z)+ip^*(z)}(\, Q(z)=e^{q(z)+iq^*(z)})$ on $R$, then 
$P(z)- \frac{ 1}z \ 
(\,Q(z)-\frac{ 1}z)$ is regular at $0$.
Then  $P, Q\in {\mathcal  S}(R)$ and $w=P(z) \, (\,Q(z))$ 
is a circular (radial) slit mapping with 
$\log |c(P)|=\alpha$ $(\,\log|c(Q)|= {\beta 
}$) and $E_{\log}(P)=E_{\log}(Q)=0 $.
We see in \S 13, Chap.\,III in \cite{ahlfors-sario} that 
 {\it  $P$ maximizes $2\pi 
 \log |c(f)| +E_{\log}(f)$, while $Q$ minimizes 
$ 2\pi 
 \log |c(f)| -E_{\log}(f)$ among ${\mathcal  S}(R)$.
} 
 
On the other hand,  M.\,Nakai  expected that the quantity
\begin{align}
 \label{nakai}
 s(R):= \alpha-\beta
 \end{align}
will have some gemetric meaning.  
 In  \cite{nakai-sario} he  named $s(R)$ 
 the {\it  harmonic span} for $( R, 0, \xi)$.
 Hereafter in this paper we shall show that $s(R)$ has 
some remarkable properties 
not only in one complex variable but also in several complex variables.

We precisely write 
\begin{eqnarray}
\label{eqn:to-mapping}
\begin{array}{lll}
&P(z)= e^{\alpha+i\theta_1}\ (z-\xi)+ \mbox{$\sum_{n=2}^\infty
 a_n(z-\xi)^n$ \quad at } \ \xi;\\[2mm]
&Q(z)= e^{\beta+i\theta_0}\ (z-\xi)+ 
 \mbox{$\sum_{n=2}^\infty b_n(z-\xi)^n $ \quad at } \ \xi,
\end{array}
\end{eqnarray}
where $\theta_1,\theta_0$ are certain constants. 
 We put 
\begin{align*}
 {  D}_1&:= P(R)=
\mathbb{P}_{w} \setminus 
\cup_{j=1}^\nu P(C_j)=
 \mathbb{P}_w \setminus 
\cup_{j=1}^\nu \ \mbox{{arc}$\{A_j^{(1)},A_j^{(2)}\}$}
;\\[2mm]
{ D}_0&:= Q(R)=\mathbb P_w \setminus \cup_{j=1}Q(C_j)=
\mathbb{P}_w\setminus \cup_{j=1}^\nu \mbox{{ segment}$\{
B_j^{(1)}, B_j^{(2)}\}$}.
\end{align*}
Here 
\begin{eqnarray}
\label{eqn:tjrj-notation}
\qquad \begin{array}{lll}
  &{\rm arc}\{A_j^{(1)},A_j^{(2)}\}&=
\{r_j e^{i\theta}:  
\theta_j^{(1)}\le \theta \le 
\theta_j^{(2)}\};
\\[2mm]
&{\rm segment}\{B_j^{(1)},B_j^{(2)}\}&=
 \{r e^{i\theta_j}:  0<r_j^{(1)}\le r \le 
 r_j^{(2)}<\infty\},
\end{array}
\end{eqnarray}
where $0<\theta_j^{(2)}-\theta_j^{(1)}<2\pi$ and 
 $r_j, \theta_j^{(k)},  \theta_j, r_j^{(k)} \, (j=1,\ldots , \nu; k=1,2)$ are
constants. We take the points  
  $a_j^{(k)}, b_j^{(k)}\in C_j\  $ such 
that 
\begin{eqnarray}
 \label{eqn:ajk-bjk}
P(a_j^{(k)})=A_j^{(k)}, \qquad  \quad 
Q(b_j^{(k)})=B_j^{(k)}.
\end{eqnarray}

By conditions $(L_1)$ and $(L_0)$ for  $p(z)$ and $q(z)$,  $\sqrt{P(z)Q(z)} $
 consists of two single-valued branches 
$H(z)$ 
and  
$-H(z)$ on $R$ where  
$H(z)$ 
has only 
one pole at $z=0$ such 
that   $H(z)-\frac{ 1}z$ is regular at $0$, and $H(z)$ has $0$ only at $z=\xi$. We write
 \begin{align*}
  H(z)=\sqrt{P(z)Q(z)} \quad \mbox{ on $R$}.
\end{align*}
 Each branch of $\log P(z)(\log Q(z))$ 
is also single-valued univalent on $R \setminus l$, while 
$\log H(z)$ is single-valued but not univalent so far. 
 We choose  three branches  in $R \setminus l$ such that  
$$
 \tau=\log H(z)= \frac{ 1}2 (\log P(z) +\log Q(z)).
$$
We fix a tubular neighborhood $V_j$ of 
each contour $C_j$ with $V_i\cap V_j=\emptyset\ (i\ne j)$ and $V_j \not\ni 0, \xi$, where 
$\log H(z)$ on $V_j$ is single-valued. 
Then we have the following geometric meaning of the harmonic span $s(R)$:
\begin{theorem} \label{lem:P0P1} ${}$
\begin{enumerate}
\item [1.]
 Each $ -(\log H)(C_j), \,j=1,\ldots , \nu$ 
 is a  $C^\omega $ convex 
curve in $\mathbb{ C}_{\tau}$, and $-H(C_j)$ is a $C^\omega$ simple closed curve in 
      $\mathbb{ C}_w$.
 \item [2.]  $H\in {\mathcal  S}(R)$ and 
$E_{\log} (H)={\mathcal  E}_{\log}(R)= \displaystyle{ \frac{ \pi}2}\,  s(R).$
 \item [3.]  Assume that $R$ is simply connected, and  
let $d(0, \xi)$ denote the Poincar\'e distance between $0$ and 
$\xi$ on $R$.
Then 
$$
s(R)= 4 \log  \cosh  d(0,\xi).
$$
\end{enumerate}
\end{theorem}
}
The method in the  proofs in Chp.\,III in \cite{ahlfors-sario} of M.\,Schiffer's results seems to  have some gaps  to 
prove 1. and 2. in Theorem \ref{lem:P0P1}. 
 We get over them by the idea of using the Schottky 
double (compact) Riemann surface $\widehat{R}$ 
 of $R$. We also apply this idea to prove Corollary 
\ref{cor:rigidity} for the variations of Riemann surfaces.

\vspace{2mm}  
\noindent {\it  Proof of Theorem \ref{lem:P0P1}.} \ 
Similarly  to $F:=\frac{ df_1}{df_0}$ used in p.\,182 in 
\cite{ahlfors-sario} (cf: (25) in \cite{schiffer}), we
 consider  the following function 
\begin{equation}
  \label{eqn:gur-function}
W=F(z):= 
\frac{ d\log Q}{d\log P}, \qquad \mbox{ $z\in R\cup \partial R$},
\end{equation}
 which 
is a single-valued meromorphic 
function on $R$  such that  $\Re\  F=0 $ on $\partial R$, since $\log P(C_j)$ is 
 a vertical segment and $\log Q(C_j)$ is a
 horizontal segment in $\mathbb{ C}_\tau$. It follows from the Schwarz
reflexion principle that $F$ is meromorphically extended to the Schottoky double 
Riemann surface 
$\widehat{ R}=
R\cup \partial R \cup R^*$ of $R$ such that  $F(z^*)=-\overline{ F(z)}$, 
where $z^*\in R^*$ is the reflexion point of $z\in R$.
Fix $C_j, j=1,\ldots ,\nu$. Since 
 $\Re\, \log P(z)= p(z)$ and $\Re\, \log 
Q(z)= q(z)$ on $R$, 
we have
\begin{align}\label{eqn:PQ}
 \log P(z)= u_1(z)+i v _1(z), \quad \log Q(z)= u_0(z)+iv_0(z), \ \ z\in V_j, 
\end{align}
where $u_1(z) (v_0(z))= \mbox{ const.}\,c_1\,(c_0)$ on $C_j$.
Then $\mathfrak{ C}_j:=\log H(C_j)$ is a closed (not necessarily simple 
 so far) curve in 
$\mathbb{ C}_\tau $:
 \begin{align}
  \label{eqn:logH}
\tau = \frac {1}2\,(\,c_1+u_0(z)) + \frac{ i}2(\,c_0 + v_1(z)), \quad 
z\in C_j.
\end{align}
Using  notation (\ref{eqn:ajk-bjk}), we shall 
show:
\begin{enumerate}
      \item [i)] $\{a_j^{(k)},b_j^{(k)}\}_{k=1,2}$ are 
$4$ distinct points, which necessarily line cyclically, for example, $(a_j^{(1)}, b_j^{(1)}, 
	    a_j^{(2)}, b_j^{(2)})$
on $C_j$;
 \item [ii)] the zeros  of $F(z)$ are 
$\{b_j^{(k)}\}_{j=1, \ldots , \nu;\,k=1,2}$ of order one, and 
 the poles are  $\{a_j^{(k)}\}_{j=1,\ldots , \nu;\,k=1,2}$ of 
order one;
 \item [iii)]  the curve $\mathfrak{ C}_j $ is locally non-singular
  in $\mathbb{ C}_\tau$;
 \item [iv)] $\Re\,F(z)>0$ on $R$; 
  \item[v)] at any $\tau\in \mathfrak{ C}_j$, the
 curvature $\frac{ 1}{\rho_j(\tau)}$ of $\mathfrak{ C}_j$ is 
  negative.
\end{enumerate}

We divide the proof into two steps.

\vspace{1mm} 
{\it  $1^{st}$ step. \  If 
we admit { i)},  then { ii)} $\sim$ {v)}  hold.}

  In fact, i) clearly implies iii). Since $P(z)(Q(z))$ is univalent on $R$ 
with the  circular (radial) 
 slit boundary condition,
 we have  $F(z)\ne 0, \infty$
on $R \cup R^*$ and 
$F(z)$ has zeros at most $b_j^{(k)}$ and 
poles at most $a_j^{(k)}$, of order one.
Thus, i) implies ii). 
Further,   i) implies  that $W=F(z)$ 
is locally one-to-one in a neighborhood of 
at any $z\in C_j$ even at $a_j^{(k)}, b_j^{(k)} (k=1,2)$,
so that 
$F$ is a meromorphic function on $\widehat{ R}$ 
of degree $2\nu$.
Hence, for each fixed $j=1,\ldots, \nu$, if $z$ travels   $C_j$ all once, 
then $F(z)$ travels the imaginary axis all just {\it twice}. 
It follows  that $F(\widehat{ R})$ is a 
$2\nu$  sheeted compact Riemann surface  over $\mathbb{ P}_W$  
 with $2(2\nu +g-1)$ branch points lying on $\mathbb{ P}_W  \setminus 
 \{\Re\, W=0\}$, and is divided by $\nu$ closed curves $F(C_j), j=1,\ldots ,\nu$
 into two connected parts over $\Re\,  W>0$  and $\Re\,  W<0$. 
Since $F(0)=1$, we have 
 $\Re\,F(z)>0 $ on $R$ and $\Re\,F(z)<0$ 
on $R^*$, which is iv).
To prove v), fix  $p_0\in C_j$  and 
take a local parameter 
$z=x+iy$ of a neighborhood $V$ of $p_0$ such that  $p_0$ 
corresponds to $z=0$ and the oriented arc $C_j \cap V$ corresponds to 
 $I:=(-\rho, \rho)$ on the $x$-axis. Then using this parameter, we 
 see from  $\Re\,F(z)>0 $ on $R$ that 
\begin{align}
 \label{eqn:partialF}
 \Im F'(x)=\Im \frac{ \partial 
 F(x)}{\partial x}<0 \ \  \ \mbox{ on  $I$}.
\end{align}
\noindent By (\ref{eqn:logH})
the subarc $\Gamma _j:= 
\log H(I)$ of  $\mathfrak{ C}_j$ in $\mathbb{ C}_\tau$ is 
of the form:
$$
\tau =u(x)+iv(x)= \frac{ 1}2[(c_1+ u_0(x))+i(c_0+v_1(x))], 
\ \ x\in I.
$$ 
Since  
the arc $\Gamma _j$ is locally  non-singular by iii),  we calculate 
the curvature ${ 1}/{\rho_j(x)}$ at the point $(u(x),v(x))$ of $\Gamma _j$\,:
$$
\frac{ 1}{\rho_j(x)}= \frac{ v''(x)u'(x)-v'(x)u''(x)}
{(v'(x)^2+u'(x)^2)^{3/2}}= 
\frac{ v_1''(x)u_0'(x)
-v_1'(x)u_0''(x)}
{(v_1'(x)^2+u_0'(x)^2)^{3/2}}.
$$
On the other hand, by (\ref{eqn:PQ}) we have, for $x\in I \subset C_j$,  
\begin{align*}
\Im\, F'(x)&= \Im \
\bigl\{
\frac{ d }{dx }
\bigl( 
\frac{ 
\frac{d u_0(x)}{dx}+i\frac{ dc_0}{dx}
}
{ 
\frac{ dc_1}{dx}+i \frac{ dv_1(x)}{dx}
}
\bigr)
\bigr\}
=
\frac{v_1''(x) u_0'(x)-v_1'(x)u_0''(x) }{v_1'(x)^2}.\\
\therefore \ \ \frac{ 1}{\rho_j(x)}&=   \frac{v'_1(x)^2}{
(v_1'(x)^2+ u_0'(x)^2)^{3/2}} \cdot \Im F'(x).
\end{align*}
Since $v'_1(0)=0$ iff $x=a_j^{(k)}$, 
(\ref{eqn:partialF}) proves v) for $p_0 \ne a_j^{(k)}$.
For $p_0=a_j^{(k)}$, since $v'_1(0)=0, \ v''_1(0),\ u_0'(0)
\ne 0$ under i),  $ v_1'(x)^2 \cdot \Im F'(x)$ is regular and 
$\ne 0$. Hence $\frac{ 1}{\rho_j(p_0)}<0$, which proves v). 

{\it $2^{nd}$ step. \ i) is true.}

In fact, assume that
  $R$ does not satisfy i). 
Clearly it does not occur   $\{a_j^{(1)}, a_j^{(2)}\}=
\{b_j^{(1)}, b_j^{(2)}\}$ for any $j$, so that 
$\{a_j^{(1)}, a_j^{(2)}\}\cap
\{b_j^{(1)}, b_j^{(2)}\}$ consisits of one point 
 for some $j$, say 
$j=1,\ldots, \nu'(\le \nu)$. We denote by $o_j$ such one point on $C_j$. 
Hence each  $\mathfrak{ C}_j:= \log 
H(C_j), j=1,\ldots, \nu'$  is a closed curve in $\mathbb{ C}_\tau$
 with only one singular point 
at $\mathfrak{ o}_j:= \log H( o_j)$ and $F$ is  a meromorphic function of degree 
$2\nu-\nu'$ on $\widehat{ R}$. By the  same reasoning as in the $1$st 
step, if $z$ travels $C_j, j=1,\ldots, \nu'$ all once, then 
$F(z)$ travels the imaginary axis all just  {\it once} in $\mathbb{ C}_\tau$,
and  $\Re\,F(z)>0$ on $R$ and  $\Re\,F(z)<0$ on 
$R^*$. This fact implied that 
 $\frac{ 1}{\rho_j(\tau)} <0$ for $\tau\in \mathfrak{ C}_j \setminus \{\mathfrak{ o}_j
\}$. 
 To reach a contradiction, we focus to $C_1$. 
We may assume $\mathfrak{ o}_1=0$ of $\mathfrak{ C}_1 (\subset 
\mathbb{ C}_\tau)$ and $a_1^{(1)}=b_1^{(1)}=o_1$ on  $C_1 (\subset 
\mathbb{ C}_z)$.
If we take 
a small subarc $C_1'$ centered at $ o_1$ 
of $C_1$ and identify $C'_1$ with $I=(-r,r)$ on the $x$-axis 
such that  $ o_1$ corresponds to $0 \in I$, then 
the subarc $\Gamma: =\log 
H(C_1')$ of $\mathfrak{ C}_1$ is written
$$
\tau = \frac{ 1}2[( a_2 x^2+ a_3 x^3+ \ldots)+i
(b_2x^2+ b_3x^3+ \ldots)], \quad x\in I,
$$
where all $a_k,b_k  $ 
are real and $a_2, b_2 \ne 0$. The other cases being similar, we  
assume $a_2, b_2>0$. 
 We put $ \Gamma'(\Gamma'')= \{\log H(x) \in \Gamma: 
\mbox{$x$ travels from $0 $}$ $\mbox{to $r \, (\, -r\,)$}\}$, 
so that $\Gamma=- \Gamma''+\Gamma '$. 
Since $1/\rho_1(\tau)<0$ for $ \tau\in 
 \mathfrak{ C}_1 \setminus \{\mathfrak{ o}_1\}$, $\mathfrak{ 
 C}_1$ has a cusp singularity at $\mathfrak{ o}_1$ such that  
$\Gamma '\ (\Gamma '')$ starts  at $\mathfrak{ o}_1$ 
whose tangent decreases 
from  $b_2/a_2>0$ as  $x $ increases (decreases) from $0$ to $r\ (-r)$.
We put $\mathfrak{ a}=\log H(a_1^{(2)})$ and 
$\mathfrak{ b}=\log H(b_1^{(2)})$. Since the 
tangent $T(\tau)$ of 
 $\mathfrak{ C}_1$ at $\tau=\log H(z)$ is $T(\tau)= {v_1'(z)}/{ u_0'(z)}$, 
we have $T(\mathfrak{ a}) = 0, \ |T(\mathfrak{ b})|=\infty$  and {\it  vise versa}. 
 This contradicts that $\mathfrak{ C}_1$ is a closed curve with 
$1/\rho_1(\tau)<0$ for any $ \tau\in 
 \mathfrak{ C}_1 \setminus \{\mathfrak{ o}_1\}$, which proves i).
 
The first assertion 1. in Theorem \ref{lem:P0P1} follows v). 
Using notation (\ref{eqn:tjrj-notation}), we have for  $j=1,\ldots, \nu$, 
$${\rm Max}_{z\in C_j} \{ \Im \log H(z)\} 
 -{\rm Min}_{z\in C_j} \{ \Im \log H(z)\}  \le \frac{ 
1}2(\theta_j^{(2)}-\theta_j^{(1)}) 
< \pi,
$$
so that $-H(C_j)$  in $\mathbb{ C}_w$ as well as $-\log H(C_j)$ in 
$\mathbb{ C}_\tau$ 
is a simple closed curve which bounds a bounded domain 
in $\mathbb{ C}_w$. The second assertion in 1. is proved. 
 To prove 2., given 
 $w'\in \mathbb{ C}_w \setminus \cup_{j=1}^\nu H(C_j)$, 
we write $N(w')$ for the number of $z$ in $R$ such that  $H(z)=w'$. 
If we denote by $W_j(w')$ 
the winding number of  $H(C_j)$  about $w'$, then we  have 
$W_j(w')\le 0$ by the second assertion in 1. Since $H(z)$ has only one pole at $z=0$  
of order one on $R$, we have by the   argument principle
$$
  N(w')-1 =  \sum_ {j=1}^\nu W_j(w') \ \le 0, 
$$
so that  $N(w')=0 \ \mbox{or} \ 1$. 
 Hence, $H(z)$ is univalent on $R$, which is the first assertion 
in 2.  To prive the other ones in 2., let
 $f\in 
 {\mathcal  S}(R) $. We put 
$u(z):=\log|f(z)|$ and $h(z):= \log |H(z)|= 
 \frac{ 1}2( p(z)+q(z))$.
Then $u(z)-h(z)$ is harmonic on the whole $R$.
 Consider 
the Dirichlet 
integral of $u-h$ on $R$,
$D_R(u-h):=\iint_{ R}\bigl[ (\frac{\partial (u-h) }{\partial x})^2
+(\frac{\partial (u-h)}{\partial y })^2  \bigr]dxdy \ge 0.
$ By Green's formula we have 
\begin{align}\label{eqn:D-u-h}
\nonumber 
 D_R(u-h)&=\int_{ \partial R}udu^*-\int_{ \partial R} udh^*- \int_{ 
 \partial R} hdu^* + \int_{ \partial R} h dh^*.
\end{align}
By $\int_{ C_j}du^*=0, j=1,\ldots ,\nu$  
and condition $(L_1)((L_0))$  for $p(z)(q(z))$,
we have 
\begin{alignat*}{3}
 \int_{ \partial R} udh^* & 
= \frac{ 1}2 \int_{ \partial R} udp^* -pdu^* 
=\pi (\,\log|c(f)|- \alpha\,);\\[2mm]
 \int_{ \partial R} hdu^* 
&= \frac{ 1}2 \int_{ \partial R} qdu^* -udq^* 
=\pi (\,\beta-\log|c(f)|\,).\\
 \therefore \ \ 
 D_R(u-h)&=
\int_{ \partial R} 
udu^* + \pi(\alpha-\beta) +\int_{ 
\partial R} hdh^*.
\end{alignat*}
We put $u=h$, in particular,  to obtain 
\begin{align*}
E_{\log}(H)\! =\!-\! \int_{ \partial R}\! h dh^*= \frac{ \pi}2 (\alpha-\beta), 
 \ \ 
 E_{\log}(H) -E_{\log}(f)=D_R(u-h)\ge 0,
\end{align*}
which are desired.

{To prove 3., we calculate the harmonic span for 
the disk 
$D=\{|z|<r\}$ in $\mathbb{ C}_z$. For $\xi\in D$, we 
denote by $p(z)(q(z)) $ the $L_1$-($L_0$-)principal function 
and by $\alpha\, (\,\beta\,)$ the $L_1$-($L_0$-)constant 
for $(D, 0,\xi)$. We write $P(z)\,(\,Q(z)\,)$ the corresponding 
circular (radial) slit mapping on $D$, where $p(z)(q(z))
=\log|P(z)|\,(\,\log|Q(z)|\,)$.
We have in \S\,5 in \cite{hamano-2}
\begin{eqnarray} \label{eqn:P}
 P(z) &=& \frac{ -1}{\xi}\cdot \frac{ z-\xi}{z}\cdot
\bigl({1- \frac{ z}{r}\ \frac{ \overline{ \xi}}{r}}\ \bigr)^{-1}, 
\quad z\in D;\\
\alpha &=& \log| \frac{d P}{d z}(\xi)|=  -2 \log|\xi|
- \log \bigl[1- \bigl(\frac{ |\xi|}{r}\bigr)^2\bigr]. \nonumber 
\end{eqnarray}
Putting $\theta_\xi = \arg \xi$, we also have 
 \begin{eqnarray}\label{eqn:Q}
 Q(z) &=& 
\frac{1 }{re^{\theta_\xi}}\ \bigl[\bigl(
\frac{ z}{re^{i \theta_\xi}}+ \frac{ re^{i\theta_\xi}} {z}
\bigr)   - \bigl(
\frac{ |\xi|}{r}+\frac{ r}{|\xi|}\bigr)\bigr] \\ \nonumber 
&=& 
 \frac{ -1}{\xi}\cdot \frac{ z-\xi}{z}
\cdot \bigl({1- \frac{ z}{r}\ \frac{ \overline{ \xi}}{r}}\ \bigr), \quad z\in D;
\\
\beta &=&  \log | \frac{d Q }{dz}(\xi)|
=  -2 \log|\xi|+ \log \bigl[1- \bigl(
\frac{ |\xi|}{r}\bigr)^2\bigr].
\nonumber 
\end{eqnarray}
Hence, the harmonic span $s(D)= \alpha-\beta$ 
for $(D, 0, \xi)$ is 
\begin{eqnarray}
 \label{eqn:ex-harmonic-s}
s(D)= 2 \log \frac{ 1}{1- \bigl(\frac{ |\xi|}{r}\bigr)^2}.
\end{eqnarray}

Now let $R$ be any simply connected domain over $ \mathbb{  C}_z$ with $R\ni 0, \xi$.
We consider
the Riemann's conformal mapping   $w =\varphi (z)$ 
from $R$ onto a disk $\widetilde D:=\{|w |<r\}$ in $\mathbb{ C}_w$ such that  
$
\varphi (0)=0$ and $\varphi '(0)=1$.
We put $\Xi:= \varphi (\xi) \in \widetilde D$.
We consider the circular (radial) slit mapping 
$\widetilde P(w)\, (\widetilde Q(w))$ of $\widetilde D$ such that  
$\widetilde P(w)-\frac{ 
1}{w}\,  (\,\widetilde Q(w)- \frac{ 1}{w})$ is regular  at $0$ and $\widetilde P( \Xi)(\widetilde Q(\Xi))=0$;
 the $L_1$-($L_0$-)constant $\widetilde \alpha \,
(\,\widetilde \beta\,)$,  and 
the harmonic span $s(\widetilde D)= \widetilde \alpha -\widetilde \beta 
$ for $(\widetilde D, 0, \Xi)$. 
By (\ref{eqn:ex-harmonic-s}), we have $ s
(\widetilde D)=-2  \log [1- \bigl(
\frac{ |\Xi|}{r}\bigr)^2]$.
 Since 
$P(z):= \widetilde P(\varphi (z))$ ($Q(z):=\widetilde 
 Q(\varphi (z))$) 
becomes a circular (radial)  slit mapping 
on $R$  such that  $P(z)-\frac{ 1}{z}\ (Q(z)-\frac{ 1}{z})$ 
is regular at $0$ and $P(\xi)(Q(\xi))=0$. 
Thus, $\log |P(z)|$ ($\log 
|Q(z)|$) is the $L_1$-\,($L_0$-)principal function for $(R, 0, \xi)$, so that the $L_1$($L_0$-)constant 
$\alpha \,(\,\beta \,)$ for $(R, 0, \xi )$ 
 is 
\begin{align*}
 \alpha&=\log| \frac{d P}{dz}(\xi)|=
\log 
\bigl( 
|\frac{d \widetilde P}{d w}( \Xi)|\cdot |\frac{d \varphi }{d z}(\xi)| 
\bigr) 
=\widetilde \alpha+\log |\frac{d \varphi }{d z}(\xi)|; 
\\[1mm]
 \beta&=\log |\frac{dQ}{dz}(\xi)|=  \log 
\bigl( 
|\frac{d \widetilde  Q}{d w}(\Xi)|\cdot |\frac{d\varphi }{d z}(\xi)| \bigr)
= \widetilde \beta+\log |\frac{d \varphi }{d z}(\xi)|.
\end{align*}
Hence, the  harmonic span $s(R)=\alpha -\beta $ 
for $(R, 0, \xi)$ is
\begin{align*} 
s(R)= \widetilde \alpha-\widetilde \beta=s(\widetilde D)= 2 \log 
\frac{ 1}{1- \bigl(
\frac{ |\Xi|}{r}\bigr)^2}.
\end{align*}
Since the Poincar\'e distance 
$d(0, \xi)$ between $0$ and $\xi$  in $R$  
is equal to $
\frac{ 1}2 \log \frac{1+\frac{ |\Xi|}{r} }{1- \frac{ |\Xi|}{r}}$,
we have 
$s(R)=4 \log\,\cosh {d(0,\xi)}$, which proves 3.
\hfill $\Box$ 
}

\begin{example}
{\rm  We certify 1. and 2. in Theorem \ref{lem:P0P1}  for the case $D=\{|z|<r\}$ and $\xi\in D$. By (\ref{eqn:P}) and 
(\ref{eqn:Q}) we have 
$$
H(z)= \sqrt{P(z)Q(z)}= \frac{ 1}z - \frac { 1}{\xi}, \quad z \in D.
$$ 
Thus $H(z)$ is  univalent on $D$. Since $C:=\partial D=\{r e^{i\theta}: 
0\le \theta \le 2\pi\}$, the closed curve 
$-H( C)= \{ \frac{e^{i\theta}}r - \frac{  1}{\xi}: 0\le\theta\le 2\pi\} $
 is simple  and $-\log H(C) $  is a convex curve. Further, we have 
$E_{\log}(H)= \pi \log \frac{ 1}{1-|\xi/r|^2}$. 

In fact,  we  prove 
it in case $r=1$ and $|\xi|<1$ for simplicity. Since each branch of 
$\log\, (\frac{ 1}z- \frac{ 1}{\xi})$ is holomorphic in $\mathbb{C}_z 
\setminus D$, we have 
\begin{eqnarray*}
 E_{\log}(H)&=& \frac{i}2 \int_{-C} \log (\frac{ 1}z- \frac{ 1}{\xi}) \, 
  \,d \ \overline{ \log (\frac{ 1}z- \frac{ 1}{\xi}) }\\
&=& \frac{ -i}2 \int_{ C}\log (\frac{ 1}z- \frac{ 1}{\xi})\ \frac{dz}{
z- 1/\overline{ \xi}} \qquad \mbox{ since $z\overline{z}=1 $ on $C$}\\
&=& \frac{ -i}2 \cdot 2\pi i\bigl[( -\log (\frac{1}{1/\overline{ \xi}}-\frac{ 1}{\xi} ) + 
 \log (\frac{ -1}{\xi})\bigr] \quad \mbox{ by Cauchy theorem} \\
&=& \pi \log \frac{ 1}{1-|\xi|^2}, 
\end{eqnarray*}
which is desired. By (\ref{eqn:ex-harmonic-s}) we thus have $E_{\log}(H)={\pi s(D)}/2 $.
 }
\end{example}
 \begin{remark}\label{span-function} {\rm 
  (1)\ Let $R_i, i=1,2$ be a planar Riemann surface
 such that   $R_i \ni 0, \xi$. If we denote by $s_i$ the harmonic span 
 for $(R_i, 0, \xi)$, then we have by 2. in   Theorem \ref{lem:P0P1} 
  that $R_1 \subset R_2$ induces $s_1 \ge s_2$, even when $R_1$ and $R_2$ are not homeomorphic to each other.

(2)\ Let $R$ be a planar Riemann surface.  By the similar proof of 3.,  the harmonic span $s_R(0, \xi)$ for $(R, 0, \xi)$ 
 is invariant under the holomorphic transformations.
 Thus the harmonic span $s_R(\xi, \eta)$ for 
$(R, \xi, \eta)$ is  a $C^\omega$ positive 
function for $(\xi, \eta)\in( R 
 \times R) \setminus \cup_{\xi\in R}(\xi,\xi)$. 
It is clear that 
$s_R(\xi, \eta)=s_R(\xi,\eta)$ and, for a fixed $\xi_0\in R$, $\lim_{\eta\to \partial R} s_R(\xi_0,\eta)=+\infty$. If we put 
$s_R(\xi, \xi)=0$ for $\xi\in R$, then $s_R(\xi,\zeta)$ is $C^2$ function on $R\times R$ which satisfies, 
for a fixed $\xi_0\in R$, 
there exist $K>0$ and $\delta>0$ such that 
\begin{eqnarray} \label{eqn:xixi}\quad
\frac{1}{K} |\eta-\xi_0|^2 \le s(\xi_0, \eta)\le K|\eta-\xi_0|^2 \quad \mbox{ for $|\eta-\xi_0|<\delta$}.
\end{eqnarray}

In fact, we may assume $R$ is a bounded domain in $\mathbb{C}_z$  and $\xi_0=0\in R$. 
We take $D_a:= \{|z|<a\} \Subset  R \Subset \{|z|<b\} :=D_b $ in ${\mathbb C}_z $.
By (1) and (\ref{eqn:ex-harmonic-s}) we have, for $\eta\in D_a$,
$$
2\log \frac{1}{1- |\eta/b|^2}= s_{D_b}(0, \eta) \le s_{R}(0,\eta) \le s_{D_a}(0,\eta)= 2\log \frac{1}{1- |\eta/a|^2},
$$
which implies (\ref{eqn:xixi}). 

  We call the function $s_R(\xi,\eta)$ on $R\times R$
 the {\it $S$-function for $R$}.}
\end{remark}
\section{ Variation formulas for the harmonic spans}
 
\noindent We return to the variation of Riemann surfaces. In this section, as in section 2.,
we assume that $\widetilde {\mathcal  R}=\cup_{t\in B}(t, \widetilde R(t))$ is  an 
unramified domain over $B \times \mathbb{ C}_z$ and  ${\mathcal  
R}=\cup _{t\in B}(t, R(t))$ satisfies conditions {\bf  1.} and {\bf  2.} in the 
beginning of  section 2.
 For a fixed $t\in B$, let  $p(t,z)\ (q(t,z))$;
 $\alpha(t)\ (\beta (t)) $ and  $s(t)$ denote 
 the 
$L_1$-($L_0$-)principal function;  the $L_1$-($L_0$-)constant
 and  the harmonic span,  
for $(R(t), 0, \xi(t))$.
 Then Lemmas 
\ref{vari-form:u} and  \ref{vari-form:L_0} imply the following variation formulas:
\begin{lemma} \label{lem:vari-formula-h-span}
\begin{align*}
 \frac{\partial s(t) }{\partial t}
&= \frac{1}{\pi} \int_{\partial R(t)} k_1 (t,z) \ \bigl(\
\bigl|
\frac{\partial p(t,z)}{\partial z} 
\bigr|^2+ \bigl|
\frac{\partial q(t,z)}{\partial z} 
\bigr|^2\ \bigr) \, ds_z;\\
 \frac{\partial^2 s(t) }{\partial t \partial \overline{ t} }
&= \frac{1}{\pi} \int_{\partial R(t)} k_2 (t,z) \ \bigl(\
\bigl|
\frac{\partial p(t,z)}{\partial z} 
\bigr|^2+ \bigl|
\frac{\partial q(t,z)}{\partial z} 
\bigr|^2\ \bigr) \, ds_z \nonumber \\[2mm]
&\quad + \frac{ 4}\pi \iint_ {R(t)}\bigl(\
\bigl|
\frac{\partial^2 p(t,z)}{\partial\overline{t} \partial{z}}
\bigr|^2 +
\bigl|
\frac{\partial^2 q(t,z)}{\partial\overline{t} \partial{z}}
\bigr|^2\ \bigr) \,dxdy\\
&\quad +\, \frac{ 2}{\pi}\, \Im\, { \,\sum_ {k=1}^g
 \bigl(\frac{\partial }{\partial t} 
\int_{ A_k(t)}\!\! * dq(t,z)\bigr)\cdot\bigl(
\frac{\partial }{\partial \overline{ t}}
 \int_{ B_k(t)}\!\! * dq(t,z) \bigr)
 }. 
\end{align*}
\end{lemma}

We say, in general, that ${\mathcal  R}: t\in B \to R(t)$   is  {\it equivalent to a 
trivial variation}, if there exists a 
biholmorphic transformation $T$ from the total 
space ${\mathcal  R}$ onto 
 a product space $B \times D$ (where $D$ is a Riemann surface) of the form 
$T:(t,z)\in {\mathcal  
R}\mapsto (t,w)=(t, f(t,z)) \in B \times D$.   

 In case $R(t)$ is planar, following 
(\ref{eqn:to-mapping}), on each $R(t), t\in B$ we construct 
 the circular and radial 
slit mappings:
\begin{align*}
P(t,z)&= 
e^{p(t,z)+ip(t,z)^*} \quad \mbox{ and } \quad Q(t,z)= 
e^{q(t,z)+iq(t,z)^*}
\end{align*}
such that  $P(t,z)- \frac{ 1}z$ and $Q(t,z)-\frac{ 1}z$ are 
regular  at $z=0$. 
We put $D_1(t)=P(t, R(t))$ and $D_0(t)=Q(t,R(t))$, so 
that 
\begin{align*}
 D_1(t)&=
\mathbb{P}_{w} \setminus 
\cup_{j=1}^\nu P(t,C_j(t))=
 \mathbb{P}_w \setminus 
\cup_{j=1}^\nu \ \mbox{{arc}$\{A_j^{(1)}(t),A_j^{(2)}(t)\}$}
;\\[2mm] D_0(t)&= \mathbb P_w \setminus \cup_{j=1}Q(t,C_j(t))=
\mathbb{P}_w\setminus \cup_{j=1}^\nu \mbox{{ segment}$\{
B_j^{(1)}(t), B_j^{(2)}(t)\}$}.
\end{align*}

\begin{theorem}
 \label{cor:rigidity}
 Assume that ${\mathcal R}=\cup_{t\in B}(t, R(t))$ is 
pseudoconvex in $\widetilde {\mathcal  R}$ and each  $R(t), \ t\in B$ 
is planar.
 Then 

 1. \ $s(t)$ is $C^\omega $ subharmonic on $B$; \\
\ \ \  2. \ if  $s(t)$ is harmonic on $B$, then 
\begin{enumerate}
 \item [(i)] \  $s(t)$ is constant on $B$;
 \item [(ii)] \ ${\mathcal  
       R}:t\in B \to R(t)$ is equivalent to a trivial variation. More concretely,
\begin{enumerate}
 \item [$(\diamond)$]  \   ${\mathcal  R}$ is biholomorphic to the  product 
       domain $B \times \widetilde D_1$, 
 where $\widetilde D_1$ is a circular 
       slit domain in $\mathbb{ P}_w$ such that  
$\widetilde D_1=\mathbb{ P}_w \setminus \cup _{j=1}^\nu 
\{\widetilde 
       A_j e^{i \theta}: 0\le \theta \le \Theta_j\}$, 
where $ \widetilde A_1=1$ and   each $\widetilde A_j (\ne 0), j= 2,\ldots ,\nu$ 
is constant,  by the holomorphic transformation $T_0: \ (t,z)\in {\mathcal  R} 
       \mapsto (t,w)=(t,\widetilde P(t,z))\in B \times \widetilde D_1$, 
where $\widetilde P(t,z)=P(t,z)/A_1^{(1)}(t) $.
\end{enumerate}
\end{enumerate}
\end{theorem}

The concrete $(\diamond)$ will be used in the proof of Corollary 
\ref{cor:rigidity-2}.

\vspace{1mm} 
\noindent {\it  Proof.}  Lemma \ref{lem:vari-formula-h-span} implies 1.
 To prove 2., we may  assume that 
${\mathcal  R}=\cup_{t\in B}(t, R(t))$ is an unramified 
domain over $B \times \mathbb{ C}_z$ such that  each $R(t), \ t\in B$
 is contained in an unramified planar domain $\widetilde R$ over 
$\mathbb{ C}_z$ and the holomorphic section $\xi $ is constant: $t\in B 
\to \xi (t)= 1\in R(t)$.
  Assume that $s(t)$ is harmonic on $B$.  
 By Lemma 
 \ref{lem:vari-formula-h-span}, we have 
\begin{enumerate}
 \item [a)]\  $k_2(t,z) \equiv 0$ on $\partial 
 {\mathcal  R}$, i.e., $\partial {\mathcal  R}$ is a  Levi flat surface over $B \times 
 \mathbb{ C}_z$;
 \item [b)] \ both $\frac{\partial p(t,z) }{\partial z}$ and $\frac{\partial  q(t,z)}{\partial z}$ are 
holomorphic for $t\in B$.
\end{enumerate}
By b) and the normalization at $z=0$,  both $w=P(t,z)$ and $w=Q(t,z)$ are holomorphic 
for two complex variables $(t,z)$ in ${\mathcal  R}$ except 
$B \times \{0\}$. We put $D_1(t)=P(t,R(t)) \subset \mathbb{ P}_w$ for 
$t\in B$, and ${\mathcal  D}_1=\cup_{t\in B}(t, D_1(t)).$ Since ${\mathcal  D}_1$ as well as ${\mathcal  R}$ over $B \times \mathbb{ C}_z$ 
is a pseudoconvex (univalent) 
domain in $B \times \mathbb{ P}_w$, it follows from Kanten Satz (p.\,352 in \cite{behnke}) that  
each edge point $A_j^{(k)}(t)$ is holomorphic for $t\in B$, 
and $A_j^{(2)}(t)= A_j^{(1)}(t)e ^{i \Theta_j}$, where $\Theta_j$ 
is constant for $t\in B$. 
We consider the map $(t,w)\in {\mathcal  D}_1 \mapsto 
(t, \widetilde w)=(t, L(t,w))\in B \times \mathbb{ P}_{\widetilde w}$, 
where $L(t,w)= w/A^{(1)}_1
(t)$, and put 
$\widetilde {\mathcal  D}_1=\cup_{t\in B}(t, \widetilde D_1(t))$ where $\widetilde D_1(t)=
L(t, D_1(t))$. Each $\widetilde D_1(t),\ t\in B$ is circular slit domain 
$\mathbb{ P}_{\widetilde w} \setminus \cup_{j=1}^\nu \widetilde C_j(t)$ 
such that the first circular slit 
$\widetilde C_1(t)
= \{e^{i \theta}:  0\le \theta\le \Theta_1\} $ is
independent of $t\in B$, say $\widetilde C_1:=\widetilde 
C_1(t)$. 
Since ${\mathcal  R}$
 is biholomorphic to $\widetilde {\mathcal D}_1 $, 
 and each $\widetilde D_1(t), t\in B$ has no ramification points, it 
 suffices for $(\diamond)$ in 2.\,(ii) to prove 
that
{ the edge point 
$\widetilde A_j^{(1)}(t):= A_j^{(1)}(t)/ A_1^{(1)}(t)$
of each arc  $\widetilde C_j^{(1)}(t), 
 j=2,\ldots , \nu$  does not depend on $t\in B$.}
 
In fact, we see from b) that the function $F(t,z)$ defined 
in (\ref{eqn:gur-function}): 
$$
 W=F(t,z)= \frac{ d_z\log Q(t,z)}{d_z\log P(t,z)},   \quad 
\mbox{ $z\in R(t)\cup \partial R(t)$}
$$
is holomorphic for  
$t\in B$ such that $F(t,0)=1$ and $\Re  F(t, z)=0$ on  $\partial 
R(t)$, i.e., 
$F(t,z)$ is meromorphic function for 
two complex variables $(t,z) \in {\mathcal  R}$ such that  
$\Re\,F(t,z)=0 $ on $\partial {\mathcal  R}$. We put 
$K_j(t)= F(t, C_j(t)), j=1,\ldots , \nu$ in $\mathbb{ P}_W$. In the 
1st step of the proof of 1. in Theorem \ref{lem:P0P1}
we proved that $K_j(t) $ rounds just twice on the imaginary axis in $\mathbb{ P}_W$.
We put $W(t)=F(t, R(t))$ and ${\mathcal  W}=
\cup _{t\in B}(t, W(t))$, so that $\partial {\mathcal  W}
=\cup_{t\in B}(t, \cup _{j=1}^\nu K_j(t))$, and ${\mathcal  R}\approx {\mathcal  W}$ 
(biholomorphic)
 by $T: (t,z)\in {\mathcal  R} \mapsto (t,W)=(t,F(t,z))\in {\mathcal  W}$.
Thus, $W(t) $  has $2\nu +g-1$ ramification points.
 Consider the following biholomorphic mapping 
$
(t,W) \in {\mathcal  W} \to (t,\widetilde w)=(t, \widetilde G(t,W))\in 
 \widetilde {\mathcal  D}_1,
$
where $\widetilde G(t,W):=L(t, P(t, F^{-1}(t,W)))$. 
We use  the following elementary fact:

\vspace{1mm}
 
$(\ast)$ \  {\it  Let $B=\{|t|<\rho\}$ in $\mathbb{ C}_t$ and 
$E=\{|z|< r\}\cap \{\Re\, z\ge 0\}$ in $\mathbb{ C}_z$. If $f(t,z)$
is a holomorphic function for two complex variables $(t, z)$ on $B \times E$ such that  
$|f(t,z)|=1$ on $B \times (E \cap \{\Re\, z=0\})$, then 
$f(t,z)$ does not depend on $t\in B$.}

\vspace{1mm}
We choose a point $W_0$ on $\partial K_1(0) \subset \partial {\mathcal  W}$
 such that  
$\widetilde G(0,W_0)= e^{i\theta_0} \in \widetilde C_1$ 
with $0<\theta_0< \Theta_1$ and the direction of $\widetilde C_1$ at 
$e^{i \theta_0}$ follows as $\theta_0$ increases. Then  we have  a small disk $B_0\subset B$ of center  $0$ and 
 a small half-disk $E=\{|W-W_0|<r\} \cap \{\Re\, W\ge 0\}$ 
in $\mathbb{ C}_W$ 
such that $|\widetilde G(t, W)|\le 1\ (=1)$ on $B_0 \times E \ (B_0 
\times (E \cap \{\Re\, W=0\}))$. By $(\ast)$,  $\widetilde G(t,W)$ 
for $W\in E \cap\{\Re\, W\ge 0\}$ 
does not depend on $t\in B_0$. By 
the  analytic continuation,  $\widetilde G(t,W)$ on ${\mathcal  W}\cup 
\partial {\mathcal  W}$ 
does not depend on $t\in B$. 

Now assume that some  $\widetilde A^{(1)}_j(t)$, $ 2\le \exists \ j \le 
\nu$ is not constant for $t\in B$.
We take a point $W_0\in \mathbb{ C}_W$ with $\Re\, W_0=0$.
Since the component $K_j(t)$ of $\partial { W}(t)$ 
winds  twice around the imaginary axis in 
$\mathbb{ P}_W$, for each $t\in B$ we find 4 points of $K_j(t)$ over $W_0$.
We fix one of them, say $W_0(t) \in K_j(t)$, to whom the corresponding point 
$z_j(t) \in C_j(t)$ continuously varies in $\partial {\mathcal  R}$ 
with $t\in B$. 
 Since  $\widetilde C_j(t)= 
 \widetilde G(t, K_j(t))= \{ \widetilde A_j^{(1)}(t) e ^{i \theta}: 0\le 
 \theta \le \Theta_j\}$, where $\Theta_j$ is constant for $t\in B$,
 we have 
 $\widetilde G(t, W_0)= \widetilde A_j^{(1)}(t) e^{i \theta(t)}$, 
  where $\theta (t) \ (0<\theta(t) < \Theta_j)$, continuously varies with 
 $t\in B$. Since $|\widetilde A_j^{(1)}(t)|$ 
as well as $\widetilde A_j^{(1)}(t)$ is not constant for $t\in B$, 
$\widetilde G(t, W_0)$ {\it does} depend on $t\in B$,  a contradiction, and  2.\,(ii) is proved.

From (2) in  Remark \ref{span-function} 
 the harmonic span $s(t)$ for $(R(t), 0,1)$ 
is equal to that for $(\widetilde D_1(t), \infty,0)$. Since $\widetilde 
D_1(t)= \widetilde D_1(0)$ for any $t\in B$, $s(t)$ is constant on $B$, 
which proves 2.\,(i).
\hfill $\Box$ 

\vspace{1mm}
For 2.\,(ii) in Theorem \ref{cor:rigidity} we cannot replace the 
condition of the harmonicity of $s(t)$ on $B$ by that of $\alpha (t)$ or 
$\beta (t)$ on $B$, in general.  However, when $R(t), t\in B$ is simply connected, such
 replacement is 
 possible by the same idea of the proof of 2.\,(ii).  
\begin{corollary} \label{cor:szz} Assume that ${\mathcal R}=\cup_{t\in B}(t, R(t))$ is 
pseudoconvex over $B \times \mathbb{ C}_z$ and each  $R(t), \ t\in B$ 
is planar.  Then  the  $S$-function $s(t,\xi,\eta)$ for  $R(t), t\in B$  is 
$C^2$ plurisubharmonic on ${\mathcal  R}^2:=\cup_{t\in B}(t, R(t)\times R(t))$.
In particular, for a fixed $t_0\in B$,  we simply put $R(t_0)=R$ and $s(t_0,\xi,\eta)$ $=s(\xi,\eta)$. 
Then $s(\xi,\eta)$ is $C^2$  plurisubharmonic  on $R \times R$ such that, for  
any complex line $l $ except $\xi=\eta$ in $R \times R$, the restriction of $s(\xi, \eta) $ 
on $l\cap (R\times R)$ is strictly subharmonic. 
\end{corollary}

\noindent {\it Proof.} \  Let 
 $t\in B \to (\xi(t),\eta (t))\in R(t) 
\times R(t)$ be any holomorphic mapping from $B$ into ${\mathcal R}^2$. 
We put   $s(t):=s(t,\xi(t), \eta(t))$ for $t\in B$, and 
$B'=B \setminus \{t\in B:\xi(t)=\eta(t)\}$.
 Consider the translation $T: \ 
(t,z)\in {\mathcal  R} \mapsto (t,w)=(t, 
z-\eta(t))$ for $t\in B'$, and  put $\widetilde {\mathcal  R}:
=T({\mathcal  R})$ and $\widetilde \xi=T\xi$.
Then $\widetilde {\mathcal  R}$ is 
 pseudoconvex   over 
$B' \times \mathbb{C}_w$  and $\widetilde \xi \in \Gamma (B', \widetilde 
{\mathcal  R})$.  By Theorem 
\ref{cor:rigidity},  the harmonic span 
$\widetilde s (t)$ for $(\widetilde R(t), 0, \widetilde \xi (t))$ 
is $C^\omega$ subharmonic on $B'$, and  so is $s(t)$ on $B'$. It follows from 
 (\ref{eqn:xixi}) that 
$s(t)$ is $C^2$ subharmonic on $B$, which proves the former part in the corollary. 
Further, by the same argument we can prove the latter part  under the second variation formula in Lemma \ref{lem:vari-formula-h-span} and  (\ref{eqn:xixi}). \hfill $\Box$

\vspace{1mm}
Theorem \ref{cor:rigidity} with 3. in Theorem \ref{lem:P0P1} directly implies
\begin{corollary} \label{cor:mdp}
Assume that ${\mathcal R}=\cup_{t\in B}(t, R(t))$ is 
pseudoconvex over $B \times \mathbb{ C}_z$ and   $R(t), \ t\in B$ 
is simply connected. Let  $\xi_i \in \Gamma(B, {\mathcal R}), \ i=1,2$ and 
let $d(t)$ denote the Poincar\'e distance between $\xi_1(t)$ 
and  $\xi_2(t)$ on $R(t)$. Then $\delta(t):=\log\,\cosh d(t)$ 
is subharmonic on $B$. Moreover, $\delta (t)$ is harmonic on $B$ if and only if $\mathcal{ R}$ is equivalent to the trivial variation.
\end{corollary}
 Prof. M.\,Brunella said to us that he could prove  the stronger 
 fact:\ \  {\it "\,$\log d (t)$ is subharmonic on $B$''} than   {\it "$\delta  (t)$ is subharmonic on $B$"}
 by the same idea  in p.\,139 in 
 \cite{brunella} which is based on \cite{b-g}, (though there was not its exact statement).
\begin{remark}\ 
{\rm  (1)  In \S 2 and \S 3,  ${\mathcal  R}=\cup_{t\in 
 B}(t, R(t))$ is assumed to be  a subdomain of an unramified domain 
 $\widetilde {\mathcal  R}=\cup_{t\in B}(t, \widetilde R(t))$ over $B \times \mathbb{ 
 C}_z$ which satisfies conditions {\bf  1.} and {\bf  2.} stated in \S 2. By the standard use of the immersion 
theorem for open  Riemann surfaces in 
\cite{g-n} or 
\cite{nishimura},
 the results in \S 2 and \S 3 hold for the following ${\mathcal  R}$: 
{let $B=\{t\in \mathbb{ C}: |t|<\rho\}$ and let $\pi: \widetilde {\mathcal  R} \to B$ be 
a two-dimensional holomorphic family (namely, $\widetilde {\mathcal  R}$  
 is a complex two-dimensional manifold and $\pi$ is a holomorphic 
 projection from $\widetilde {\mathcal  R}$ onto $B$)  
such that  each fiber $\widetilde R(t)=\pi^{-1}(t),\, t\in B$ is 
irreducible and 
 non-singular in $\widetilde {\mathcal  R}$. Putting 
$\widetilde {\mathcal  R}= \cup_{t\in 
 B}(t, \widetilde R(t))$, our ${\mathcal  R}$ is  a subdomain 
of $\widetilde {\mathcal R} $ defined by 
  ${\mathcal  R}=\cup_{t\in B}(t, R(t)) \subset 
 \widetilde {\mathcal  R}$ 
which satisfies the corresponding conditions 
 ${\bf  1.}$ and ${\bf  2.}$

(2) In conditions {\bf  1.} and {\bf  2.},  if we replace 
 $C^\omega$ smooth by $C^\infty$ smooth, i.e., 
 ${\mathcal  R}: t\in B \to R(t) \Subset \widetilde R(t)$ is a variation such that  
 $\partial  R(t), 
 t\in B$ is $C^\infty$ smooth in $\widetilde R(t)$ 
and $\partial {\mathcal  R}$ is $C^\infty$ 
 smooth in $\widetilde {\mathcal  R}$, then the results in \S 2 and 
 \S 3  hold by replacing $C^\omega$ by $C^\infty$.
In fact, Lemmas 
\ref{vari-form:u} and \ref{vari-form:L_0}, 
on which all results are based,  hold 
for the  $C^\infty$ category  by a little not essentially 
 change of  the proofs for the $C^\omega$ category (cf: \S 2 in \cite {l-y}).}}
\end{remark}

 \vspace{-5mm}
\section{Approximation theorem for general variations of  planar Riemann surfaces}
\noindent We 
consider 
the general variation of  Riemann surfaces 
$
{\mathcal  R}: t\in \Delta \to R(t)
$
defined as follows: let   $\Delta$  be an open or a compact Riemann 
surface and $\pi: {\mathcal  R} \to \Delta $ 
be a two-dimensional holomorphic family such that  
each fiber $R(t)=\pi^{-1}(t), \ t\in \Delta $ is irreducible 
and non-singular in ${\mathcal  R}$ and is {\it planar}.
In case $\Delta $ is open, we assume that ${\mathcal  R}$ is Stein.
 We call such ${\mathcal  R}$ the variation {\it  of type} ({\bf A}).
In case $\Delta $ is compact, we assume that,  
for any disk $B \subset \Delta $, ${\mathcal  R}|_B$ is of type $({\bf  
A})$, i.e., $\pi^{-1}(B)=\cup_{t\in B}(t,R(t))$ is Stein.
We call such ${\mathcal  R}$ the variation {\it  of type}  
({\bf B}).  In general, $R(t)$ might be infinite ideal boundary components 
and ${\mathcal R}: t\in \Delta \to R(t)$ might not be topologically trivial.  
 To state the approximation theorem for these variations ${\mathcal  R} $ 
we make the following

\vspace{2mm}

\noindent {\bf Preparation.}   Let $\Delta $ 
and $\pi: {\mathcal  R}\to \Delta $ 
be of type ({\bf  A}). 
 Due to Oka-Grauert  (cf: Theorem 8.22 in \cite{nishino}),  ${\mathcal R}$ admits  a $C^\omega$ strtictly plurisubharmonic 
 exhaustion function $\psi(t,z)$.
 Let $\xi, \ \eta \in \Gamma (\Delta, {\mathcal  R})$ 
such that  $\xi \cap \eta=\emptyset$. 
 Let $B \Subset \Delta $ be a small disk such that  
we find a continuous  curve $g(t)$ connecting $\xi (t)$ and 
$\eta(t)$ on $R(t), \,t\in B$  which continuously varies in ${\mathcal  R}$ 
with $t\in B$.
 We put ${\mathcal  
  R}|_{B} 
 =\cup _{t\in B}(t, R(t))$; $\xi |_B=\cup_{t\in B}(t, \xi (t)); \ 
 \eta|_B=\cup_{t\in B}(t, \eta(t))$, and  $g|_B= \cup _{t\in B}(t, 
g(t))$. We take  so large $a \gg 1$ that
 ${\mathcal  R}(a)|_B:= \{ (t,z) 
\in {\mathcal  R}|_B: \psi (t, z)<a\} \supset g|_B$.  
Then we find an increasing 
sequence $\{a_n\}_n$ such that $\lim_{ n \to 
\infty}a_n=\infty$ and if we put 
\begin{equation}
\label{eqn:Rn}  {\mathcal  R}_n = \mbox{the conn. comp. of ${\mathcal  
 R(a_n)|_B}$} \mbox{ which contains $g|_B$,}
\end{equation}
  then  1) each ${\mathcal  R}_n, n=1,2,\ldots $ is a connected domain 
with real three-dimensional $C^\omega$ 
 surfaces $\partial {\mathcal  R}_n$ in ${\mathcal  R}|_B$ (but each 
       $R_n(t), t\in B$ is not always connected);

2) if we consider 
       the set ${\mathcal  L}$ of points $t\in B$ such that there exists a point $(t,z(t)) 
       \in \partial {\mathcal  R}_n$ with $\frac{\partial \psi}{\partial 
       z}(t, z(t))=0$, then  ${\mathcal  L}$ consists of two kind 
       of  families ${\mathcal  L}', \ {\mathcal  L}''$ of 
finite $C^\omega$ arcs in $B$:  

\vspace{2mm}
\qquad $${\mathcal  L}'=\{l_1', \ldots , l'_{m}\}, \qquad {\mathcal  L}''=
\{l_1'', \ldots , l''_{\mu}\}, $$
which have  the following property:

\underline {for ${\mathcal  L}'$}: \ 
for any $t_0\in {\mathcal  L}'$ (except a finite set 
at which some $l'_i$ and $l'_j$ or $l'_i$ itself 
       intersect transversally), 
 say $t_0\in l_i'$, $\partial R_n(t_0)$ (consisting of a finite number of closed curves) 
has only one singular point at $z(t_0)$, 
 and we find a small bi-disk $B_0 \times V$
centered at $(t_0, z(t_0))$ in ${\mathcal  R}_{n+1} $ 
such that  $B_0 \Subset B $
 and  $l'_i\cap B_0$ divides $B_0$ into two connected 
       domains $B_0' $ and $B_0''$ in the manner that 
\begin{enumerate}
\item [i)] each $\partial R_n(t), t\in B_0'\cup B_0''$ has no singular points;
 \item [ii)] each $\partial R_n(t), t\in l'_i\cap B_0 $ has one 
singular point $ z(t) $ at which two subarcs of $\partial 
      R_n(t)$ transversally intersect;
 \item [iii)] each $R_n(t)\cap V, t\in B'_0\cup (l_i'\cap B_0)$ 
      consists of two (connected) domains, while each $R_n(t)\cap V, t\in 
      B_0''$ consists of one domain;  
\end{enumerate}

\underline {for ${\mathcal  L}''$:}\ for any $t_0\in {\mathcal  L}''$ 
(except a finite point set), say 
$t_0\in l''_i$, we find a unique point $(t_0,z(t_0)) \in \partial 
       {\mathcal  R}_n$ with 
$\frac{\partial \psi}{\partial z}(t_0, z(t_0))=0 $,
and a small bi-disk $B_0 \times V$ centered at 
$(t_0,z(t_0))$ in ${\mathcal  R}_{n+1}$ such that $B_0 \Subset B$ and  
 $l''_i \cap B_0$  divides $B_0$ into two connected 
       domains $B_0' $ and $B_0''$; a $C^\omega$ mapping ${\mathfrak z}$: 
$t\in l''_i \cap B_0 \to z(t)$ such that  $(t, z(t))\in \partial 
       {\mathcal  R}_n$ with $\frac{\partial \psi}{\partial z}(t, z(t))=0 $
 in the manner that 
\begin{enumerate}
\item [i)]  $[R_n(t)\cup \partial R_n(t)] \cap V= \emptyset $ for $t\in B_0'\cup ( l''_i\cap B_0)$;
 \item [ii)] $ R_n (t) \cap V $ for $t\in B_0''$ is a simply connected 
      domain $\delta_n(t)$
      such that, for a given $t^0\in  l''_i\cap B_0$, $\delta _n(t)$ 
      shrinkingly  approaches  the point $ z(t^0)$ as $t\in B_0'' \to t^0$.  
\end{enumerate}

For the singular point $z(t), t\in l'_i \subset {\mathcal  L}'$, 
we have the connected component $C(t)$ of 
$ \partial R_n(t)$ passing through $z(t)$. Then $C(t)$ consists of one 
closed curve,   or two 
closed curves $C_i(t), i=1,2$ such that $C(t)=C_1(t)\cup C_2(t)$ 
and $C_1(t)\cap C_2(t)=z(t)$. 
 For example, in figure (FIII) below, $C(t)$ consists of one closed 
curve, and in figures (FI), (FII), $C(t)$ consists of two closed curves.

For the singular point $z(t), t\in l''_i  \subset {\mathcal  L}''$, 
  we have 
$(t, z(t))\in \partial 
{\mathcal  R}_n$ but $z(t) \not\in \partial R_n(t)$.

 Fix  $t\in B$ and $n \ge 1$ and consider 
the connected component $R_n'(t)$ of $R_n(t)$ which contains $g(t)$.
We put ${\mathcal  
R}_n'=\cup_{t\in B}(t, R_n'(t))$ and $\partial {\mathcal  R}_n'= \cup 
_{t\in B} (t, \partial R'_n(t))$. The variation
$$
{\mathcal  R}_n': 
t\in B \to R_n'(t)
$$ is no longer smooth variation of $R'_n(t)$ with 
$t\in B$, i.e., ${\mathcal  R}_n'$ satisfies  neither  corresponding 
condition {\bf  1.} nor {\bf 2.} of ${\mathcal  R}$ in \S 2. Since $R(t)$ is irreducible in $\mathcal R$, we have 
$R_n'(t) \Subset R'_{n+1}(t)$; ${\mathcal R}'_n\to 
{\mathcal R}|_B \ (n\to \infty)$, and $R_n'(t) \to R(t) \ (n\to \infty)$ for $t\in B$. 

By i), ii) for ${\mathcal  L}''$, there exists 
a neighborhood ${\mathcal V}$ of $\cup_{t\in {\mathcal  
L}'' }(t,  z(t))$ in ${\mathcal  R}_{n+1}$ such that  
 $[{\mathcal  R}'_n \cup 
\partial {\mathcal R}'_n ]\cap {\mathcal  V}= \emptyset$, so that ${\mathcal  L}''$ does not give any 
 influence for the variation ${\mathcal  R}_n'$ (contrary to for 
 ${\mathcal  R}_n$). 

Each $R(t), t\in \Delta $
is assumed \, {\it planar}. We separate 
the singular point $z(t)$ of $\partial R_n(t), t\in l'_i 
 \subset {\mathcal  L}'$ such that  $z(t)\in \partial R_n'(t)$ 
into the following two cases $({\bf c1})$ and $({\bf  c2})$: let
$C(t)$ denote the connected component 
of $\partial R_n(t)$ passing through $z(t)$. Then 
\begin{enumerate}
 \item [({\bf  c1})]  $C(t)$ consists of two closed curves 
       $C_i(t),i=1,2$, and  one of them, say $C_1(t)$, is one of 
        boundary 
       components of $R_n'(t)$, so that $[C_2(t) \setminus \{z(t)\}] \cap 
\partial R_n'(t)=\emptyset$; 
 \item [({\bf  c2})]  $C(t)$ is one of the boundary components of 
 $R_n'(t)$. 
\end{enumerate}

For example, if we take the shadowed part 
in figure (FI) (resp. (FII), (FIII)) as $R'_n(t)$, then 
the singular point $z(t)$ 
is of case 
({\bf c1}) (resp. ({\bf c2})).

\vspace{3mm}

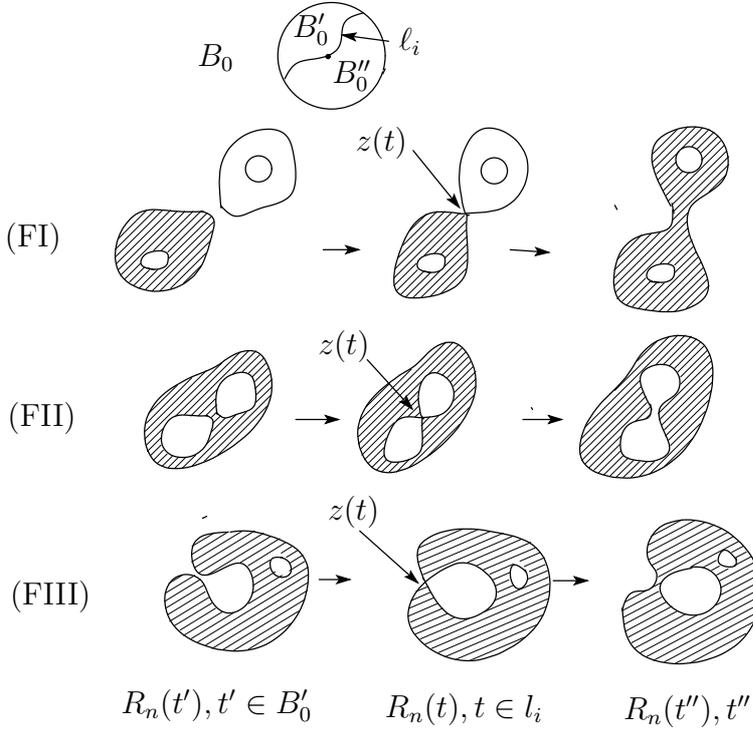
\unitlength 0.1in
\begin{picture}( 11.1900,  5.5600)( 8.9000,-15.6500)
%
\special{pn 8}%
\special{ar 2614 1288 278 278  0.3063240 6.2831853}%
\special{ar 2614 1288 278 278  0.0000000 0.2662520}%
%
\special{pn 8}%
\special{sh 1}%
\special{ar 2596 1292 10 10 0  6.28318530717959E+0000}%
\special{sh 1}%
\special{ar 2596 1292 10 10 0  6.28318530717959E+0000}%
%
\special{pn 8}%
\special{pa 2778 1064}%
\special{pa 2744 1074}%
\special{pa 2712 1086}%
\special{pa 2688 1104}%
\special{pa 2674 1128}%
\special{pa 2668 1160}%
\special{pa 2666 1194}%
\special{pa 2660 1226}%
\special{pa 2646 1254}%
\special{pa 2622 1274}%
\special{pa 2592 1290}%
\special{pa 2560 1298}%
\special{pa 2530 1304}%
\special{pa 2498 1308}%
\special{pa 2466 1314}%
\special{pa 2436 1324}%
\special{pa 2412 1344}%
\special{pa 2392 1370}%
\special{pa 2376 1400}%
\special{pa 2372 1408}%
\special{sp}%
\put(25.0900,-11.4800){\makebox(0,0){$B'_0$}}%
\put(27.2000,-13.9000){\makebox(0,0){$B''_0$}}%
\put(30.2100,-12.1500){\makebox(0,0){$\ell_i$}}%
%
\special{pn 8}%
\special{pa 2910 1212}%
\special{pa 2672 1180}%
\special{fp}%
\special{sh 1}%
\special{pa 2672 1180}%
\special{pa 2736 1208}%
\special{pa 2724 1188}%
\special{pa 2740 1170}%
\special{pa 2672 1180}%
\special{fp}%
%
\special{pn 8}%
\special{pa 2892 1356}%
\special{pa 2880 1382}%
\special{fp}%
\put(20.1500,-12.9800){\makebox(0,0){$B_0$}}%
\end{picture}%

\vspace{2mm} 

\unitlength 0.1in
\begin{picture}( 37.8000,  9.6600)(  5.6300,-17.1600)
%
\special{pn 8}%
\special{pa 1622 1192}%
\special{pa 1590 1194}%
\special{pa 1558 1194}%
\special{pa 1526 1194}%
\special{pa 1494 1192}%
\special{pa 1460 1190}%
\special{pa 1428 1186}%
\special{pa 1394 1184}%
\special{pa 1362 1182}%
\special{pa 1330 1186}%
\special{pa 1300 1194}%
\special{pa 1272 1208}%
\special{pa 1246 1230}%
\special{pa 1224 1256}%
\special{pa 1204 1284}%
\special{pa 1188 1318}%
\special{pa 1176 1352}%
\special{pa 1170 1388}%
\special{pa 1166 1424}%
\special{pa 1170 1458}%
\special{pa 1178 1490}%
\special{pa 1192 1520}%
\special{pa 1212 1544}%
\special{pa 1236 1566}%
\special{pa 1266 1582}%
\special{pa 1298 1596}%
\special{pa 1330 1604}%
\special{pa 1366 1610}%
\special{pa 1400 1608}%
\special{pa 1434 1604}%
\special{pa 1466 1596}%
\special{pa 1496 1584}%
\special{pa 1526 1568}%
\special{pa 1552 1548}%
\special{pa 1576 1528}%
\special{pa 1594 1504}%
\special{pa 1612 1480}%
\special{pa 1626 1452}%
\special{pa 1636 1424}%
\special{pa 1648 1392}%
\special{pa 1656 1360}%
\special{pa 1666 1326}%
\special{pa 1674 1292}%
\special{pa 1682 1256}%
\special{pa 1676 1226}%
\special{pa 1654 1206}%
\special{pa 1624 1194}%
\special{pa 1614 1192}%
\special{sp}%
%
\special{pn 8}%
\special{pa 2954 1204}%
\special{pa 2920 1198}%
\special{pa 2886 1196}%
\special{pa 2854 1194}%
\special{pa 2822 1196}%
\special{pa 2794 1202}%
\special{pa 2768 1214}%
\special{pa 2744 1230}%
\special{pa 2722 1250}%
\special{pa 2702 1274}%
\special{pa 2684 1304}%
\special{pa 2668 1334}%
\special{pa 2654 1368}%
\special{pa 2640 1404}%
\special{pa 2628 1442}%
\special{pa 2618 1478}%
\special{pa 2612 1516}%
\special{pa 2610 1550}%
\special{pa 2612 1580}%
\special{pa 2622 1606}%
\special{pa 2638 1624}%
\special{pa 2662 1638}%
\special{pa 2692 1644}%
\special{pa 2726 1644}%
\special{pa 2762 1642}%
\special{pa 2800 1636}%
\special{pa 2836 1628}%
\special{pa 2870 1616}%
\special{pa 2900 1602}%
\special{pa 2928 1586}%
\special{pa 2950 1566}%
\special{pa 2968 1544}%
\special{pa 2982 1518}%
\special{pa 2988 1488}%
\special{pa 2992 1456}%
\special{pa 2992 1424}%
\special{pa 2992 1390}%
\special{pa 2990 1354}%
\special{pa 2988 1320}%
\special{pa 2986 1288}%
\special{pa 2986 1256}%
\special{pa 2988 1224}%
\special{pa 2988 1212}%
\special{sp}%
%
\special{pn 8}%
\special{pa 2954 1204}%
\special{pa 2980 1204}%
\special{fp}%
%
\special{pn 8}%
\special{pa 2244 1392}%
\special{pa 2428 1392}%
\special{fp}%
\special{sh 1}%
\special{pa 2428 1392}%
\special{pa 2362 1372}%
\special{pa 2376 1392}%
\special{pa 2362 1412}%
\special{pa 2428 1392}%
\special{fp}%
%
\special{pn 8}%
\special{pa 3212 1392}%
\special{pa 3412 1406}%
\special{fp}%
\special{sh 1}%
\special{pa 3412 1406}%
\special{pa 3346 1382}%
\special{pa 3358 1402}%
\special{pa 3344 1422}%
\special{pa 3412 1406}%
\special{fp}%
%
\special{pn 8}%
\special{pa 4044 1164}%
\special{pa 4026 1130}%
\special{pa 4008 1094}%
\special{pa 3990 1058}%
\special{pa 3974 1024}%
\special{pa 3962 992}%
\special{pa 3952 960}%
\special{pa 3946 928}%
\special{pa 3946 900}%
\special{pa 3950 872}%
\special{pa 3958 846}%
\special{pa 3974 822}%
\special{pa 3996 802}%
\special{pa 4024 784}%
\special{pa 4056 770}%
\special{pa 4090 758}%
\special{pa 4128 752}%
\special{pa 4164 750}%
\special{pa 4200 754}%
\special{pa 4234 762}%
\special{pa 4264 778}%
\special{pa 4290 800}%
\special{pa 4312 826}%
\special{pa 4328 858}%
\special{pa 4338 892}%
\special{pa 4342 928}%
\special{pa 4340 962}%
\special{pa 4332 996}%
\special{pa 4318 1026}%
\special{pa 4298 1050}%
\special{pa 4272 1072}%
\special{pa 4244 1092}%
\special{pa 4216 1110}%
\special{pa 4188 1128}%
\special{pa 4164 1148}%
\special{pa 4144 1170}%
\special{pa 4134 1198}%
\special{pa 4128 1230}%
\special{pa 4130 1262}%
\special{pa 4136 1296}%
\special{pa 4144 1326}%
\special{pa 4156 1356}%
\special{pa 4168 1386}%
\special{pa 4180 1414}%
\special{pa 4196 1444}%
\special{pa 4210 1476}%
\special{pa 4224 1510}%
\special{pa 4234 1542}%
\special{pa 4240 1572}%
\special{pa 4238 1602}%
\special{pa 4228 1630}%
\special{pa 4208 1654}%
\special{pa 4182 1674}%
\special{pa 4152 1690}%
\special{pa 4118 1704}%
\special{pa 4086 1712}%
\special{pa 4054 1716}%
\special{pa 4022 1716}%
\special{pa 3992 1716}%
\special{pa 3960 1712}%
\special{pa 3926 1708}%
\special{pa 3890 1702}%
\special{pa 3856 1696}%
\special{pa 3824 1686}%
\special{pa 3796 1672}%
\special{pa 3774 1652}%
\special{pa 3760 1626}%
\special{pa 3752 1596}%
\special{pa 3748 1562}%
\special{pa 3750 1528}%
\special{pa 3754 1494}%
\special{pa 3762 1460}%
\special{pa 3772 1430}%
\special{pa 3786 1400}%
\special{pa 3804 1376}%
\special{pa 3826 1354}%
\special{pa 3850 1336}%
\special{pa 3880 1322}%
\special{pa 3912 1312}%
\special{pa 3948 1306}%
\special{pa 3984 1300}%
\special{pa 4016 1292}%
\special{pa 4040 1276}%
\special{pa 4052 1250}%
\special{pa 4054 1218}%
\special{pa 4054 1184}%
\special{sp}%
%
\special{pn 8}%
\special{pa 3766 1182}%
\special{pa 3754 1160}%
\special{fp}%
%
\special{pn 8}%
\special{pa 4060 1202}%
\special{pa 4050 1150}%
\special{fp}%
%
\special{pn 8}%
\special{pa 1716 1174}%
\special{pa 1712 1140}%
\special{pa 1710 1108}%
\special{pa 1708 1074}%
\special{pa 1706 1042}%
\special{pa 1708 1010}%
\special{pa 1710 978}%
\special{pa 1716 948}%
\special{pa 1726 918}%
\special{pa 1740 890}%
\special{pa 1756 864}%
\special{pa 1778 840}%
\special{pa 1802 818}%
\special{pa 1830 798}%
\special{pa 1862 782}%
\special{pa 1898 770}%
\special{pa 1932 764}%
\special{pa 1968 762}%
\special{pa 2000 766}%
\special{pa 2028 776}%
\special{pa 2050 792}%
\special{pa 2068 816}%
\special{pa 2080 846}%
\special{pa 2088 880}%
\special{pa 2094 914}%
\special{pa 2100 952}%
\special{pa 2102 988}%
\special{pa 2102 1022}%
\special{pa 2098 1054}%
\special{pa 2090 1084}%
\special{pa 2076 1108}%
\special{pa 2056 1126}%
\special{pa 2032 1140}%
\special{pa 2002 1152}%
\special{pa 1970 1162}%
\special{pa 1936 1170}%
\special{pa 1902 1182}%
\special{pa 1868 1196}%
\special{pa 1836 1210}%
\special{pa 1806 1218}%
\special{pa 1778 1216}%
\special{pa 1752 1200}%
\special{pa 1728 1178}%
\special{pa 1716 1164}%
\special{sp}%
%
\special{pn 8}%
\special{pa 2986 1202}%
\special{pa 3020 1200}%
\special{pa 3054 1198}%
\special{pa 3086 1196}%
\special{pa 3118 1190}%
\special{pa 3150 1184}%
\special{pa 3180 1172}%
\special{pa 3206 1156}%
\special{pa 3232 1136}%
\special{pa 3254 1112}%
\special{pa 3272 1084}%
\special{pa 3288 1054}%
\special{pa 3298 1020}%
\special{pa 3304 986}%
\special{pa 3304 952}%
\special{pa 3300 918}%
\special{pa 3288 884}%
\special{pa 3270 854}%
\special{pa 3248 824}%
\special{pa 3220 800}%
\special{pa 3190 780}%
\special{pa 3158 766}%
\special{pa 3126 758}%
\special{pa 3094 756}%
\special{pa 3066 762}%
\special{pa 3040 776}%
\special{pa 3016 798}%
\special{pa 2998 824}%
\special{pa 2982 856}%
\special{pa 2968 888}%
\special{pa 2958 924}%
\special{pa 2952 956}%
\special{pa 2950 990}%
\special{pa 2948 1020}%
\special{pa 2950 1052}%
\special{pa 2954 1082}%
\special{pa 2960 1114}%
\special{pa 2964 1146}%
\special{pa 2972 1178}%
\special{pa 2982 1210}%
\special{pa 3000 1212}%
\special{pa 2996 1212}%
\special{sp}%
\put(7.4300,-13.3600){\makebox(0,0){(FI)}}%
%
\special{pn 8}%
\special{pa 1330 1488}%
\special{pa 1300 1474}%
\special{pa 1300 1444}%
\special{pa 1320 1416}%
\special{pa 1350 1404}%
\special{pa 1380 1398}%
\special{pa 1414 1398}%
\special{pa 1442 1412}%
\special{pa 1442 1444}%
\special{pa 1436 1474}%
\special{pa 1408 1490}%
\special{pa 1376 1496}%
\special{pa 1346 1492}%
\special{pa 1338 1488}%
\special{sp}%
%
\special{pn 8}%
\special{pa 2762 1504}%
\special{pa 2732 1490}%
\special{pa 2730 1458}%
\special{pa 2750 1432}%
\special{pa 2780 1420}%
\special{pa 2810 1412}%
\special{pa 2844 1412}%
\special{pa 2874 1426}%
\special{pa 2876 1456}%
\special{pa 2870 1488}%
\special{pa 2842 1506}%
\special{pa 2812 1512}%
\special{pa 2780 1508}%
\special{pa 2770 1504}%
\special{sp}%
%
\special{pn 8}%
\special{pa 3950 1558}%
\special{pa 3920 1544}%
\special{pa 3918 1514}%
\special{pa 3938 1486}%
\special{pa 3968 1472}%
\special{pa 3998 1466}%
\special{pa 4032 1466}%
\special{pa 4062 1478}%
\special{pa 4064 1510}%
\special{pa 4058 1540}%
\special{pa 4032 1558}%
\special{pa 4000 1564}%
\special{pa 3968 1560}%
\special{pa 3958 1558}%
\special{sp}%
%
\special{pn 8}%
\special{pa 1330 1488}%
\special{pa 1368 1496}%
\special{fp}%
%
\special{pn 8}%
\special{pa 2762 1504}%
\special{pa 2800 1512}%
\special{fp}%
%
\special{pn 8}%
\special{pa 3934 1558}%
\special{pa 3980 1564}%
\special{fp}%
%
\special{pn 8}%
\special{pa 2710 874}%
\special{pa 2948 1170}%
\special{fp}%
\special{sh 1}%
\special{pa 2948 1170}%
\special{pa 2922 1106}%
\special{pa 2914 1128}%
\special{pa 2890 1130}%
\special{pa 2948 1170}%
\special{fp}%
\put(25.4800,-8.3600){\makebox(0,0){$z(t)$}}%
%
\special{pn 8}%
\special{ar 1910 970 68 68  0.0000000 6.2831853}%
%
\special{pn 8}%
\special{ar 3130 980 68 68  0.0000000 6.2831853}%
%
\special{pn 8}%
\special{ar 4138 922 68 68  0.0000000 6.2831853}%
%
\special{pn 4}%
\special{pa 4154 1342}%
\special{pa 4030 1464}%
\special{fp}%
\special{pa 4160 1376}%
\special{pa 4064 1472}%
\special{fp}%
\special{pa 4174 1404}%
\special{pa 4064 1512}%
\special{fp}%
\special{pa 4186 1430}%
\special{pa 3908 1708}%
\special{fp}%
\special{pa 4010 1566}%
\special{pa 3874 1702}%
\special{fp}%
\special{pa 3970 1566}%
\special{pa 3846 1688}%
\special{fp}%
\special{pa 3936 1560}%
\special{pa 3814 1682}%
\special{fp}%
\special{pa 3914 1540}%
\special{pa 3792 1662}%
\special{fp}%
\special{pa 4132 1280}%
\special{pa 3766 1648}%
\special{fp}%
\special{pa 4126 1246}%
\special{pa 3758 1614}%
\special{fp}%
\special{pa 4132 1200}%
\special{pa 3752 1580}%
\special{fp}%
\special{pa 3990 1300}%
\special{pa 3752 1540}%
\special{fp}%
\special{pa 3942 1308}%
\special{pa 3752 1498}%
\special{fp}%
\special{pa 3894 1314}%
\special{pa 3766 1444}%
\special{fp}%
\special{pa 4344 948}%
\special{pa 4050 1240}%
\special{fp}%
\special{pa 4344 906}%
\special{pa 4058 1192}%
\special{fp}%
\special{pa 4330 880}%
\special{pa 4050 1158}%
\special{fp}%
\special{pa 4322 846}%
\special{pa 4030 1138}%
\special{fp}%
\special{pa 4140 988}%
\special{pa 4018 1110}%
\special{fp}%
\special{pa 4106 982}%
\special{pa 4004 1084}%
\special{fp}%
\special{pa 4086 960}%
\special{pa 3990 1056}%
\special{fp}%
\special{pa 4072 934}%
\special{pa 3976 1028}%
\special{fp}%
\special{pa 4086 880}%
\special{pa 3970 996}%
\special{fp}%
\special{pa 4174 750}%
\special{pa 3956 968}%
\special{fp}%
\special{pa 4132 750}%
\special{pa 3950 934}%
\special{fp}%
\special{pa 4078 764}%
\special{pa 3950 892}%
\special{fp}%
\special{pa 4208 756}%
\special{pa 4098 866}%
\special{fp}%
\special{pa 4242 764}%
\special{pa 4154 852}%
\special{fp}%
\special{pa 4268 778}%
\special{pa 4180 866}%
\special{fp}%
\special{pa 4290 798}%
\special{pa 4200 886}%
\special{fp}%
%
\special{pn 4}%
\special{pa 4302 824}%
\special{pa 4208 920}%
\special{fp}%
\special{pa 4330 1002}%
\special{pa 4242 1090}%
\special{fp}%
\special{pa 4140 1314}%
\special{pa 3990 1464}%
\special{fp}%
\special{pa 4200 1458}%
\special{pa 3950 1708}%
\special{fp}%
\special{pa 4214 1484}%
\special{pa 3982 1716}%
\special{fp}%
\special{pa 4222 1518}%
\special{pa 4024 1716}%
\special{fp}%
\special{pa 4234 1546}%
\special{pa 4064 1716}%
\special{fp}%
\special{pa 4242 1580}%
\special{pa 4126 1696}%
\special{fp}%
%
\special{pn 4}%
\special{pa 1570 1192}%
\special{pa 1366 1396}%
\special{fp}%
\special{pa 1528 1192}%
\special{pa 1194 1526}%
\special{fp}%
\special{pa 1298 1464}%
\special{pa 1216 1546}%
\special{fp}%
\special{pa 1318 1484}%
\special{pa 1242 1560}%
\special{fp}%
\special{pa 1344 1498}%
\special{pa 1262 1580}%
\special{fp}%
\special{pa 1386 1498}%
\special{pa 1290 1594}%
\special{fp}%
\special{pa 1678 1246}%
\special{pa 1324 1600}%
\special{fp}%
\special{pa 1670 1294}%
\special{pa 1358 1608}%
\special{fp}%
\special{pa 1658 1348}%
\special{pa 1398 1608}%
\special{fp}%
\special{pa 1638 1410}%
\special{pa 1454 1594}%
\special{fp}%
\special{pa 1616 1472}%
\special{pa 1522 1566}%
\special{fp}%
\special{pa 1664 1220}%
\special{pa 1446 1436}%
\special{fp}%
\special{pa 1644 1200}%
\special{pa 1440 1404}%
\special{fp}%
\special{pa 1610 1192}%
\special{pa 1406 1396}%
\special{fp}%
\special{pa 1488 1192}%
\special{pa 1182 1498}%
\special{fp}%
\special{pa 1454 1186}%
\special{pa 1174 1464}%
\special{fp}%
\special{pa 1412 1186}%
\special{pa 1168 1430}%
\special{fp}%
\special{pa 1372 1186}%
\special{pa 1168 1390}%
\special{fp}%
\special{pa 1330 1186}%
\special{pa 1182 1336}%
\special{fp}%
%
\special{pn 4}%
\special{pa 2950 1200}%
\special{pa 2610 1540}%
\special{fp}%
\special{pa 2726 1464}%
\special{pa 2610 1580}%
\special{fp}%
\special{pa 2732 1498}%
\special{pa 2622 1608}%
\special{fp}%
\special{pa 2766 1504}%
\special{pa 2644 1628}%
\special{fp}%
\special{pa 2800 1512}%
\special{pa 2670 1640}%
\special{fp}%
\special{pa 2848 1504}%
\special{pa 2712 1640}%
\special{fp}%
\special{pa 2990 1404}%
\special{pa 2752 1640}%
\special{fp}%
\special{pa 2990 1444}%
\special{pa 2800 1634}%
\special{fp}%
\special{pa 2984 1492}%
\special{pa 2854 1620}%
\special{fp}%
\special{pa 2990 1362}%
\special{pa 2874 1478}%
\special{fp}%
\special{pa 2990 1322}%
\special{pa 2874 1436}%
\special{fp}%
\special{pa 2984 1288}%
\special{pa 2854 1416}%
\special{fp}%
\special{pa 2984 1246}%
\special{pa 2820 1410}%
\special{fp}%
\special{pa 2976 1212}%
\special{pa 2766 1424}%
\special{fp}%
\special{pa 2908 1200}%
\special{pa 2616 1492}%
\special{fp}%
\special{pa 2874 1192}%
\special{pa 2630 1436}%
\special{fp}%
\special{pa 2834 1192}%
\special{pa 2658 1368}%
\special{fp}%
\special{pa 2780 1206}%
\special{pa 2690 1294}%
\special{fp}%
\end{picture}%

\vspace{2mm} 

\unitlength 0.1in
\begin{picture}( 36.9800,  7.7500)(  1.0500,-11.2500)
%
\special{pn 8}%
\special{pa 1186 800}%
\special{pa 1154 802}%
\special{pa 1122 800}%
\special{pa 1088 798}%
\special{pa 1056 796}%
\special{pa 1024 800}%
\special{pa 996 814}%
\special{pa 974 838}%
\special{pa 956 870}%
\special{pa 948 906}%
\special{pa 948 942}%
\special{pa 960 972}%
\special{pa 982 996}%
\special{pa 1010 1012}%
\special{pa 1044 1018}%
\special{pa 1078 1018}%
\special{pa 1112 1010}%
\special{pa 1140 994}%
\special{pa 1164 974}%
\special{pa 1182 950}%
\special{pa 1194 922}%
\special{pa 1204 890}%
\special{pa 1214 856}%
\special{pa 1218 822}%
\special{pa 1194 804}%
\special{pa 1182 800}%
\special{sp}%
%
\special{pn 8}%
\special{pa 2280 810}%
\special{pa 2246 806}%
\special{pa 2214 806}%
\special{pa 2186 814}%
\special{pa 2162 830}%
\special{pa 2144 856}%
\special{pa 2128 888}%
\special{pa 2114 922}%
\special{pa 2104 960}%
\special{pa 2100 996}%
\special{pa 2106 1022}%
\special{pa 2126 1040}%
\special{pa 2158 1044}%
\special{pa 2194 1040}%
\special{pa 2230 1030}%
\special{pa 2260 1016}%
\special{pa 2284 996}%
\special{pa 2296 972}%
\special{pa 2300 940}%
\special{pa 2300 906}%
\special{pa 2298 872}%
\special{pa 2298 840}%
\special{pa 2298 814}%
\special{sp}%
%
\special{pn 8}%
\special{pa 2280 810}%
\special{pa 2296 810}%
\special{fp}%
%
\special{pn 8}%
\special{pa 3478 748}%
\special{pa 3460 714}%
\special{pa 3444 678}%
\special{pa 3432 646}%
\special{pa 3426 614}%
\special{pa 3430 588}%
\special{pa 3446 564}%
\special{pa 3474 546}%
\special{pa 3508 534}%
\special{pa 3544 532}%
\special{pa 3578 538}%
\special{pa 3606 558}%
\special{pa 3626 586}%
\special{pa 3634 620}%
\special{pa 3630 656}%
\special{pa 3616 684}%
\special{pa 3590 706}%
\special{pa 3562 724}%
\special{pa 3538 744}%
\special{pa 3524 770}%
\special{pa 3524 804}%
\special{pa 3532 836}%
\special{pa 3542 864}%
\special{pa 3556 894}%
\special{pa 3570 926}%
\special{pa 3580 958}%
\special{pa 3578 988}%
\special{pa 3558 1012}%
\special{pa 3528 1028}%
\special{pa 3496 1038}%
\special{pa 3464 1040}%
\special{pa 3432 1036}%
\special{pa 3398 1032}%
\special{pa 3364 1024}%
\special{pa 3338 1008}%
\special{pa 3326 980}%
\special{pa 3324 946}%
\special{pa 3328 912}%
\special{pa 3338 882}%
\special{pa 3356 856}%
\special{pa 3380 838}%
\special{pa 3412 826}%
\special{pa 3448 820}%
\special{pa 3476 808}%
\special{pa 3484 778}%
\special{pa 3484 760}%
\special{sp}%
%
\special{pn 8}%
\special{pa 2870 768}%
\special{pa 2864 758}%
\special{fp}%
%
\special{pn 8}%
\special{pa 3486 770}%
\special{pa 3480 740}%
\special{fp}%
%
\special{pn 8}%
\special{pa 1236 792}%
\special{pa 1232 758}%
\special{pa 1230 724}%
\special{pa 1232 692}%
\special{pa 1238 662}%
\special{pa 1252 634}%
\special{pa 1272 610}%
\special{pa 1300 590}%
\special{pa 1334 578}%
\special{pa 1368 572}%
\special{pa 1398 580}%
\special{pa 1420 598}%
\special{pa 1430 628}%
\special{pa 1436 664}%
\special{pa 1438 700}%
\special{pa 1436 734}%
\special{pa 1424 758}%
\special{pa 1400 774}%
\special{pa 1368 784}%
\special{pa 1334 794}%
\special{pa 1302 808}%
\special{pa 1272 814}%
\special{pa 1248 798}%
\special{pa 1236 786}%
\special{sp}%
%
\special{pn 8}%
\special{pa 2296 810}%
\special{pa 2330 808}%
\special{pa 2364 804}%
\special{pa 2394 796}%
\special{pa 2420 780}%
\special{pa 2442 756}%
\special{pa 2458 726}%
\special{pa 2464 692}%
\special{pa 2460 658}%
\special{pa 2446 626}%
\special{pa 2420 598}%
\special{pa 2388 580}%
\special{pa 2358 574}%
\special{pa 2328 584}%
\special{pa 2306 606}%
\special{pa 2290 636}%
\special{pa 2282 670}%
\special{pa 2278 704}%
\special{pa 2280 734}%
\special{pa 2286 764}%
\special{pa 2290 798}%
\special{pa 2308 816}%
\special{pa 2302 816}%
\special{sp}%
\put(3.3000,-8.1000){\makebox(0,0){(FII)}}%
%
\special{pn 8}%
\special{pa 1134 586}%
\special{pa 1104 598}%
\special{pa 1076 612}%
\special{pa 1048 628}%
\special{pa 1020 644}%
\special{pa 992 662}%
\special{pa 966 682}%
\special{pa 942 704}%
\special{pa 920 728}%
\special{pa 902 754}%
\special{pa 886 782}%
\special{pa 872 812}%
\special{pa 864 846}%
\special{pa 860 880}%
\special{pa 862 914}%
\special{pa 866 948}%
\special{pa 876 980}%
\special{pa 890 1008}%
\special{pa 906 1032}%
\special{pa 928 1052}%
\special{pa 952 1064}%
\special{pa 980 1074}%
\special{pa 1010 1076}%
\special{pa 1042 1076}%
\special{pa 1074 1072}%
\special{pa 1110 1064}%
\special{pa 1144 1052}%
\special{pa 1180 1038}%
\special{pa 1216 1022}%
\special{pa 1250 1002}%
\special{pa 1284 982}%
\special{pa 1318 960}%
\special{pa 1348 936}%
\special{pa 1378 912}%
\special{pa 1404 888}%
\special{pa 1428 862}%
\special{pa 1450 834}%
\special{pa 1470 806}%
\special{pa 1486 778}%
\special{pa 1500 748}%
\special{pa 1510 718}%
\special{pa 1518 686}%
\special{pa 1520 656}%
\special{pa 1520 624}%
\special{pa 1516 590}%
\special{pa 1508 558}%
\special{pa 1494 526}%
\special{pa 1478 496}%
\special{pa 1456 474}%
\special{pa 1428 462}%
\special{pa 1396 460}%
\special{pa 1364 466}%
\special{pa 1330 478}%
\special{pa 1302 494}%
\special{pa 1274 512}%
\special{pa 1250 532}%
\special{pa 1222 548}%
\special{pa 1194 562}%
\special{pa 1164 574}%
\special{pa 1134 586}%
\special{pa 1134 586}%
\special{sp}%
%
\special{pn 8}%
\special{pa 2180 590}%
\special{pa 2154 610}%
\special{pa 2128 628}%
\special{pa 2104 648}%
\special{pa 2080 668}%
\special{pa 2056 692}%
\special{pa 2034 716}%
\special{pa 2014 742}%
\special{pa 1998 770}%
\special{pa 1984 798}%
\special{pa 1974 830}%
\special{pa 1966 862}%
\special{pa 1964 896}%
\special{pa 1968 930}%
\special{pa 1974 964}%
\special{pa 1986 996}%
\special{pa 2000 1026}%
\special{pa 2020 1052}%
\special{pa 2040 1072}%
\special{pa 2064 1086}%
\special{pa 2092 1094}%
\special{pa 2120 1098}%
\special{pa 2150 1096}%
\special{pa 2180 1088}%
\special{pa 2214 1078}%
\special{pa 2246 1064}%
\special{pa 2278 1046}%
\special{pa 2310 1024}%
\special{pa 2342 1002}%
\special{pa 2374 976}%
\special{pa 2404 950}%
\special{pa 2432 922}%
\special{pa 2458 894}%
\special{pa 2482 866}%
\special{pa 2502 836}%
\special{pa 2522 806}%
\special{pa 2538 774}%
\special{pa 2552 744}%
\special{pa 2564 712}%
\special{pa 2572 680}%
\special{pa 2576 650}%
\special{pa 2576 618}%
\special{pa 2574 586}%
\special{pa 2568 554}%
\special{pa 2556 524}%
\special{pa 2542 492}%
\special{pa 2524 464}%
\special{pa 2500 438}%
\special{pa 2474 422}%
\special{pa 2446 416}%
\special{pa 2414 420}%
\special{pa 2382 432}%
\special{pa 2354 450}%
\special{pa 2328 472}%
\special{pa 2306 496}%
\special{pa 2284 518}%
\special{pa 2260 540}%
\special{pa 2234 558}%
\special{pa 2206 574}%
\special{pa 2180 590}%
\special{sp}%
%
\special{pn 8}%
\special{pa 3272 686}%
\special{pa 3250 708}%
\special{pa 3228 730}%
\special{pa 3208 754}%
\special{pa 3188 780}%
\special{pa 3170 808}%
\special{pa 3152 836}%
\special{pa 3138 868}%
\special{pa 3126 902}%
\special{pa 3118 936}%
\special{pa 3114 968}%
\special{pa 3116 1000}%
\special{pa 3124 1028}%
\special{pa 3140 1054}%
\special{pa 3162 1074}%
\special{pa 3188 1092}%
\special{pa 3220 1104}%
\special{pa 3254 1114}%
\special{pa 3288 1122}%
\special{pa 3324 1124}%
\special{pa 3360 1126}%
\special{pa 3392 1124}%
\special{pa 3426 1118}%
\special{pa 3458 1110}%
\special{pa 3488 1100}%
\special{pa 3518 1086}%
\special{pa 3548 1070}%
\special{pa 3574 1052}%
\special{pa 3600 1032}%
\special{pa 3626 1010}%
\special{pa 3650 986}%
\special{pa 3672 958}%
\special{pa 3692 930}%
\special{pa 3712 900}%
\special{pa 3730 868}%
\special{pa 3746 834}%
\special{pa 3760 798}%
\special{pa 3774 762}%
\special{pa 3784 724}%
\special{pa 3792 686}%
\special{pa 3800 648}%
\special{pa 3802 610}%
\special{pa 3804 574}%
\special{pa 3802 540}%
\special{pa 3796 508}%
\special{pa 3786 478}%
\special{pa 3772 450}%
\special{pa 3756 428}%
\special{pa 3734 410}%
\special{pa 3708 396}%
\special{pa 3676 386}%
\special{pa 3644 382}%
\special{pa 3610 382}%
\special{pa 3576 388}%
\special{pa 3544 398}%
\special{pa 3514 410}%
\special{pa 3486 428}%
\special{pa 3462 448}%
\special{pa 3438 472}%
\special{pa 3418 498}%
\special{pa 3402 524}%
\special{pa 3386 552}%
\special{pa 3368 578}%
\special{pa 3350 604}%
\special{pa 3328 628}%
\special{pa 3306 652}%
\special{pa 3284 674}%
\special{pa 3272 686}%
\special{sp}%
%
\special{pn 4}%
\special{pa 1260 526}%
\special{pa 862 922}%
\special{fp}%
\special{pa 958 866}%
\special{pa 870 954}%
\special{fp}%
\special{pa 946 916}%
\special{pa 876 986}%
\special{fp}%
\special{pa 946 954}%
\special{pa 894 1006}%
\special{fp}%
\special{pa 958 980}%
\special{pa 908 1032}%
\special{fp}%
\special{pa 978 1000}%
\special{pa 926 1050}%
\special{fp}%
\special{pa 1004 1012}%
\special{pa 952 1064}%
\special{fp}%
\special{pa 1036 1018}%
\special{pa 984 1070}%
\special{fp}%
\special{pa 1074 1018}%
\special{pa 1016 1076}%
\special{fp}%
\special{pa 1126 1006}%
\special{pa 1054 1076}%
\special{fp}%
\special{pa 1388 782}%
\special{pa 1106 1064}%
\special{fp}%
\special{pa 1516 692}%
\special{pa 1170 1038}%
\special{fp}%
\special{pa 1496 750}%
\special{pa 1240 1006}%
\special{fp}%
\special{pa 1522 648}%
\special{pa 1438 730}%
\special{fp}%
\special{pa 1522 610}%
\special{pa 1438 692}%
\special{fp}%
\special{pa 1510 584}%
\special{pa 1432 660}%
\special{fp}%
\special{pa 1502 552}%
\special{pa 1432 622}%
\special{fp}%
\special{pa 1490 526}%
\special{pa 1420 596}%
\special{fp}%
\special{pa 1478 500}%
\special{pa 1400 578}%
\special{fp}%
\special{pa 1458 482}%
\special{pa 1368 570}%
\special{fp}%
\special{pa 1438 462}%
\special{pa 1318 584}%
\special{fp}%
\special{pa 1400 462}%
\special{pa 1068 794}%
\special{fp}%
\special{pa 1240 660}%
\special{pa 1100 802}%
\special{fp}%
\special{pa 1228 712}%
\special{pa 1138 802}%
\special{fp}%
\special{pa 1234 744}%
\special{pa 1176 802}%
\special{fp}%
\special{pa 1234 782}%
\special{pa 1208 808}%
\special{fp}%
\special{pa 1254 802}%
\special{pa 1214 840}%
\special{fp}%
\special{pa 1278 814}%
\special{pa 1208 884}%
\special{fp}%
\special{pa 1336 794}%
\special{pa 1182 948}%
\special{fp}%
\special{pa 1356 468}%
\special{pa 1022 802}%
\special{fp}%
%
\special{pn 4}%
\special{pa 1182 564}%
\special{pa 862 884}%
\special{fp}%
\special{pa 1106 602}%
\special{pa 870 840}%
\special{fp}%
\special{pa 1030 642}%
\special{pa 882 788}%
\special{fp}%
%
\special{pn 4}%
\special{pa 2460 416}%
\special{pa 1968 910}%
\special{fp}%
\special{pa 2294 622}%
\special{pa 1968 948}%
\special{fp}%
\special{pa 2280 672}%
\special{pa 1980 974}%
\special{fp}%
\special{pa 2152 838}%
\special{pa 1986 1006}%
\special{fp}%
\special{pa 2114 916}%
\special{pa 2000 1030}%
\special{fp}%
\special{pa 2102 966}%
\special{pa 2018 1050}%
\special{fp}%
\special{pa 2102 1006}%
\special{pa 2038 1070}%
\special{fp}%
\special{pa 2114 1030}%
\special{pa 2064 1082}%
\special{fp}%
\special{pa 2140 1044}%
\special{pa 2088 1094}%
\special{fp}%
\special{pa 2178 1044}%
\special{pa 2128 1094}%
\special{fp}%
\special{pa 2230 1030}%
\special{pa 2166 1094}%
\special{fp}%
\special{pa 2544 756}%
\special{pa 2230 1070}%
\special{fp}%
\special{pa 2568 692}%
\special{pa 2300 960}%
\special{fp}%
\special{pa 2576 646}%
\special{pa 2300 922}%
\special{fp}%
\special{pa 2576 608}%
\special{pa 2460 724}%
\special{fp}%
\special{pa 2576 570}%
\special{pa 2466 678}%
\special{fp}%
\special{pa 2562 544}%
\special{pa 2460 646}%
\special{fp}%
\special{pa 2550 518}%
\special{pa 2448 622}%
\special{fp}%
\special{pa 2536 494}%
\special{pa 2428 602}%
\special{fp}%
\special{pa 2524 468}%
\special{pa 2402 590}%
\special{fp}%
\special{pa 2504 448}%
\special{pa 2376 576}%
\special{fp}%
\special{pa 2486 430}%
\special{pa 2332 582}%
\special{fp}%
\special{pa 2384 800}%
\special{pa 2300 884}%
\special{fp}%
\special{pa 2338 806}%
\special{pa 2294 852}%
\special{fp}%
\special{pa 2492 846}%
\special{pa 2326 1012}%
\special{fp}%
\special{pa 2280 710}%
\special{pa 2166 826}%
\special{fp}%
\special{pa 2280 750}%
\special{pa 2224 806}%
\special{fp}%
\special{pa 2288 782}%
\special{pa 2262 806}%
\special{fp}%
\special{pa 2422 416}%
\special{pa 1968 870}%
\special{fp}%
\special{pa 2242 558}%
\special{pa 1974 826}%
\special{fp}%
%
\special{pn 4}%
\special{pa 2128 634}%
\special{pa 2006 756}%
\special{fp}%
%
\special{pn 4}%
\special{pa 3328 752}%
\special{pa 3116 962}%
\special{fp}%
\special{pa 3444 674}%
\special{pa 3116 1002}%
\special{fp}%
\special{pa 3456 700}%
\special{pa 3124 1034}%
\special{fp}%
\special{pa 3462 732}%
\special{pa 3142 1052}%
\special{fp}%
\special{pa 3328 906}%
\special{pa 3162 1072}%
\special{fp}%
\special{pa 3322 950}%
\special{pa 3188 1084}%
\special{fp}%
\special{pa 3328 982}%
\special{pa 3212 1098}%
\special{fp}%
\special{pa 3340 1008}%
\special{pa 3238 1110}%
\special{fp}%
\special{pa 3366 1020}%
\special{pa 3270 1116}%
\special{fp}%
\special{pa 3392 1034}%
\special{pa 3302 1122}%
\special{fp}%
\special{pa 3424 1040}%
\special{pa 3340 1122}%
\special{fp}%
\special{pa 3462 1040}%
\special{pa 3380 1122}%
\special{fp}%
\special{pa 3500 1040}%
\special{pa 3424 1116}%
\special{fp}%
\special{pa 3738 842}%
\special{pa 3476 1104}%
\special{fp}%
\special{pa 3764 778}%
\special{pa 3584 956}%
\special{fp}%
\special{pa 3782 720}%
\special{pa 3572 930}%
\special{fp}%
\special{pa 3796 668}%
\special{pa 3558 906}%
\special{fp}%
\special{pa 3802 624}%
\special{pa 3552 874}%
\special{fp}%
\special{pa 3802 586}%
\special{pa 3540 848}%
\special{fp}%
\special{pa 3802 546}%
\special{pa 3526 822}%
\special{fp}%
\special{pa 3796 514}%
\special{pa 3520 790}%
\special{fp}%
\special{pa 3788 482}%
\special{pa 3636 636}%
\special{fp}%
\special{pa 3776 458}%
\special{pa 3636 598}%
\special{fp}%
\special{pa 3756 438}%
\special{pa 3622 572}%
\special{fp}%
\special{pa 3738 418}%
\special{pa 3604 554}%
\special{fp}%
\special{pa 3718 400}%
\special{pa 3578 540}%
\special{fp}%
\special{pa 3692 386}%
\special{pa 3546 534}%
\special{fp}%
\special{pa 3660 380}%
\special{pa 3508 534}%
\special{fp}%
\special{pa 3622 380}%
\special{pa 3142 860}%
\special{fp}%
\special{pa 3424 618}%
\special{pa 3124 918}%
\special{fp}%
%
\special{pn 4}%
\special{pa 3430 650}%
\special{pa 3328 752}%
\special{fp}%
\special{pa 3578 386}%
\special{pa 3296 668}%
\special{fp}%
\special{pa 3514 412}%
\special{pa 3412 514}%
\special{fp}%
\special{pa 3700 918}%
\special{pa 3552 1066}%
\special{fp}%
\special{pa 3482 752}%
\special{pa 3404 828}%
\special{fp}%
\special{pa 3482 790}%
\special{pa 3456 816}%
\special{fp}%
%
\special{pn 8}%
\special{pa 2014 508}%
\special{pa 2246 764}%
\special{fp}%
\special{sh 1}%
\special{pa 2246 764}%
\special{pa 2216 700}%
\special{pa 2210 724}%
\special{pa 2186 728}%
\special{pa 2246 764}%
\special{fp}%
\put(18.9400,-4.3500){\makebox(0,0){$z(t)$}}%
%
\special{pn 8}%
\special{pa 1644 820}%
\special{pa 1864 820}%
\special{fp}%
\special{sh 1}%
\special{pa 1864 820}%
\special{pa 1798 800}%
\special{pa 1812 820}%
\special{pa 1798 840}%
\special{pa 1864 820}%
\special{fp}%
%
\special{pn 8}%
\special{pa 2820 820}%
\special{pa 3020 820}%
\special{fp}%
\special{sh 1}%
\special{pa 3020 820}%
\special{pa 2954 800}%
\special{pa 2968 820}%
\special{pa 2954 840}%
\special{pa 3020 820}%
\special{fp}%
\end{picture}%

\vspace{2mm} 

\unitlength 0.1in
\begin{picture}( 39.4800, 10.4000)(  2.3000,-18.0500)
%
\special{pn 8}%
\special{pa 1420 956}%
\special{pa 1382 958}%
\special{pa 1344 962}%
\special{pa 1310 968}%
\special{pa 1280 978}%
\special{pa 1256 996}%
\special{pa 1238 1020}%
\special{pa 1230 1054}%
\special{pa 1230 1090}%
\special{pa 1238 1126}%
\special{pa 1252 1158}%
\special{pa 1274 1178}%
\special{pa 1300 1182}%
\special{pa 1330 1176}%
\special{pa 1362 1162}%
\special{pa 1396 1144}%
\special{pa 1430 1128}%
\special{pa 1460 1118}%
\special{pa 1488 1118}%
\special{pa 1512 1132}%
\special{pa 1532 1158}%
\special{pa 1544 1190}%
\special{pa 1550 1228}%
\special{pa 1548 1266}%
\special{pa 1536 1302}%
\special{pa 1516 1334}%
\special{pa 1488 1358}%
\special{pa 1458 1374}%
\special{pa 1428 1380}%
\special{pa 1398 1376}%
\special{pa 1370 1360}%
\special{pa 1346 1336}%
\special{pa 1328 1308}%
\special{pa 1314 1276}%
\special{pa 1302 1246}%
\special{pa 1288 1218}%
\special{pa 1266 1196}%
\special{pa 1236 1182}%
\special{pa 1202 1176}%
\special{pa 1170 1182}%
\special{pa 1140 1196}%
\special{pa 1118 1218}%
\special{pa 1102 1246}%
\special{pa 1092 1278}%
\special{pa 1090 1312}%
\special{pa 1094 1348}%
\special{pa 1102 1380}%
\special{pa 1116 1412}%
\special{pa 1132 1438}%
\special{pa 1154 1464}%
\special{pa 1178 1486}%
\special{pa 1204 1504}%
\special{pa 1234 1522}%
\special{pa 1264 1536}%
\special{pa 1296 1548}%
\special{pa 1328 1558}%
\special{pa 1360 1568}%
\special{pa 1394 1572}%
\special{pa 1426 1576}%
\special{pa 1458 1578}%
\special{pa 1492 1576}%
\special{pa 1524 1572}%
\special{pa 1554 1566}%
\special{pa 1584 1556}%
\special{pa 1614 1544}%
\special{pa 1642 1528}%
\special{pa 1668 1510}%
\special{pa 1694 1490}%
\special{pa 1718 1468}%
\special{pa 1740 1442}%
\special{pa 1760 1416}%
\special{pa 1778 1388}%
\special{pa 1794 1358}%
\special{pa 1808 1326}%
\special{pa 1818 1294}%
\special{pa 1826 1260}%
\special{pa 1832 1226}%
\special{pa 1834 1192}%
\special{pa 1834 1158}%
\special{pa 1828 1126}%
\special{pa 1816 1098}%
\special{pa 1800 1072}%
\special{pa 1780 1050}%
\special{pa 1754 1030}%
\special{pa 1726 1014}%
\special{pa 1696 1000}%
\special{pa 1664 986}%
\special{pa 1632 976}%
\special{pa 1602 966}%
\special{pa 1570 958}%
\special{pa 1540 952}%
\special{pa 1508 950}%
\special{pa 1476 950}%
\special{pa 1444 954}%
\special{pa 1430 956}%
\special{sp}%
%
\special{pn 8}%
\special{pa 3820 1144}%
\special{pa 3788 1150}%
\special{pa 3756 1158}%
\special{pa 3726 1168}%
\special{pa 3698 1182}%
\special{pa 3672 1204}%
\special{pa 3656 1230}%
\special{pa 3654 1260}%
\special{pa 3666 1290}%
\special{pa 3686 1320}%
\special{pa 3714 1348}%
\special{pa 3744 1368}%
\special{pa 3776 1382}%
\special{pa 3808 1390}%
\special{pa 3842 1390}%
\special{pa 3874 1384}%
\special{pa 3904 1370}%
\special{pa 3930 1350}%
\special{pa 3952 1324}%
\special{pa 3966 1292}%
\special{pa 3972 1260}%
\special{pa 3968 1228}%
\special{pa 3954 1196}%
\special{pa 3932 1170}%
\special{pa 3906 1152}%
\special{pa 3876 1140}%
\special{pa 3844 1138}%
\special{pa 3814 1138}%
\special{sp}%
%
\special{pn 8}%
\special{pa 1308 874}%
\special{pa 1286 884}%
\special{fp}%
%
\special{pn 8}%
\special{pa 2430 1198}%
\special{pa 2448 1232}%
\special{pa 2466 1266}%
\special{pa 2484 1298}%
\special{pa 2504 1328}%
\special{pa 2526 1354}%
\special{pa 2548 1376}%
\special{pa 2572 1392}%
\special{pa 2600 1404}%
\special{pa 2630 1410}%
\special{pa 2662 1408}%
\special{pa 2696 1398}%
\special{pa 2728 1382}%
\special{pa 2760 1362}%
\special{pa 2784 1338}%
\special{pa 2802 1310}%
\special{pa 2808 1280}%
\special{pa 2806 1248}%
\special{pa 2794 1218}%
\special{pa 2776 1190}%
\special{pa 2752 1168}%
\special{pa 2728 1148}%
\special{pa 2698 1130}%
\special{pa 2668 1120}%
\special{pa 2638 1116}%
\special{pa 2606 1118}%
\special{pa 2576 1130}%
\special{pa 2546 1144}%
\special{pa 2516 1160}%
\special{pa 2490 1178}%
\special{pa 2464 1196}%
\special{pa 2440 1216}%
\special{pa 2416 1240}%
\special{pa 2396 1264}%
\special{pa 2378 1292}%
\special{pa 2364 1322}%
\special{pa 2354 1354}%
\special{pa 2350 1386}%
\special{pa 2352 1416}%
\special{pa 2360 1446}%
\special{pa 2374 1476}%
\special{pa 2394 1502}%
\special{pa 2418 1528}%
\special{pa 2442 1550}%
\special{pa 2470 1570}%
\special{pa 2498 1586}%
\special{pa 2528 1600}%
\special{pa 2558 1612}%
\special{pa 2590 1620}%
\special{pa 2622 1626}%
\special{pa 2654 1628}%
\special{pa 2686 1628}%
\special{pa 2716 1624}%
\special{pa 2748 1616}%
\special{pa 2778 1606}%
\special{pa 2810 1596}%
\special{pa 2840 1582}%
\special{pa 2868 1566}%
\special{pa 2898 1550}%
\special{pa 2924 1532}%
\special{pa 2950 1512}%
\special{pa 2974 1490}%
\special{pa 2998 1466}%
\special{pa 3018 1442}%
\special{pa 3036 1414}%
\special{pa 3052 1386}%
\special{pa 3064 1354}%
\special{pa 3074 1322}%
\special{pa 3082 1288}%
\special{pa 3084 1254}%
\special{pa 3086 1218}%
\special{pa 3082 1184}%
\special{pa 3076 1150}%
\special{pa 3066 1120}%
\special{pa 3054 1090}%
\special{pa 3036 1064}%
\special{pa 3016 1040}%
\special{pa 2992 1020}%
\special{pa 2966 1004}%
\special{pa 2936 992}%
\special{pa 2906 980}%
\special{pa 2874 972}%
\special{pa 2840 964}%
\special{pa 2806 958}%
\special{pa 2774 952}%
\special{pa 2742 946}%
\special{pa 2710 942}%
\special{pa 2678 940}%
\special{pa 2646 940}%
\special{pa 2614 940}%
\special{pa 2582 944}%
\special{pa 2550 950}%
\special{pa 2518 958}%
\special{pa 2486 970}%
\special{pa 2456 986}%
\special{pa 2430 1004}%
\special{pa 2410 1026}%
\special{pa 2400 1052}%
\special{pa 2398 1082}%
\special{pa 2402 1114}%
\special{pa 2412 1148}%
\special{pa 2424 1182}%
\special{pa 2432 1204}%
\special{sp}%
%
\special{pn 8}%
\special{pa 3460 1352}%
\special{pa 3462 1386}%
\special{pa 3466 1420}%
\special{pa 3470 1452}%
\special{pa 3478 1484}%
\special{pa 3488 1512}%
\special{pa 3502 1538}%
\special{pa 3520 1562}%
\special{pa 3542 1582}%
\special{pa 3566 1600}%
\special{pa 3594 1614}%
\special{pa 3626 1626}%
\special{pa 3658 1634}%
\special{pa 3694 1640}%
\special{pa 3730 1644}%
\special{pa 3768 1644}%
\special{pa 3806 1640}%
\special{pa 3844 1636}%
\special{pa 3880 1628}%
\special{pa 3918 1616}%
\special{pa 3952 1604}%
\special{pa 3986 1588}%
\special{pa 4018 1568}%
\special{pa 4048 1548}%
\special{pa 4074 1524}%
\special{pa 4098 1498}%
\special{pa 4118 1470}%
\special{pa 4136 1440}%
\special{pa 4150 1410}%
\special{pa 4162 1378}%
\special{pa 4170 1346}%
\special{pa 4176 1312}%
\special{pa 4178 1278}%
\special{pa 4178 1246}%
\special{pa 4174 1212}%
\special{pa 4168 1180}%
\special{pa 4158 1148}%
\special{pa 4146 1118}%
\special{pa 4132 1088}%
\special{pa 4114 1062}%
\special{pa 4094 1036}%
\special{pa 4072 1012}%
\special{pa 4048 990}%
\special{pa 4022 968}%
\special{pa 3996 950}%
\special{pa 3968 934}%
\special{pa 3938 922}%
\special{pa 3908 910}%
\special{pa 3878 902}%
\special{pa 3846 896}%
\special{pa 3812 892}%
\special{pa 3780 892}%
\special{pa 3746 896}%
\special{pa 3714 904}%
\special{pa 3682 914}%
\special{pa 3654 930}%
\special{pa 3630 952}%
\special{pa 3610 980}%
\special{pa 3594 1012}%
\special{pa 3588 1044}%
\special{pa 3588 1078}%
\special{pa 3598 1110}%
\special{pa 3616 1136}%
\special{pa 3636 1162}%
\special{pa 3640 1192}%
\special{pa 3630 1226}%
\special{pa 3608 1252}%
\special{pa 3578 1264}%
\special{pa 3546 1264}%
\special{pa 3514 1268}%
\special{pa 3486 1290}%
\special{pa 3466 1338}%
\special{pa 3460 1368}%
\special{pa 3464 1362}%
\special{sp}%
%
\special{pn 8}%
\special{pa 1700 1090}%
\special{pa 1666 1086}%
\special{pa 1640 1096}%
\special{pa 1634 1128}%
\special{pa 1646 1166}%
\special{pa 1674 1190}%
\special{pa 1710 1194}%
\special{pa 1740 1178}%
\special{pa 1748 1148}%
\special{pa 1736 1112}%
\special{pa 1710 1094}%
\special{pa 1694 1090}%
\special{sp}%
%
\special{pn 8}%
\special{pa 2960 1160}%
\special{pa 2928 1136}%
\special{pa 2902 1124}%
\special{pa 2884 1136}%
\special{pa 2878 1172}%
\special{pa 2884 1212}%
\special{pa 2902 1244}%
\special{pa 2926 1256}%
\special{pa 2952 1240}%
\special{pa 2970 1208}%
\special{pa 2968 1176}%
\special{pa 2954 1154}%
\special{sp}%
%
\special{pn 8}%
\special{pa 4050 1100}%
\special{pa 4024 1068}%
\special{pa 3998 1048}%
\special{pa 3974 1054}%
\special{pa 3956 1086}%
\special{pa 3956 1118}%
\special{pa 3982 1134}%
\special{pa 4020 1140}%
\special{pa 4050 1132}%
\special{pa 4050 1102}%
\special{pa 4048 1092}%
\special{sp}%
%
\special{pn 8}%
\special{pa 1890 1204}%
\special{pa 2048 1204}%
\special{fp}%
\special{sh 1}%
\special{pa 2048 1204}%
\special{pa 1982 1184}%
\special{pa 1996 1204}%
\special{pa 1982 1224}%
\special{pa 2048 1204}%
\special{fp}%
%
\special{pn 8}%
\special{pa 3102 1210}%
\special{pa 3278 1210}%
\special{fp}%
\special{sh 1}%
\special{pa 3278 1210}%
\special{pa 3212 1190}%
\special{pa 3226 1210}%
\special{pa 3212 1230}%
\special{pa 3278 1210}%
\special{fp}%
\put(5.0000,-12.9000){\makebox(0,0){(FIII)}}%
%
\special{pn 4}%
\special{pa 1490 952}%
\special{pa 1230 1112}%
\special{fp}%
\special{pa 1556 950}%
\special{pa 1248 1140}%
\special{fp}%
\special{pa 1602 962}%
\special{pa 1266 1170}%
\special{fp}%
\special{pa 1640 978}%
\special{pa 1312 1180}%
\special{fp}%
\special{pa 1686 990}%
\special{pa 1484 1114}%
\special{fp}%
\special{pa 1714 1012}%
\special{pa 1522 1132}%
\special{fp}%
\special{pa 1636 1102}%
\special{pa 1540 1160}%
\special{fp}%
\special{pa 1634 1142}%
\special{pa 1548 1196}%
\special{fp}%
\special{pa 1662 1164}%
\special{pa 1546 1236}%
\special{fp}%
\special{pa 1698 1182}%
\special{pa 1544 1276}%
\special{fp}%
\special{pa 1832 1140}%
\special{pa 1524 1330}%
\special{fp}%
\special{pa 1830 1180}%
\special{pa 1262 1530}%
\special{fp}%
\special{pa 1828 1222}%
\special{pa 1300 1548}%
\special{fp}%
\special{pa 1828 1262}%
\special{pa 1336 1564}%
\special{fp}%
\special{pa 1816 1308}%
\special{pa 1392 1570}%
\special{fp}%
\special{pa 1786 1368}%
\special{pa 1448 1574}%
\special{fp}%
\special{pa 1746 1432}%
\special{pa 1524 1568}%
\special{fp}%
\special{pa 1446 1378}%
\special{pa 1226 1514}%
\special{fp}%
\special{pa 1390 1372}%
\special{pa 1188 1496}%
\special{fp}%
\special{pa 1354 1354}%
\special{pa 1170 1468}%
\special{fp}%
\special{pa 1336 1326}%
\special{pa 1142 1444}%
\special{fp}%
\special{pa 1318 1298}%
\special{pa 1116 1422}%
\special{fp}%
\special{pa 1300 1268}%
\special{pa 1108 1386}%
\special{fp}%
\special{pa 1292 1234}%
\special{pa 1098 1352}%
\special{fp}%
\special{pa 1274 1204}%
\special{pa 1090 1318}%
\special{fp}%
\special{pa 1246 1182}%
\special{pa 1092 1276}%
\special{fp}%
\special{pa 1180 1182}%
\special{pa 1112 1224}%
\special{fp}%
\special{pa 1814 1110}%
\special{pa 1736 1158}%
\special{fp}%
\special{pa 1806 1076}%
\special{pa 1738 1118}%
\special{fp}%
\special{pa 1778 1052}%
\special{pa 1720 1088}%
\special{fp}%
%
\special{pn 4}%
\special{pa 1750 1030}%
\special{pa 1694 1066}%
\special{fp}%
\special{pa 1414 958}%
\special{pa 1232 1072}%
\special{fp}%
\special{pa 1348 960}%
\special{pa 1232 1030}%
\special{fp}%
%
\special{pn 4}%
\special{pa 2896 976}%
\special{pa 2428 1178}%
\special{fp}%
\special{pa 2842 962}%
\special{pa 2416 1146}%
\special{fp}%
\special{pa 2776 952}%
\special{pa 2402 1114}%
\special{fp}%
\special{pa 2710 944}%
\special{pa 2398 1078}%
\special{fp}%
\special{pa 2634 940}%
\special{pa 2406 1038}%
\special{fp}%
\special{pa 2940 992}%
\special{pa 2660 1114}%
\special{fp}%
\special{pa 2986 1010}%
\special{pa 2704 1132}%
\special{fp}%
\special{pa 3010 1038}%
\special{pa 2738 1154}%
\special{fp}%
\special{pa 2898 1122}%
\special{pa 2774 1176}%
\special{fp}%
\special{pa 2880 1166}%
\special{pa 2788 1206}%
\special{fp}%
\special{pa 2872 1206}%
\special{pa 2800 1238}%
\special{fp}%
\special{pa 2896 1234}%
\special{pa 2814 1270}%
\special{fp}%
\special{pa 3076 1192}%
\special{pa 2796 1314}%
\special{fp}%
\special{pa 3078 1228}%
\special{pa 2746 1372}%
\special{fp}%
\special{pa 3082 1264}%
\special{pa 2438 1542}%
\special{fp}%
\special{pa 3076 1304}%
\special{pa 2472 1564}%
\special{fp}%
\special{pa 3068 1344}%
\special{pa 2506 1586}%
\special{fp}%
\special{pa 3050 1388}%
\special{pa 2540 1608}%
\special{fp}%
\special{pa 3022 1438}%
\special{pa 2596 1622}%
\special{fp}%
\special{pa 2952 1506}%
\special{pa 2672 1626}%
\special{fp}%
\special{pa 2664 1408}%
\special{pa 2414 1514}%
\special{fp}%
\special{pa 2588 1404}%
\special{pa 2390 1488}%
\special{fp}%
\special{pa 2554 1382}%
\special{pa 2368 1462}%
\special{fp}%
\special{pa 2530 1354}%
\special{pa 2364 1426}%
\special{fp}%
\special{pa 2506 1328}%
\special{pa 2350 1396}%
\special{fp}%
\special{pa 2482 1300}%
\special{pa 2356 1354}%
\special{fp}%
\special{pa 2468 1270}%
\special{pa 2366 1314}%
\special{fp}%
\special{pa 2456 1240}%
\special{pa 2392 1266}%
\special{fp}%
\special{pa 3072 1158}%
\special{pa 2970 1202}%
\special{fp}%
\special{pa 3060 1126}%
\special{pa 2978 1162}%
\special{fp}%
%
\special{pn 4}%
\special{pa 3056 1090}%
\special{pa 2952 1136}%
\special{fp}%
\special{pa 3032 1064}%
\special{pa 2930 1108}%
\special{fp}%
%
\special{pn 4}%
\special{pa 3972 938}%
\special{pa 3612 1130}%
\special{fp}%
\special{pa 4008 956}%
\special{pa 3638 1154}%
\special{fp}%
\special{pa 4034 980}%
\special{pa 3646 1190}%
\special{fp}%
\special{pa 3662 1220}%
\special{pa 3612 1246}%
\special{fp}%
\special{pa 3648 1266}%
\special{pa 3468 1362}%
\special{fp}%
\special{pa 3664 1296}%
\special{pa 3464 1402}%
\special{fp}%
\special{pa 3690 1320}%
\special{pa 3470 1438}%
\special{fp}%
\special{pa 3716 1344}%
\special{pa 3476 1474}%
\special{fp}%
\special{pa 3742 1370}%
\special{pa 3492 1504}%
\special{fp}%
\special{pa 3788 1384}%
\special{pa 3508 1534}%
\special{fp}%
\special{pa 3844 1392}%
\special{pa 3524 1562}%
\special{fp}%
\special{pa 4178 1250}%
\special{pa 3550 1588}%
\special{fp}%
\special{pa 4174 1290}%
\special{pa 3586 1606}%
\special{fp}%
\special{pa 4170 1332}%
\special{pa 3622 1626}%
\special{fp}%
\special{pa 4166 1372}%
\special{pa 3678 1634}%
\special{fp}%
\special{pa 4142 1424}%
\special{pa 3734 1642}%
\special{fp}%
\special{pa 4100 1486}%
\special{pa 3810 1640}%
\special{fp}%
\special{pa 4172 1214}%
\special{pa 3944 1338}%
\special{fp}%
\special{pa 4166 1180}%
\special{pa 3966 1286}%
\special{fp}%
\special{pa 4160 1144}%
\special{pa 3970 1246}%
\special{fp}%
\special{pa 4144 1114}%
\special{pa 3964 1210}%
\special{fp}%
\special{pa 4008 1148}%
\special{pa 3948 1180}%
\special{fp}%
\special{pa 3962 1134}%
\special{pa 3922 1156}%
\special{fp}%
\special{pa 3946 1104}%
\special{pa 3876 1142}%
\special{fp}%
\special{pa 4060 1006}%
\special{pa 3820 1134}%
\special{fp}%
\special{pa 4086 1030}%
\special{pa 4026 1062}%
\special{fp}%
\special{pa 4102 1060}%
\special{pa 4042 1092}%
\special{fp}%
\special{pa 4128 1084}%
\special{pa 4038 1132}%
\special{fp}%
\special{pa 3582 1262}%
\special{pa 3472 1322}%
\special{fp}%
\special{pa 3936 918}%
\special{pa 3596 1100}%
\special{fp}%
%
\special{pn 4}%
\special{pa 3890 904}%
\special{pa 3590 1064}%
\special{fp}%
\special{pa 3834 896}%
\special{pa 3594 1024}%
\special{fp}%
\special{pa 3768 892}%
\special{pa 3608 978}%
\special{fp}%
\put(13.6000,-18.6000){\makebox(0,0){$R_n(t'),t'\in B_0'$}}%
\put(26.4000,-18.8000){\makebox(0,0){$R_n(t),t\in l_i$}}%
\put(40.0000,-18.9000){\makebox(0,0){$R_n(t''),t''\in B_0''$}}%
\put(20.7000,-8.5000){\makebox(0,0){$z(t)$}}%
%
\special{pn 8}%
\special{pa 2100 950}%
\special{pa 2400 1220}%
\special{fp}%
\special{sh 1}%
\special{pa 2400 1220}%
\special{pa 2364 1162}%
\special{pa 2360 1184}%
\special{pa 2338 1190}%
\special{pa 2400 1220}%
\special{fp}%
\end{picture}%

\vspace{5mm}

 For $t\in B$ we consider the $L_1$-$(L_0$-)principal function $ p_n(t,z) 
(q_n(t,z))$ and the harmonic span $s_n(t)$ for  $(R_n'(t ), \xi(t), 
\eta(t))$. Then we have 

\begin{lemma}(S.\,Hamano\, \cite{hamano-4})
\label{lem:C1-ness} 
Let ${\mathcal  R} $ be  
a Stein manifold and each $R(t), t\in \Delta $ is planar. Then
\begin{enumerate}
 \item [1.] $p_n(t,z)$ and $q_n(t,z)$ are continuous for $(t,z)$ in 
       ${\mathcal  R}_n'$, and $s_n(t)$ is continuous on $B$; 
 \item [2.]assume that at each singular point $z(t)$ of $\partial R_n(t),\, 
t\in l_i' \subset {\mathcal  L}'$ such that  $z(t)\in R_n'(t)$, 
case ({\bf c1})  only occurs. Then 
\begin{enumerate}
 \item [(i)]$p_n(t,z)$ 
       and $q_n(t,z)$ are of class $C^1$ for $(t,z)$ on 
${\mathcal   R}_n' \setminus \{\xi|_B, \eta|_B\}$;
 \item [(ii)]$s_n(t)$ is        $C^1$ subharmonic on $B$.
 \end{enumerate}
 \item [3.]there exist counter-examples for case ({\bf  c2})
such that  $p_n(t,z)$ or $q_n(t,z)$ is not of class $C^1$ on ${\mathcal  
       R}_n' \setminus \{\xi|_B, \eta|_B\}$, and $s_n(t)$ is neither of 
       class $C^1$ on $B$ nor subharmonic on $B$.
\end{enumerate} 
\end{lemma}

As an example of figure (F\,I), let $B=\{
|t|< 1/10\}$; $D=\{|z|<2\}$;  
$\psi_1=(e^{-100+|t|^2}/|z-1|^2)-1; \psi_2= 
 |z^2-1|- (1-2 \Re\,t -|t|^2); \psi_3= 
(e^{-100+|t|^2}/|z+1|^2)-1$ and ${\mathcal  R}=\{ (t,z) \in B \times D: \psi_1<0, \psi_2<0,
 \psi_3<0\}$, so that $\partial {\mathcal  R}$ is $C^\omega$ strictly pseudoconvex in $B \times D$.
 Then 
the arc $l' =\{t\in B: 2 \Re\, t\ +|t|^2=0\}$ divides $B$ into two domains 
$B'\cup B''$ such that  
 the connected 
 components of $\partial  R(t), t\in l' $ consists of two circles $ \psi_1(t,z)=0$, 
 $\psi_3(t,z)=0$ and the leminiscate $C: |z^2-1|=1$ which is singular at 
 $z(t)=0$.

 As an example of ${\mathcal  L}''$. Let $B, \ D $ be the same as above. Let $\psi(t,z):
= |z-t|^2 +|t|^2 +2  \Re\, t$ and  put ${\mathcal  R}=\{ (t,z) \in B \times D:
 \psi(t,z)<0\}$. Then the arc $l'' =\{t\in B: \phi(t)=0\}$, where $\phi(t)= 
-|t|^2- 2 \Re\, t $,  divides $B$ into two domains 
$B'(B'')=\{t\in B: \phi(t)<0 (>0)\}$ such that $R(t)=\emptyset$ 
for $t\in B'\cup l''$ and 
$ R(t)= \{|z-t|^2 < \phi(t)\}$  for $t\in B''$. 
The mapping ${\mathfrak z}: t\in l'' \to z(t)= t $ so that 
 $(t,t)\in \partial {\mathcal R} $ but $t \not \in \partial R(t)$, 
and each $R(t), t\in B''$ is a disk $\{|z-t|< \phi(t)\}$ 
which schrinkingly approaches the point $z=t^0 $ as $t\to t^0\in l''$.

\vspace{2mm} 
 Since the Stein manifold admits a $C^\omega$ strictly plurisubharmonic 
exhaustion function, we immediately have 
\begin{lemma} \label{last-remark}
  Let ${\mathcal  R}: \Delta \to R(t)$ be of type $({\bf  A})$ or $({\bf B})$, and let $\xi, 
\eta\in \Gamma (\Delta , {\mathcal  R})$ such that  $\xi \cap \eta=\emptyset$. 
Assume
\begin{enumerate} 
 \item [($\star$)] $R(t), t\in \Delta $ 
is homeomorphic to a domain $D$ in $\mathbb{ C}_w$ bounded by  $\nu$ 
boundary component such that  $1\le \nu<\infty$ and $\nu$ is independent of 
$t\in \Delta $. 
\end{enumerate}
Then,  for any $t_0\in \Delta $, there exists 
a small disk $B \Subset \Delta$ of center $t_0$ such that  
we find an increasing sequence $\{{\mathcal  R}'_n\}_n$ of case 
$({\bf  c1})$ 
 such that  $\lim_{ n\to \infty} {\mathcal  R}'_n={\mathcal  R}|_{B}$.
\end{lemma}
Now we consider  the variation 
$
{\mathcal  R}: t\in \Delta \to R(t)
$ of type $({\bf  A})$. 
 Let $\xi, \ \eta \in \Gamma (\Delta, {\mathcal  R})$ 
such that  $\xi \cap \eta=\emptyset$.
We fix so small disk $B\Subset \Delta $ that we can fix local parameters 
$(t,z)$ of $\xi |_B$ and $\eta|_B$ in ${\mathcal  R}|_B$ and $\{ {\mathcal R}_n\}_n$
satisfies conditions in {\it  Preparation } to 
 these $\Delta $ and $B$, precisely saying, similar to (\ref{eqn:Rn}) we define 
\begin{eqnarray}
\label{eqn:Rn-2} \qquad \qquad  {\mathcal  R}_n= 
\mbox{the conn. comp. of ${\mathcal  R}(a_n)|_B 
  $ which contains $g|_B$},  
\end{eqnarray}
which  satisfies  conditions 
$1)$ and $2)$ in {\it Preparation}. 
We put ${\mathcal  R}_n=\cup_{t\in B}(t, R_n(t))$, and for $t\in B$ we denote by $R'_n(t)$ 
the connected component of $R_n(t)$ 
     which contains $g|_B (\supset \xi (t), \eta(t)) $ and put ${\mathcal  R}_n'
=\cup_{t\in B}(t, 
     R_n'(t))$. We then have the $L_1$-$(L_0$-)principal function $ p_n(t,z) 
(q_n(t,z))$; 
 the $L_1$-$(L_0$-)constant $ \alpha _n\, (\beta _n)$
 and the harmonic span $s_n(t)$ for $(R_n'(t), \xi(t),\eta(t))$.
In one complex variable it is known (cf: \S\,8,\,Chap.\,III in \cite{ahlfors-sario}) 
that, for a fixed $t\in B$,  
 $p_n(t,z) (q_n(t,z))$ uniformly converges to a certain function $p(t,z) (q(t,z))$ on any 
compact set in $R(t) \setminus \{\xi (t), \eta(t)\}$. Thus 
$p(t,z ) \, (q(t,z))$ is harmonic  on $R(t) \setminus \{\xi(t), 
\eta(t)\}$
 with the same pole as $p_n(t,z)(q_n(t,z))$ at $\xi(t)$ and $\eta(t)$. Putting
$\alpha (t) (\beta (t))= \lim_{ z\to \eta(t)}(p(t,z)- \log |z-\eta(t)|) \ 
(
\lim_{ z\to \eta(t)} (q(t,z)- \log |z-\eta(t)|))$, we have $\alpha _n(t) 
\to \alpha (t) \ (\beta_n (t) \to  \beta (t))$ as $n \to 
\infty$. We call $p(t,z) (q(t,z))$
 the $L_1$-($L_0$-){\it principal function} and $s(t): =\alpha (t)- \beta (t)$ the 
{\it harmonic span} for $(R(t), \xi (t), \eta(t))$. Since $R(t)$ is planar, we have 
\begin{align}
 \label{sn-s} s_n(t) \searrow s(t) \  \mbox{ as $n\to \infty$,  } \quad \mbox{  $t\in  B$}.
\end{align}
Their proofs in \cite{ahlfors-sario} also imply that, 
given  $K \Subset {\mathcal  R}|_B \setminus  \{\xi|_B,  \eta|_B\}$, for sufficiently large $n$, 
\begin{align}
   \label{eqn:bounded}
\mbox{$p_n(t,z), q_n(t,z), p(t,z), q(t,z)$ 
are uniformly bounded on $K$.}
\end{align}
\noindent 
Note that $p(t,z), q(t,z), \alpha (t), \beta (t)$ depend on 
the choice of local parameters  of  $\xi |_B$ and $\eta|_B$ but 
$s(t)$ does not depend on them, so that $s(t)$ is a non-negative function on $\Delta $. 

\vspace{2mm} Using these notations we have the following approximation 
\begin{theorem}\label{lemma:rigidity-2}
Let ${\mathcal  R}: t\in \Delta \to R(t)$ be of 
  type $({\bf A})$  and let $\xi, \eta \in \Gamma (\Delta , {\mathcal  
 R})$ such that  $\xi \cap \eta =\emptyset$. Let $s(t)$ denote 
the harmonic span for $(R(t, \xi (t), \eta(t)), t\in \Delta 
 $. 
Assume
\begin{enumerate}
 \item [($*$)]  \,for any $t_0\in \Delta$, there exists a small disk 
$B \Subset \Delta$ of center $t_0$ such that 
we find an increasing sequence 
$\{{\mathcal  R}'_n\}_n$ of case $({\bf  c1})$ 
such that  $\lim_{ n\to \infty} {\mathcal  R}'_n={\mathcal  R}|_{B}$.
\end{enumerate}
Then
\begin{enumerate}
 \item [1.]  $s(t)$ is  subharmonic on $\Delta$; 
\item [2.] \ (Simultaneous uniformization) if 
       $s(t)$ is harmonic on $\Delta $, then  ${\mathcal  R}$ 
is a biholomprphic to a univalent 
 domain  in $\Delta   \times \mathbb{ P}$.
\end{enumerate}
\end{theorem}
\noindent {\it Proof.} \ 
To show 1. let $t_0\in \Delta$. Then we have a disk $B 
\subset \Delta $ which satisfies condition ($*$). By
 2.(ii) in Lemma \ref{lem:C1-ness}, 
$s_n(t)$ is $C^1$  subharmonic
 on $B$, hence $s(t)$ is suharmonic on $B$, and on $\Delta $.
To prove 2., 
we cover $\Delta $ by small disks $\{B_i\}_{i=1,2,\ldots }
$ with condition ($*$), i..e, for fixed $B_i$, 
 we find an increasing sequences $\{ {\mathcal R}'_n 
\}_n $ (depending on $B_i$) such that    each ${\mathcal  
R}_n', n=1,2,\ldots $ is of case ({\bf c1} ) and $ 
 \lim_{ n\to \infty} {\mathcal  R}_n'={\mathcal  
 R}|_{B_i}$. We divide the proof of 2. into two steps.

{\it  $1^{st}$ step.\ Each ${\mathcal R }|_{B_i}$, $i=1,2,\ldots $ is 
 biholmorphic to a univalent domain ${\mathcal  D}_i$ 
in $B \times \mathbb{ P}$.}

In fact, we simply write $B=B_i$. We put ${\mathcal  R}_n'=\cup_{t\in B}(t, R_n'(t)), \ n=1,2,\ldots $
and consider $p_n(t,z), q_n(t,z)$ and $s_n(t)$ for each $(R_n'(t), \xi(t), 
\eta(t)),$ $ t\in B$ as above. We put 
\begin{eqnarray} \label{eqn:nrmz}
\begin{array}{lll}
  P_n(t,z)(P(t,z))&= 
e^{p_n(t,z)+i p_n(t,z)^*} \ (e^{p(t,z)+ip(t,z)^*});\\[2mm]
Q_n(t,z)(Q(t,z))&= 
e^{q_n(t,z)+i q_n(t,z)^*}\  (e^{q(t,z)+iq(t,z)^*}),
\end{array}
\end{eqnarray}
 which are  normalized 
\begin{eqnarray}
 \label{eqn:normalization}
\textstyle{ \frac{ 1}{z-\xi(t)}}+a_0(t)+ a_1(t)(z-\xi(t)) + \ldots \quad \mbox{  
at 
 $z=\xi(t)$}.
\end{eqnarray}
For a fixed $t\in B$, $P_n(t,z)(Q_n(t,z))$   uniformly converges
to $P(t,z)(Q(t,z))$ on any compact set in $R(t)$; 
 $w=P_n(t,z) (Q_n(t,z))$ is a circular (radial) slit mapping on $R_n'(t)$, 
and hence $P(t,z)(Q(t,z))$ is  an univalent function on $R(t)$. We call 
such $P(t,z)(Q(t,z))$ the circular (radial) slit mapping for $(R(t), \xi 
(t), \eta(t))$. 
  For the 1st  step it  suffices to show 
\begin{enumerate}
 \item [$(a)$]  {the harmonicity of $s(t)$ on 
  $B$ implies that  $P(t,z)$ is holomorphic for two complex variables 
$(t,z)$ in ${\mathcal  R}|_B \setminus \{\xi|_B\}$.}
\end{enumerate}

In fact, fix a point $(t_0, z_0)$ in ${\mathcal  R}|_B \setminus 
\{\xi|_B , \eta|_B\}$ and let $B_0 \times V \Subset {\mathcal  R}|_B \setminus
 \{\xi|_B , \eta|_B\}$ 
be a bi-disk 
centered at $(t_0,z_0)$, a local coordinate of a neighborhood of 
$(t_0,z_0)$. We put $f(t,z):
= \frac{\partial p(t,z)}
{\partial z} $ for $(t,z)\in B_0 \times V$. From  (\ref{eqn:normalization})
  it suffices for $(a)$ to 
prove that $f(t,z)$ is holomorphic for $(t,z)$ in $B_0 \times V$.
Since each $f(t,z), t\in B_0$ is 
holomorphic for $z\in V$ and $f(t,z)$ is uniformly bounded in $B_0 
\times V$ by (\ref{eqn:bounded}), it thus suffices for (a) to show that, for any fixed $z' \in V$, 
it holds $\frac{\partial f(t,z') }{\partial \overline{ t }}=0 $ 
on $B_0$   
in the sense of 
distribution, i.e., it holds, for  any $\varphi (t)=\varphi (t_1+it_2)\in C_0^{\infty}(B_0)$, 
\begin{eqnarray}
 \label{eqn:last-one}
I:= \int_{ B_0} f(t,z') \frac{\partial \varphi (t)}{\partial \overline{ t}}
dt_1dt_2=0.
\end{eqnarray}

To prove this  by contradiction, assume $I\ne 0$. 
 We fix a small disk $V_0=\{|z-z'|<r_0\} \Subset V$ of center  $z'$, 
so that we have 
$R'_n(t) \Supset V_0$ for any $t\in B_0$ and  $n \ge \exists \, n_0$.
We  see from the  mean-value theorem for holomorphic functions for $z$ that 
$$
I= \frac{ 1}{\pi r_0^2} \iint_{ B_0 \times V_0} f(t,z) \frac{\partial \varphi (t)}
{\partial \overline{ 
t}}dt_1dt_2dxdy.
$$
Since $f_n(t,z) : =\frac{\partial p_n(t,z)}{\partial z}  \to f(t,z) \ (n 
\to \infty)$ uniformly on $V_0$ for a fixed $t\in B_0$ and since 
$f_n(t,z), f(t,z)$ are uniformly bounded in $B_0 \times V_0 $ by (\ref{eqn:bounded}), 
it follows
 from Lebesgue bounded theorem that 
\begin{align*}
& \quad \ \ \ I= \frac{ 1}{\pi r_0^2} \ \lim_{ n\to \infty} 
\iint_{ B_0 \times V_0} f_n(t,z) \frac{\partial \varphi (t)}
{\partial \overline{ 
t}}dt_1dt_2dxdy.
\\[1mm]
& \therefore \quad 
\bigl|\frac{ 1}{\pi r_0^2}\iint_{ B_0 \times V_0} f_n(t,z) \frac{\partial \varphi (t)}{\partial 
\overline{ t} }dt_1dt_2dxdy\bigr|\ge \frac{ |I|}{2} >0
 \ \ \mbox{ for $ \mbox{any} \  n  \ge \exists\,N$}.
\end{align*}

On the other hand, using 2.(ii) in Lemma \ref{lem:C1-ness} 
under condition $(*)$ in Theorem \ref{lemma:rigidity-2}, we see that, for a 
fixed $z\in V_0$, $p_n(t,z)$, and hence $f_n(t,z)$ is of class $C^1$ for 
$t\in B_0$. 
It follows that 
$$
\int_{ B_0}\!\! f_n(t,z)\frac{\partial \varphi (t)}{\partial \overline{ t}} 
dt_1dt_2 
=
- \int_{ B_0}\!\! \varphi (t)\frac{\partial f_n(t,z)}{\partial \overline{ t}} 
dt_1dt_2.
$$
Hence, putting $I_0=\frac{\pi r_0^2 |I|}{2} >0$, we have from Schwarz inequality 
\begin{align*}
 I_0^2&\le \bigl(\iint_{ B_0 \times V_0}|\varphi (t)|^2dt_1dt_2dxdy\bigr)\
\bigl(\iint_{ B_0 \times V_0}\bigl|
 \frac{\partial f_n(t,z)}{\partial \overline{ t}}\bigr|^2dt_1dt_2dxdy\bigr)\\
&=: C \iint_{ B_0 \times   V_0} 
\bigl|\frac{\partial f_n(t,z)}{\partial \overline{ t}}\bigr|^2 dt_1dt_2dxdy,
\end{align*}
where $C>0$ is independent of $n$. 
Lemma \ref{lem:vari-formula-h-span} 
and  ni) for ${\mathcal L}'$ in {\it  Preparation} 
 for the pseudoconvex domain $
 {\mathcal  R}_n$ imply 
 $$
0\le \frac{ 4}{\pi}\ \int_{ R'_n(t)} \bigl|
\frac{\partial f_n(t,z)}{\partial \overline{ t}}\bigr|^2 dxdy
\le \frac{\partial^2 s_n(t)}{\partial t \partial  \overline{ t} } \quad 
\mbox{ for any $t\in B \setminus {\mathcal  L}'$}.
$$
Since ${\mathcal  L}'$ (depending on $n$) consists of a finite number of 
$C^\omega$ arcs in $B$; $R_n'(t) \supset V_0$ for $n 
\ge n_0$, and $f_n \in C^1(B_0 \times V_0)$, it follows that 
$$
I_0^2 \le 
C \iint_{ (B_0\setminus {\mathcal  L}') \times   V_0} 
\bigl|\frac{\partial f_n(t,z)}{\partial \overline{ t}}\bigr|^2 dt_1dt_2dxdy
\le \frac{ C\pi}4 \int_{ B_0 \setminus {\mathcal  L}'}
\frac{\partial^2 s_n(t)}{\partial t \partial \overline{ t}}  dt_1dt_2.
$$
We fix 
a disk $B_1: B_0 \Subset B_1 \Subset B $ 
and a $C_0^\infty$ function 
$\varphi _1(t)\ge 0$ on $B_1$ such that
 $\varphi _1(t)\equiv 1$ on $B_0$. Since 
$\frac{\partial^2 s_n(t)}{\partial t \partial \overline{ t}} \ge 0$ on 
 $B_1 \setminus {\mathcal  L}' $, we have 
$$
\int_{ B_0 \setminus {\mathcal  L}'} \frac{\partial^2 s_n(t)}{\partial t \partial 
 \overline{ t}}dt_1dt_2 \le
 \int_{ B_1 \setminus  {\mathcal  L}'}  \varphi _1(t)
\frac{\partial^2 s_n(t)}{\partial t \partial 
 \overline{ t}}dt_1dt_2.
$$
Since $s_n(t)$ is of class $C^1$ on $B$ and $\varphi (t)\equiv 0$ on 
$\partial B_1$, we have 
$$
 \int_{ B_1 \setminus  {\mathcal  L}'}  \varphi _1(t)
\frac{\partial^2 s_n(t)}{\partial t \partial 
 \overline{ t}}dt_1dt_2 = \int_{ B_1}  s_n(t)
\frac{\partial^2 \varphi _1(t)}{\partial t \partial 
 \overline{ t}}dt_1dt_2,
$$
 both being equal to $
-\frac{1}{4}\int_{ B_1}(\frac{\partial \varphi _1}{\partial t_1 }
\frac{\partial s_n}{\partial t_1 }+
\frac{\partial \varphi _1}{\partial t_2} \frac{\partial s_n}{\partial t_2} )
dt_1dt_2$. We have by (\ref{sn-s})
\begin{align*}
 0<I_0^2 &\le \frac{ C \pi}4\int_{B_1} s_n(t) \frac{\partial^2 \varphi _1(t)}{\partial t \partial 
\overline{ t}} dt_1 dt_2\\
& \to \frac{ C\pi}4 \int_{ B_1} 
s(t) \frac{\partial^2 \varphi _1(t)}{\partial t \partial 
\overline{ t}} dt_1 dt_2  \quad \mbox{ as $n\to \infty$}\\
&=0 \quad \mbox{ by the harmonicity of $s(t)$ on $B$,}
\end{align*}
which is a contradiction,  and the 1st step is proved.

\vspace{2mm} {\it $2^{nd}$ step. \ Assertion 2. is true.}

In fact, fix $B_i, i=1,2,\ldots $ and let  $P_i(t,z)$ denote  
the circular slit mapping for $(R(t), \xi (t), 
\eta(t))$ used in (a) in the 1st step for ${\mathcal  R}|_{B_i}$. 
From the theory of one complex variable, for a fixed $t\in B_i\cap B_j$, 
 there exists $a_{ij}(t)\ne 0$ such that $P_i(t,z)= a_{ij}(t) P_j(t,z)$ on $R(t)$.
 Since $a_{ij}(t)$ is holomorphic on $B_i\cap B_j$ and since $\Delta $ is an open Riemann surface, 
we have nonvanishing holomorphic function $a_i(t)$ on $B_i$ such that  
$a_{ij}(t)= a_j(t)/a_i(t)$ on $B_i \cap B_j$. Thus, $ a_i(t)P_i(t,z)$ on $B_i, i=1,2,\ldots $ 
defines a holomorphic function ${\mathcal  P}(t,z)$ 
on ${\mathcal  R}$, so that $T: (t,z)\in 
{\mathcal  R}\to (t,w)=(t,{\mathcal  P}(t,z))\in B \times \mathbb{ P}_w $ proves 
the 2nd step. 
\hfill $\Box$ 
\begin{corollary} (Rigidity) \label{cor:rigidity-2}
Let ${\mathcal  R}: t\in \Delta \to R(t)$ be a variation of type $({\bf A})$ or 
$({\bf  B})$ and let $\xi , \eta \in \Gamma (\Delta , {\mathcal  R})$ 
 such that  $\xi \ne \eta$. Let $s(t)$ denote the harmonic span for each 
 $(R(t), \xi(t),\eta(t)), t\in \Delta $. Assume that $R(t), t\in \Delta $ 
satisfies condition {\rm ($\star$)} 
in Lemma \ref{last-remark}.
Then
 \begin{enumerate}
\item [1.] in case when ${\mathcal  R}$ is of type {\rm ({\bf A})}, $s(t)$ is subharmonic on $\Delta$. Moreover,  if 
	    $s(t)$ is harmonic on $\Delta $, then 
	    ${\mathcal  R}$ is equivalent to a trivial variation;
\item [2.] in case when ${\mathcal  R}$ is of type {\rm ({\bf  B})}, $s(t)$ is constant on $\Delta$. Moreover, 
	    if there exists at least one ideal boundary component  $C(t)$ of 
       $R(t), t\in \Delta $ such that
\begin{enumerate} 
\item [(i)] $C(t)$ moves homotopically  with $t\in \Delta$ in ${\mathcal R}  $;
\item [(ii)] each $C(t), t\in  \Delta$ is positive harmonic measure on $R(t)$,
\end{enumerate}
then ${\mathcal  R}$ is equivalent to a trivial variation. 
\end{enumerate}
\end{corollary}
\noindent 
{\it  Proof.} \ The 
 proofs of assertions 1. and 2. are essentially same, we give 
 the proof of 2.   By Lemma \ref{last-remark} and 1. in  Theorem \ref{lemma:rigidity-2},
$s(t)$ is subharmonic on the compact  $\Delta $, so that $s(t)$ is constant on $\Delta $.
By Lemma \ref{last-remark}
 we cover $\Delta $ by small disks $\{B_i\}_{i=1,2,\ldots }
$ which satisfies condition $(*)$ in Theorem \ref{lemma:rigidity-2}. 
  Since $s(t)\equiv {\rm const.} $ on $B_i$, it follows by
the proof of 2. in   
Theorem \ref{lemma:rigidity-2} that the circular slit mapping $P_i(t,z)$ 
for $(R(t), \xi(t), \eta(t))$ for $t\in B_i$
is holomorphic for $t\in B_i $. Since 
$D_i(t):= P_i(t,R(t))$ is a 
circular slit domain with $\nu$ circular arcs $\{A_j^{(1)}(t), A_j^{(2)}(t)\}$
(which might reduce to 
a point,\,i.e., $A_j^{(1)}(t)= A_j^{(2)}(t))$ for some $j$),
it follows from Kanten 
Satz that $A_j^{(k)}(t), k=1,2; j=1,\ldots, \nu$ is holomorphic on $B_i$.
For each $i=1,2,\ldots $  we conventionally rename  
arc $\{A_1^{(1)}(t), A_1^{(2)}(t)\}= P_i(t, C(t))$, where $C(t)$ is stated 
in 2. By the homotpy condition (i), $A_1(t)$ is single-valued on $\Delta$. By  (ii),  the arc $\{A_1^{(1)}, A_1^{(2)}\}$
 does not reduce to a point. By the same argument in 3.\,(ii) in Theorem \ref{cor:rigidity}, we see  that 
that $\widetilde P_i(t,z):= P_i(t,z)/A_1^{(1)}(t)$ is independent of $t\in B_i$, hence so is 
of $i=1,2,\ldots$. We denote it by $\widetilde {\mathcal  P}(t,z)=\widetilde {\mathcal P}(z)$ on ${\mathcal  R}$. Then 
 $T_0:(t,z)\in {\mathcal  R}\to (t,w)=(t, \widetilde {\mathcal  P}(z))$ 
biholomorphically maps ${\mathcal  R}$ onto 
$\Delta \times \widetilde D_1$ where   
$\widetilde D_1$ is the same form 
$(\diamond)$ in (ii) in Theorem \ref{cor:rigidity} (but some boundary 
components of $\widetilde D_1$ might be points $\widetilde A_j$).
 \hfill $\Box$ 

\vspace{2mm}Applying Corollary \ref{cor:rigidity-2} to the case when each $R(t), 
t\in \Delta$ is simply connected, we have 
\begin{corollary} \label{cor:last} 
 Corollary \ref{cor:mdp} holds under the weaker 
condition  
that ${\mathcal R}=\cup_{t\in B}(t, R(t))$ is a Stein manifold. 
\end{corollary}

 \

\end{document}